# A Low Order Finite Element Method for Poroelasticity with Applications to Lung Modelling

Lorenz Berger

Keble College

University of Oxford

A thesis submitted for the degree of

*Doctor of Philosophy*

Trinity Term 2015

## Abstract

In the last few decades modelling deformation and flow in porous media has been of great interest due to its application in various fields including biomechanics, soil mechanics, geophysics, physical chemistry and material sciences. Particularly in biology, virtually any application of poroelasticity implies the use of nonlinear constitutive models, irregular three-dimensional geometries, complicated boundary conditions and jumps in material coefficients, characteristics that can only be simulated numerically.

In this thesis we develop a stabilised finite element method for solving the equations of poroelasticity to enable solving complex models of biological tissues such as the human lungs. For the proposed numerical scheme, we use the lowest possible approximation order: piecewise constant approximation for the pressure, and piecewise linear continuous elements for the displacements and fluid flux. Due to the discontinuous pressure approximation, sharp pressure gradients due to changes in material coefficients or boundary layer solutions can be captured reliably. We begin by developing theoretical results for approximating the linear poroelastic equations valid in small deformations. In particular, we prove existence and uniqueness, an energy estimate and an optimal a-priori error estimate for the discretised problem. We then extend this work and construct a stabilised finite element method to solve the poroelastic equations valid in large deformations. We present the linearisation and discretisation for this nonlinear problem, and give a detailed account of the implementation. We rigorously test both the linear and nonlinear finite element method using numerous test problems to verify theoretical stability and convergence results, and the method's ability to reliably capture steep pressure gradients.

Finally, we derive a poroelastic model for lung parenchyma coupled to an airway fluid network model, and develop a stable method to solve the coupled model. Numerical simulations, on a realistic lung geometry, illustrate the coupling between the poroelastic medium and the network flow model, and simulations of tidal breathing are shown to reproduce global physiologically realistic measurements. We also investigate the effect of airway constriction and tissue weakening on the ventilation, tissue stress and alveolar pressure distribution.

# Acknowledgements

My biggest thanks goes out to my supervisors. Dr. David Kay is a living legend. His enthusiasm, energy and kindness have made my time in Oxford very enjoyable. Not only has he given up countless hours to further my mathematical understanding but also acted as an excellent research and football mentor. I am also indebted to Dr. Rafel Bordas who has guided me through the DPhil and has been a constant source of ideas and support. I would also like to express my appreciation for Professor Simon Tavener who gave me a great amount of his time and attention and invited me for a brilliant stay at Colorado State University. I would like to thank Dr. Kelly Burrowes for many inspiring conversations about lung modelling, and Professor Vicente Grau for introducing me to this project and a trip to Paris. I am also thankful to my Transfer and Confirmation examiners Dr. Jonathan Whitley and Professor Kevin Burrage for providing detailed feedback and suggestions that have shaped much of this thesis.

For all the company I would like to thank all my friends at the Computational Biology group, the DTC, and Keble College. Also thanks to the Redemption crew for organising many unforgettable trips and putting things into context. Finally, I'd like to thank my family in Swansea and Bavaria for their love and support throughout.

# Publications

Below are a list of publications which directly relate to the work described in this thesis.

- **L. Berger**, R. Bordas, D. Kay, and S. Tavener; Stabilized low-order finite element approximation for linear three-field poroelasticity *SIAM Journal on Scientific Computing (Accepted)*
- **L. Berger**, R. Bordas, D. Kay, and S. Tavener; A stabilized finite element method for finite-strain three-field poroelasticity *Computational Mechanics (Under Review)*
- **L. Berger**, R. Bordas, K. Burrowes, V. Grau, D. Kay, and S. Tavener; A poroelastic model coupled to a fluid network with applications in lung modelling *International Journal for Numerical Methods in Biomedical Engineering (Accepted)*

# Conference Presentations

The work described in this thesis was presented at the following international conferences:

- **L. Berger**, R. Bordas, K. Burrowes, C. Brightling, R. Hartley, D. Kay; Understanding The Interdependence Between Parenchymal Deformation And Ventilation In Obstructive Lung Disease, *The American Thoracic Society conference, San Diego, May 2014. (Poster)*
- **L. Berger**, R. Bordas, D. Kay; Solving the Generalised Large Deformation Poroelastic Equations for Modelling Tissue Deformation and Ventilation in the Lung, *European Numerical Mathematics and Advanced Applications conference, EPFL, Lausanne, August 2013. (Oral)*

# Contents

# Chapter 1

# Introduction

Poroelasticity is a theory in which a complex fluid-structure interaction is approximated by a superposition of the solid and fluid components. This theory can capture complex interactions between a deformable porous medium and the fluid flow within it, and has originally been developed to study numerous geomechanical applications ranging from reservoir engineering (Phillips and Wheeler, 2007a) to earthquake fault zones (White and Borja, 2008). Poroelastic models have since been used to model a variety of biological tissues and processes. Simulations using these models can help to advance the understanding of the biomechanics of the tissue under investigation. However, after many decades of research there remain numerous challenges associated with the numerical solution of these poroelastic models.

We begin this chapter with a brief overview of poroelastic models in biology. We then highlight some of the numerical challenges that will form the main motivation for the work presented in this thesis. Finally, we outline the goals and structure of the thesis.

## 1.1 Poroelastic models in biology

Poroelastic models have been proposed for a variety of biological tissues and processes. Unlike many geomechanics applications, which usually assume small deformations in the deformable porous medium, these biological poroelastic models often experience large deformations and require the more complicated nonlinear poroelastic theory.

For example, the coupling of flow in coronary vessels with the mechanical deformation of myocardial tissue is a central feature of cardiac physiology and can be accounted for using a poroelastic model of coronary perfusion (Hyde, 2013). This coupling has been shown to exist in the large epicardial coronary vessels within which flow is impeded and even reversed during contraction. This complicated interplay between the dynamics of vessel compression with resistance and pressure gradients has motivated the development of poroelastic models (Cookson et al., 2012).

Another example is modelling tissue deformation and the ventilation in the lungs. To achieve this tight coupling between the tissue deformation and the ventilation we will develop a multiscale model in Chapter 7 that approximates the lung parenchyma by a biphasic (tissue and air, ignoring blood) poroelastic model, that is then coupled to an airway fluid network model. Such an integrated model of ventilation and tissue mechanics is particularly important for understanding respiratory diseases since nearly all pulmonary diseases lead to some abnormality of lung tissue mechanics (Suki and Bates, 2011).

Other biological poroelastic applications include, protein-based hydrogels embedded within cells (Galie et al., 2011), orbital soft tissues of the eye (Luboz et al., 2004), brain oedema and hydrocephalus (Li et al., 2010; Wirth and Sobey, 2006), microcirculation of blood and interstitial fluid in the liver lobule (Leungchavaphongse, 2013), and interstitial fluid and tissue in articular cartilage and intervertebral discs (Galbusera et al., 2011; Holmes and Mow, 1990; Mow et al., 1980). Understanding the biomechanics of these tissues has a wide range of useful applications from tracking tumours (Rajagopal et al., 2010) to surgery planning (Luboz et al., 2004).

## 1.2 Numerical challenges

The method that we use for spatially discretising the equations in this work is the finite element method (FEM).

When using the finite element method to solve the poroelastic equations the main challenge is to ensure convergence of the method and prevent numerical instabilities that often manifest themselves in the form of spurious oscillations in the pressure. It has been suggested that this problem is caused by the saddle point structure in the coupled equations resulting in a violation of the famous Ladyzhenskaya-Babuska-Brezzi (LBB) condition, thus highlighting the need for a stable combination of mixed finite elements (Haga et al., 2012).

In addition to this, there has been a need for a method that is able to overcome localised pressure oscillations due to steep pressure gradients in the solution. In particular, when modelling the diseased lung, abrupt changes in tissue properties and heterogeneous airway narrowing are possible. This can result in a patchy ventilation and pressure distribution (Venegas et al., 2005). In this situation existing methods that solve the poroelastic equations using a continuous pressure approximation would struggle to capture the steep gradients in pressure, and result in localised oscillations in the pressure (Phillips and Wheeler, 2008).

Another numerical challenge in practical 3D applications is the algebraic system arising from the finite element discretisation. This can lead to a very large matrix system that has many unknowns and is severely ill-conditioned, making it difficult to solve using standard iterative solvers. Therefore low-order finite element methods that allow for efficient preconditioning are preferable (Ferronato et al., 2010; White and Borja, 2011).

The implementation of finite element codes can also be a challenge. This is especially true when using non-standard elements that are not supported in existing finite element libraries. For example, assembling and calculating higher order stress quantities on discontinuous and non-conforming finite elements in 3D can be particularly difficult. Therefore a method that uses standard and simple to implement elements is very appealing (White and Borja, 2011).

For large deformation applications, common in biology, convergence of the nonlinear coupled problem using Newton's method or other iterative methods is also nontrivial (Ün and Spilker, 2006). This problem can be especially delicate when the nonlinear poroelastic model is tightly coupled to yet another fluid model such as a fluid network model, approximating the airways in the lungs.

## 1.3 Thesis goals

The main goal of this thesis is to rigorously develop a finite element method for solving the linear and nonlinear poroelastic equations. We then plan to demonstrate this methodology by simulating the lung breathing on a realistic geometry. More specific targets are:

1. Develop a practical low-order finite element method for solving the linear poroelastic equations using a discontinuous pressure approximation. Prove theoretical results about the discretisation, including existence and uniqueness, an energy estimate and an optimal a-priori error estimate.

2. Extend the method to a non-linear finite element method to solve the poroelastic equations valid in large deformations.

3. Rigorously test the method using numerous test problems to verify theo-

retical stability and convergence results, and its ability to reliably capture steep pressure gradients.

4. Present a poroelastic model for lung parenchyma coupled to an airway fluid network model, and develop a stable method to numerically solve the coupled model.

5. Solve the computational lung model on a realistic geometry, with boundary conditions extracted from imaging data, to simulate breathing, and evaluate the effect of tissue weakening and airway narrowing on lung function.

## 1.4 Thesis structure

The contributions of each chapter to the thesis are as follows:

**Chapter 2:** We introduce the general theory of poroelasticity valid in large deformations and state the linear poroelastic equations, valid in small deformations.

**Chapter 3:** We outline the basic concepts of the standard continuous Galerkin finite element method. We then discuss mixed problems and their stability requirement. We conclude the chapter by discussing numerical methods currently available to solve the poroelastic equations.

**Chapter 4:** We present a stabilised finite element method for the linear three-field (displacement, fluid flux and pressure) poroelasticity problem. By applying a local pressure jump stabilisation term to the mass conservation equation we avoid pressure oscillations. For the fully-discretised problem we prove existence and uniqueness, an energy estimate and an optimal a-priori error estimate.

**Chapter 5:** We present numerical experiments in 2D and 3D illustrate the con-vergence of the method, and show the effectiveness of the method to overcome spurious pressure oscillations. The added mass effect of the stabilisation term is shown to be negligible in 3D.

**Chapter 6:** We modify the method developed in Chapter 4 to solve the three-field nonlinear quasi-static incompressible poroelasticity problem valid in large deformations. We present the linearisation and discretisation of the equations, and give a detailed account of the implementation. Numerical experiments in 3D verify the method and illustrate its ability to reliably capture steep pressure gradients.

**Chapter 7:** We begin by giving an overview of lung physiology and existing ventilation models. We then present the model assumptions required for the proposed poroelastic lung model, and outline its mathematical formulation and coupling to the airway fluid network. A numerical method is presented to discretise the equations in a monolithic way to ensure unconditional stability. Finally, numerical simulations on a realistic lung geometry that illustrate the coupling between the poroelastic medium and the network flow model are presented. Simulations of tidal breathing are shown to reproduce global physiologically realistic measurements. We also investigate the effect of airway constriction and tissue weakening on the ventilation, tissue stress and alveolar pressure distribution.

**Chapter 8:** We review the main contributions and propose future lines of research.

# Chapter 2

# Poroelasticity theory

Two complementary approaches have been developed for modelling a deformable porous medium. Mixture theory, also known as the Theory of Porous Media (TPM) (Boer, 2005; Bowen, 2010), has its roots in the classical theories of gas mixtures and makes use of a volume fraction concept in which the porous medium is represented by spatially superposed interacting media. An alternative, purely macroscopic approach is mainly associated with the work of Biot, a detailed description can be found in the book by Coussy (2004).

The theory developed by Biot (Biot, 1941) assumes that stress and other related concepts hold at the macro level, such as the fluid flow through the porous matrix. The constitutive equations involve well defined and measurable quantities at the macro level, as for example the permeability. The equations are generally formulated in a Lagrangian description using a macroscale Helmholtz energy potential.

Relationships between the two theories are explored by Coussy et al. (1998). As is most common in biological applications, we use the mixture theory for poroelasticity as outlined in Boer (2005).

## 2.1 Kinematics

Within continuum mixture theory, a poroelastic medium is treated as the superposition of two interacting continua simultaneously occupying the same physical space. The superscript $\alpha \in \{s, f\}$ denotes a quantity related to the solid or fluid, respectively. Before presenting the mixture theory, we give a review of solid mechanics. This will form the basis of the description of the solid skeleton. The following review of continuum mechanics closely follows Chapter 4 in Gonzalez and Stuart (2008), and the standard Poromechanics book by Coussy (2004).

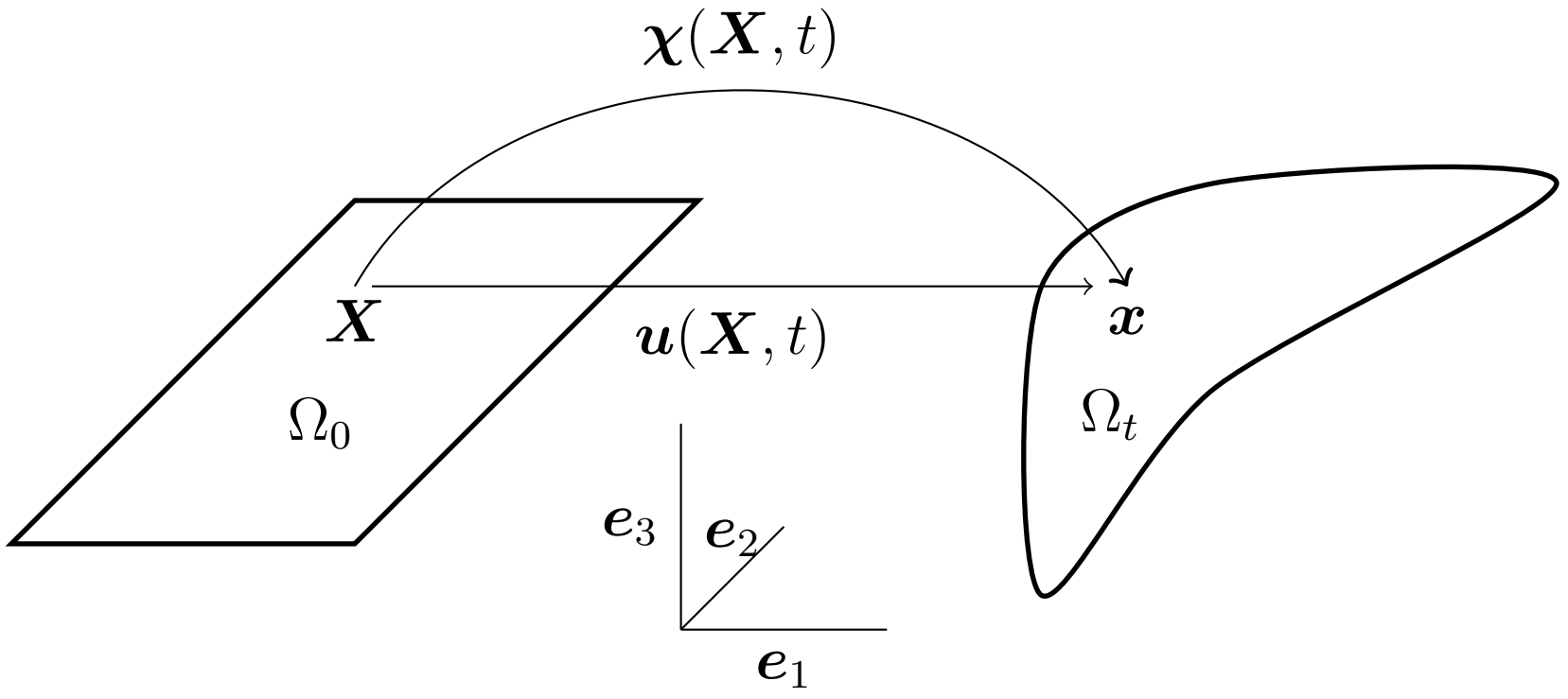

Figure 2.1: Illustration of the solid deformation.

Let the volume $\Omega_0$ be the undeformed Lagrangian (material) reference configuration and let $\boldsymbol{X} = \{X\boldsymbol{e}_1 + Y\boldsymbol{e}_2 + Z\boldsymbol{e}_3\}$ indicate the position of a solid particle in $\Omega_0$ at $t = 0$, where $X, Y$ and $Z$ are the components of the position with respect to the standard orthonormal basis $\{\boldsymbol{e}_1, \boldsymbol{e}_2, \boldsymbol{e}_3\}$ for $\mathbb{R}^3$. The position of a solid particle in the current Eulerian (spatial) configuration $\Omega_t$ is given by $\boldsymbol{x} = \{x\boldsymbol{e}_1 + y\boldsymbol{e}_2 + z\boldsymbol{e}_3\}$, with $\boldsymbol{x} = \boldsymbol{\chi}(\boldsymbol{X}, t)$, shown in Figure 2.1. The deformation map, $\boldsymbol{\chi}(\boldsymbol{X}, t)$, is a continuously differentiable, invertible mapping from $\Omega_0$ to $\Omega_t$. Thus the inverse of the deformation map, $\boldsymbol{\chi}^{-1}(\boldsymbol{x}, t)$, is such that $\boldsymbol{X} = \boldsymbol{\chi}^{-1}(\boldsymbol{x}, t)$. The displacement field is given by

$$\boldsymbol{u}(\boldsymbol{X}, t) = \boldsymbol{\chi}(\boldsymbol{X}, t) - \boldsymbol{X}. \tag{2.1}$$

The deformation gradient tensor is

$$\boldsymbol{F} = \frac{\partial \boldsymbol{\chi}(\boldsymbol{X}, t)}{\partial \boldsymbol{X}} = \boldsymbol{I} + \frac{\partial \boldsymbol{u}(\boldsymbol{X}, t)}{\partial \boldsymbol{X}}, \tag{2.2}$$

and maps a material line element in the reference configuration $d\boldsymbol{X}$, to a line element $d\boldsymbol{x}$ in the current configuration, i.e. $d\boldsymbol{x} = \boldsymbol{F} d\boldsymbol{X}$. The symmetric right Cauchy-Green deformation tensor is given by

$$\boldsymbol{C} = \boldsymbol{F}^T \boldsymbol{F}. \tag{2.3}$$

The Jacobian is defined as

$$J = \det(\boldsymbol{F}), \tag{2.4}$$

and represents the change in an infinitesimal small volume from a reference volume element $d\Omega_0$ to a current configuration volume element $d\Omega_t$

$$d\Omega_t = J d\Omega_0. \tag{2.5}$$

Note that $J > 0$, to avoid self penetration of the body. Also, $\boldsymbol{F}$ is invertible, and it is easy to see that the inverse of the deformation gradient is the deformation gradient of the inverse map

$$\boldsymbol{F}^{-1} = \frac{\partial \boldsymbol{\chi}^{-1}(\boldsymbol{x}, t)}{\partial \boldsymbol{x}} = \frac{\partial \boldsymbol{X}}{\partial \boldsymbol{x}}. \tag{2.6}$$

We denote by $\boldsymbol{V}(\boldsymbol{X}, t)$ the velocity at time $t$ of the material (fixed) solid particle $\boldsymbol{X}$. By definition we have

$$\boldsymbol{V}(\boldsymbol{X}, t) = \frac{\partial}{\partial t} \boldsymbol{\chi}(\boldsymbol{X}, t). \tag{2.7}$$

Similarly, we denote by $\boldsymbol{A}(\boldsymbol{X}, t)$ the acceleration of the material solid particle,

given by

$$\boldsymbol{A}(\boldsymbol{X},t) = \frac{\partial^2}{\partial t^2}\boldsymbol{\chi}(\boldsymbol{X},t). \tag{2.8}$$

We see that the velocity and acceleration of material particles are material fields. Also note that $\frac{\partial}{\partial t}\boldsymbol{u}(\boldsymbol{X},t) = \frac{\partial}{\partial t}\boldsymbol{\chi}(\boldsymbol{X},t)$. We will also require a spatial description of these fields. We denote by $\boldsymbol{v}^s(\boldsymbol{x},t)$ the spatial description of the material solid velocity field, such that

$$\boldsymbol{v}^s(\boldsymbol{x},t) = \left[\frac{\partial}{\partial t}\boldsymbol{\chi}(\boldsymbol{X},t)\right]\bigg|_{\boldsymbol{X}=\chi^{-1}(\boldsymbol{x},t)}. \tag{2.9}$$

Due to the definition of $\boldsymbol{v}^s$ in (2.9) we also have (see section 4.4.4 in Gonzalez and Stuart (2008))

$$\boldsymbol{v}^s(\boldsymbol{x},t)|_{\boldsymbol{x}=\boldsymbol{\chi}(\boldsymbol{X},t)} = \frac{\partial}{\partial t}\boldsymbol{\chi}(\boldsymbol{X},t). \tag{2.10}$$

To simplify the notation we will follow Ateshian et al. (2010) and write

$$\boldsymbol{v}^s(\boldsymbol{x},t) = \frac{\partial}{\partial t}\boldsymbol{\chi}(\boldsymbol{X},t). \tag{2.11}$$

Similarly, for the spatial description of the solid acceleration, we have

$$\boldsymbol{a}^s(\boldsymbol{x},t) = \left[\frac{\partial^2}{\partial t^2}\boldsymbol{\chi}(\boldsymbol{X},t)\right]\bigg|_{\boldsymbol{X}=\chi^{-1}(\boldsymbol{x},t)}. \tag{2.12}$$

Notice that $\boldsymbol{v}^s(\boldsymbol{x},t)$ and $\boldsymbol{a}^s(\boldsymbol{x},t)$ correspond to the velocity and acceleration of the solid material particle whose current coordinates are $\boldsymbol{x}$ at time $t$. The acceleration of the fluid is given by (see section 3.1 in Boer (2005)),

$$\boldsymbol{a}^f = \frac{d^f\boldsymbol{v}^f}{dt} = \frac{\partial}{\partial t}\boldsymbol{v}^f + (\nabla\boldsymbol{v}^f)\boldsymbol{v}^f. \tag{2.13}$$

The **particle derivative of a field** $\mathcal{G}(\boldsymbol{x},t)$ with respect to the particle $\alpha$ ($s$ or

$f$) is given by

$$\frac{d^{\alpha}}{dt}\mathcal{G} = \frac{\partial \mathcal{G}}{\partial t} + (\nabla \mathcal{G})\boldsymbol{v}^{\alpha}, \tag{2.14}$$

where $\nabla(\cdot) = \partial/\partial \boldsymbol{x}(\cdot)$ denotes the partial derivative with respect to the *deformed* configuration. We will use $\nabla$ to denote the spatial gradient in $\Omega_t$ rather than the more explicit $\nabla_{\boldsymbol{x}=\boldsymbol{\chi}(\boldsymbol{X},t)}$. The latter more clearly indicates the dependency of the gradient operator on the deformation $\boldsymbol{\chi}(\boldsymbol{X}, t)$ and highlights the inherent nonlinearity that arises due to the fact that the deformation $\boldsymbol{\chi}(\boldsymbol{X}, t)$ is one of the unknowns. Similarly the deformed domain $\Omega_t$, is a function of the deformation map $\boldsymbol{\chi}$, and therefore incorporates another important nonlinearity.

The **particle derivative of a material volume** with respect to the $\alpha$-constituent is given by (see section 1.3.1 in Coussy (2004))

$$\frac{d^{\alpha}}{dt}\int_{\Omega_t} d\Omega_t = \int_{\Omega_t} \nabla \cdot \boldsymbol{v}^{\alpha} d\Omega_t. \tag{2.15}$$

The particle derivative also applies to a volume integral. Thus, for any quantity $\mathcal{G}$, associated with the $\alpha$ constituent, we have

$$\frac{d^{\alpha}}{dt}\int_{\Omega_t} \mathcal{G} d\Omega_t = \int_{\Omega_t} \left(\frac{d^{\alpha}\mathcal{G}}{dt} + \mathcal{G}\nabla \cdot \boldsymbol{v}^{\alpha}\right) d\Omega_t = \int_{\Omega_t} \left(\frac{\partial \mathcal{G}}{\partial t} + \nabla \cdot \mathcal{G}\boldsymbol{v}^{\alpha}\right) d\Omega_t. \tag{2.16}$$

This is commonly known as the Reynolds transport theorem. In the last step of (2.16) we have used the identity $\nabla \cdot (\psi \boldsymbol{s}) = \boldsymbol{s} \cdot \nabla\psi + \psi\nabla \cdot \boldsymbol{s}$ for some scalar $\psi$ and vector $\boldsymbol{s}$.

## 2.2 Volume fractions

We restrict our attention to saturated porous media which are assumed to consist of solid and fluid parts. The fluid accounts for volume fractions $\phi_0(\boldsymbol{X}, t = 0)$ and $\phi(\boldsymbol{x}, t)$ of the total volume in the reference and the current and deformed configurations respectively, where $\phi$ is known as the porosity. The fractions for the solid are therefore $1 - \phi_0$ and $1 - \phi$ in the reference and the current configuration respectively. For a mixture the density in the current configuration is given by

$$\rho = \rho^s(1 - \phi) + \rho^f \phi \quad \text{in } \Omega_t, \tag{2.17}$$

where $\rho^s$ and $\rho^f$ are the densities of the fluid and solid, respectively. We assume that both the solid and the fluid are incompressible so that $\rho^s = \rho_0^s$ and $\rho^f = \rho_0^f$. For notational convenience we also define

$$\hat{\rho}^s = \rho^s(1 - \phi), \tag{2.18}$$

and

$$\hat{\rho}^f = \rho^f \phi. \tag{2.19}$$

Due to mass conservation and the incompressibility of both the solid and the fluid phase we have

$$J = \frac{1 - \phi_0}{1 - \phi}, \tag{2.20}$$

where $J$ represents the change in volume of the solid skeleton. The solid skeleton includes the solid (tissue in biological applications) and the voids occupied by the fluid. Note that although the solid is assumed to be incompressible the solid skeleton is able to change in volume, since fluid can enter or leave the solid skeleton.

## 2.3 Conservation of mass

When no mass change occurs, neither for the solid skeleton or the fluid contained in $\Omega_t$, using the Reynolds transport theorem (2.16), the balance of mass, for a volume $V(t)$ that moves with the deforming poroelastic medium, can be

expressed as

$$\frac{d^s}{dt}\int_{V(t)}(1-\phi)\rho^s d\Omega_t = \int_{V(t)}\left(\frac{\partial(1-\phi)\rho^s}{\partial t} + \nabla\cdot((1-\phi)\rho^s\boldsymbol{v}^s)\right)d\Omega_t,$$

$$\frac{d^f}{dt}\int_{V(t)}\phi\rho^f d\Omega_t = \int_{V(t)}\left(\frac{\partial\phi\rho^f}{\partial t} + \nabla\cdot(\phi\rho^f\boldsymbol{v}^f)\right)d\Omega_t.$$

Thus, the balance of mass for the solid is given by

$$\frac{\partial(1-\phi)\rho^s}{\partial t} + \nabla\cdot((1-\phi)\rho^s\boldsymbol{v}^s) = 0 \quad \text{in } \Omega_t, \tag{2.21}$$

where $\boldsymbol{v}^s$ is the velocity vector of the solid. Similarly, the balance of mass for the fluid is given by

$$\frac{\partial\phi\rho^f}{\partial t} + \nabla\cdot(\phi\rho^f\boldsymbol{v}^f) = \rho^f g \quad \text{in } \Omega_t, \tag{2.22}$$

where $\boldsymbol{v}^f$ is the velocity vector of the fluid and $g$ is a general source or sink term. Noting that $\rho^s$ and $\rho^f$ are constants (in space and time), these can be factored out of equations (2.21) and (2.22). Adding these two equations then provides the mass balance or continuity equation of the mixture (see section 8.3 in Boer (2005)),

$$\nabla\cdot((1-\phi)\boldsymbol{v}^s) + \nabla\cdot(\phi\boldsymbol{v}^f) = g \quad \text{in } \Omega_t. \tag{2.23}$$

## 2.4 Conservation of momentum

The balance law of linear momentum for each individual constituent is given by

$$\frac{d^\alpha}{dt}\int_{V(t)}\hat{\rho}^\alpha\boldsymbol{v}^\alpha d\Omega_t = \int_{V(t)}\nabla\cdot\boldsymbol{\sigma}^\alpha + \hat{\rho}^\alpha\boldsymbol{f} + \hat{\boldsymbol{p}}^\alpha + \Theta^\alpha\boldsymbol{v}^\alpha\, d\Omega_t. \tag{2.24}$$

Here $\boldsymbol{\sigma}^\alpha$ is the Cauchy stress tensor of the $\alpha$ constituent, $\boldsymbol{f}$ is a volume force acting on the constituents, $\hat{\boldsymbol{p}}^\alpha$ are interaction forces representing frictional interactions between the solid and fluid, defined later in section 7.5.1, and $\Theta^\alpha \boldsymbol{v}^\alpha$ is the variation of momentum due to the $\alpha$ constituent source term (Chapelle and Moireau, 2014). Note that from (2.21) and (2.22) that we have $\Theta^s = 0$ and $\Theta^f = \rho^f g$. Using the first step of the Reynolds transport theorem (2.16), and the chain rule, we obtain

$$\nabla \cdot \boldsymbol{\sigma}^\alpha + \hat{\rho}^\alpha \boldsymbol{f} + \hat{\boldsymbol{p}}^\alpha + \Theta^\alpha \boldsymbol{v}^\alpha = \hat{\rho}^\alpha \boldsymbol{a}^\alpha + \boldsymbol{v}^\alpha \left( \frac{d^\alpha \hat{\rho}^\alpha}{dt} + \hat{\rho}^\alpha \nabla \cdot \boldsymbol{v}^\alpha \right) \quad \text{in } \Omega_t, \tag{2.25}$$

where $\boldsymbol{a}^\alpha$ are acceleration vectors of the constituents. Since each constituent exerts an equal and opposite interaction force on the other,

$$\hat{\boldsymbol{p}}^s + \hat{\boldsymbol{p}}^f = 0. \tag{2.26}$$

## 2.5 Constitutive relations

The interaction force is given by (see (Coussy, 2004, eqn. (3.49)))

$$\hat{\boldsymbol{p}}^s = -\hat{\boldsymbol{p}}^f = -p\nabla\phi + \phi^2 \boldsymbol{k}^{-1} \cdot (\boldsymbol{v}^f - \boldsymbol{v}^s), \tag{2.27}$$

where $\boldsymbol{k}$ is the (dynamic) permeability tensor. The first term, $p\nabla\phi$, accounts for the pressure effect resulting from the variation of the section offered to the fluid flow, and the second term, $\phi^2 \boldsymbol{k} \cdot (\boldsymbol{v}^f - \boldsymbol{v}^s)$, describes the viscous resistance opposed by the shear stress to the fluid flow from the drag at the internal walls of the porous network (Coussy, 2004). This particular choice for the interaction force means that the momentum balance for the fluid flow can later be reduced to the well known Darcy law.

The permeability tensor in the current configuration is given by

$$\boldsymbol{k} = J^{-1}\boldsymbol{F}\boldsymbol{k}_0(\boldsymbol{\chi})\boldsymbol{F}^T, \tag{2.28}$$

where $\boldsymbol{k}_0(\boldsymbol{\chi})$ is the permeability in the reference configuration, which may be chosen to be some (nonlinear) function dependent on the deformation. Examples of deformation dependent permeability tensors for biological tissues can be found in Holmes and Mow (1990); Kowalczyk and Kleiber (1994); Lai and Mow (1980).

The solid stress tensor is given by the effective stress principle (see eqn. (8.62) in Boer (2005)),

$$\boldsymbol{\sigma}^s = \boldsymbol{\sigma}^s_e - (1-\phi)\boldsymbol{I}p, \tag{2.29}$$

where $\boldsymbol{\sigma}^s_e$ is the effective stress tensor given by

$$\boldsymbol{\sigma}^s_e = \frac{1}{J}\boldsymbol{F} \cdot 2\frac{\partial W(\boldsymbol{\chi})}{\partial \boldsymbol{C}} \cdot \boldsymbol{F}^T. \tag{2.30}$$

Here $W(\boldsymbol{\chi})$ denotes a strain-energy law (hyperelastic Helmholtz energy functional) dependent on the deformation of the solid. The fluid stress tensor can be written as (see (Boer, 2005, eqn. (8.63)))

$$\boldsymbol{\sigma}^f = \boldsymbol{\sigma}^f_{vis} - \phi\boldsymbol{I}p, \tag{2.31}$$

where $\boldsymbol{\sigma}^f_{vis}$ denotes the viscous stress tensor of the fluid, given by (see Boer (2005, eqn. (6.145)))

$$\boldsymbol{\sigma}^f_{vis} = \mu_f\phi(\nabla\boldsymbol{v}_f + (\nabla\boldsymbol{v}_f)^T - \frac{2}{3}\nabla\cdot\boldsymbol{v}_f), \tag{2.32}$$

where $\mu_f$ is the dynamic viscosity of the fluid.

Summing the conservation laws (2.25) for its constituents and applying the constitutive relations, the conservation of linear momentum for the mixture is

$$\hat{\rho}^s \boldsymbol{a}^s + \hat{\rho}^f \boldsymbol{a}^f + \boldsymbol{v}^s \left( \frac{d^s \hat{\rho}^s}{dt} + \hat{\rho}^s \nabla \cdot \boldsymbol{v}^s \right) + \boldsymbol{v}^f \left( \frac{d^f \hat{\rho}^f}{dt} + \hat{\rho}^f \nabla \cdot \boldsymbol{v}^f \right)$$
$$= \nabla \cdot (\boldsymbol{\sigma}_e + \boldsymbol{\sigma}_{vis} - p\boldsymbol{I}) + \rho \boldsymbol{f} + g\boldsymbol{v}^f \quad \text{in } \Omega_t. \quad (2.33)$$

Applying (2.21) and (2.22), along with applications of (2.14), we get

$$\hat{\rho}^s \boldsymbol{a}^s + \hat{\rho}^f \boldsymbol{a}^f = \nabla \cdot (\boldsymbol{\sigma}_e + \boldsymbol{\sigma}_{vis} - p\boldsymbol{I}) + \rho \boldsymbol{f} \quad \text{in } \Omega_t. \quad (2.34)$$

The momentum equation for the fluid flow can be identified from (2.25) with $\alpha = f$ as

$$\hat{\rho}^f \boldsymbol{a}^f = \nabla \cdot (\boldsymbol{\sigma}^f_{vis} - \phi p \boldsymbol{I}) + \hat{\rho}^f \boldsymbol{f} + p \nabla \phi - \phi^2 \boldsymbol{k}^{-1} (\boldsymbol{v}^f - \boldsymbol{v}^s) \, \text{in} \, \Omega_t. \quad (2.35)$$

## 2.6 Summary of the general poroelasticity model

We consider $\Omega_t$ to be a bounded domain in $\mathbb{R}^2$ or $\mathbb{R}^3$, and for the purpose of defining boundary conditions, $\partial\Omega_t = \Gamma_D \cup \Gamma_N$ for displacement and stress boundary conditions and $\partial\Omega_t = \Gamma_P \cup \Gamma_F$ for pressure and flux boundary conditions, with outward pointing unit normal $\boldsymbol{n}$. The strong problem for the full mixture theory model is to find $\boldsymbol{\chi}(\boldsymbol{X},t)$, $\boldsymbol{v}^f(\boldsymbol{x},t)$ and $p(\boldsymbol{x},t)$ such that

$$\hat{\rho}^s\boldsymbol{a}^s + \hat{\rho}^f\boldsymbol{a}^f = \nabla\cdot(\boldsymbol{\sigma}_e + \boldsymbol{\sigma}_{vis} - p\boldsymbol{I}) + \rho\boldsymbol{f} \quad \text{in } \Omega_t, \tag{2.36a}$$

$$\hat{\rho}^f\boldsymbol{a}^f = \nabla\cdot(\boldsymbol{\sigma}^f_{vis} - \phi p\boldsymbol{I}) + p\nabla\phi - \phi\boldsymbol{k}^{-1}(\boldsymbol{v}^f - \boldsymbol{v}^s) + \hat{\rho}^f\boldsymbol{f} \quad \text{in } \Omega_t, \tag{2.36b}$$

$$\nabla\cdot((1-\phi)\boldsymbol{v}^s) + \nabla\cdot(\phi\boldsymbol{v}^f) = g \quad \text{in } \Omega_t, \tag{2.36c}$$

$$\boldsymbol{\chi}(\boldsymbol{X},t)|_{\boldsymbol{X}=\boldsymbol{\chi}^{-1}(\boldsymbol{x},t)} = \boldsymbol{X} + \boldsymbol{u}_D \quad \text{on } \Gamma_D, \tag{2.36d}$$

$$(\boldsymbol{\sigma}_e + \boldsymbol{\sigma}_{vis} - p\boldsymbol{I})\boldsymbol{n} = \boldsymbol{t}_N \quad \text{on } \Gamma_N, \tag{2.36e}$$

$$\boldsymbol{v}^f = \boldsymbol{v}^f_D \quad \text{on } \Gamma_F, \tag{2.36f}$$

$$(\boldsymbol{\sigma}_{vis} - \phi p\boldsymbol{I})\boldsymbol{n} = \boldsymbol{s}_P \quad \text{on } \Gamma_P, \tag{2.36g}$$

$$\boldsymbol{\chi}(0) = \boldsymbol{X}, \quad \boldsymbol{v}^s(0) = \boldsymbol{v}^{s0}, \quad \boldsymbol{v}^f(0) = \boldsymbol{v}^{f0} \quad \text{in } \Omega_0. \tag{2.36h}$$

We have also summarised all the variables and corresponding equations in Table 2.1.

| Unknown | Notation | Equation | |
|---|---|---|---|
| **Primary variables** | | **Primary equations (general model)** | |
| Motion of the solid | $\boldsymbol{\chi}$ | $\hat{\rho}^s \boldsymbol{a}^s + \hat{\rho}^f \boldsymbol{a}^f = \nabla \cdot (\boldsymbol{\sigma}_e + \boldsymbol{\sigma}_{vis} - p\boldsymbol{I}) + \rho \boldsymbol{f}$ | (2.34) |
| Fluid velocity | $\boldsymbol{v}^f$ | $\hat{\rho}^f \boldsymbol{a}^f = \nabla \cdot (\boldsymbol{\sigma}^f_{vis} - \phi p \boldsymbol{I}) + p\nabla\phi - \phi^2 \boldsymbol{k}^{-1}(\boldsymbol{v}^f - \boldsymbol{v}^s) + \hat{\rho}^f \boldsymbol{f}$ | (2.35) |
| Pressure of the fluid | $p$ | $\nabla \cdot ((1-\phi)\boldsymbol{v}^s) + \nabla \cdot (\phi \boldsymbol{v}^f) = g$ | (2.23) |
| **Secondary variables** | | **Secondary equations** | |
| Deformation gradient tensor | $\boldsymbol{F}$ | $\boldsymbol{F} = \frac{\partial}{\partial \boldsymbol{X}} \boldsymbol{\chi}(\boldsymbol{X}, t)$ | (2.2) |
| Right Cauchy-Green tensor | $\boldsymbol{C}$ | $\boldsymbol{C} = \boldsymbol{F}^T \boldsymbol{F}$ | (2.3) |
| Jacobian | $J$ | $J = \det(\boldsymbol{F})$ | (2.4) |
| Velocity of the solid | $\boldsymbol{v}^s$ | $\boldsymbol{v}^s(\boldsymbol{x}, t)\vert_{\boldsymbol{x}=\boldsymbol{\chi}(\boldsymbol{X},t)} = \frac{\partial}{\partial t} \boldsymbol{\chi}(\boldsymbol{X}, t)$ | (2.10) |
| Acceleration of the solid | $\boldsymbol{a}^s$ | $\boldsymbol{a}^s(\boldsymbol{x}, t)\vert_{\boldsymbol{x}=\chi(\boldsymbol{X},t)} = \frac{\partial^2}{\partial t^2} \boldsymbol{\chi}(\boldsymbol{X}, t)$ | (2.12) |
| Acceleration of the fluid | $\boldsymbol{a}^f$ | $\boldsymbol{a}^f = \frac{\partial}{\partial t} \boldsymbol{v}^f + (\nabla \boldsymbol{v}^f)\boldsymbol{v}^f$ | (2.13) |
| Porosity | $\phi$ | $\phi = 1 - \frac{1-\phi_0}{J}$ | (2.20) |
| Mixture density | $\rho$ | $\rho = \rho^s(1-\phi) + \rho^f \phi$ | (2.17) |
| Eulerian solid density | $\hat{\rho}_s$ | $\hat{\rho}^s = \rho^s(1-\phi)$ | (2.18) |
| Eulerian fluid density | $\hat{\rho}_f$ | $\hat{\rho}^f = \rho^f \phi$ | (2.19) |
| **Constitutive variables** | | **Constitutive equations** | |
| Solid elastic stress tensor | $\boldsymbol{\sigma}_e$ | $\boldsymbol{\sigma}^s_e = \frac{1}{J} \boldsymbol{F} \cdot 2 \frac{\partial W(\boldsymbol{\chi})}{\partial \boldsymbol{C}} \cdot \boldsymbol{F}^T$ | (2.30) |
| Fluid viscous stress tensor | $\boldsymbol{\sigma}_{vis}$ | $\boldsymbol{\sigma}^f_{vis} = \mu_f \phi(\nabla \boldsymbol{v}_f + (\nabla \boldsymbol{v}_f)^T - \frac{2}{3}\nabla \cdot \boldsymbol{v}_f)$ | (2.32) |
| Permeability tensor | $\boldsymbol{k}$ | $\boldsymbol{k} = J^{-1} \boldsymbol{F} \boldsymbol{k}_0(\boldsymbol{\chi}) \boldsymbol{F}^T$ | (2.28) |

Table 2.1: Recapitulating the unknowns and equations of the general poroelasticity model.

## 2.7 Simplification and reformulation of the model

To arrive at the quasi-static, fully saturated, incompressible three-field large deformation poroelasticity model, we will now ignore inertia forces (left hand side of (2.34) and (2.35)), and ignore the viscous shear stress in the fluid ($\boldsymbol{\sigma}^f_{vis}$ in (2.35)). Justifications for making these modelling assumptions with respect to the proposed lung model will be given in section 7.4. After making these assumptions, and rewriting the equations in terms of the fluid flux, given by

$$\boldsymbol{z} = \phi(\boldsymbol{v}^f - \boldsymbol{v}^s), \tag{2.37}$$

the resulting problem is to find $\boldsymbol{\chi}(\boldsymbol{X},t)$, $\boldsymbol{z}(\boldsymbol{x},t)$ and $p(\boldsymbol{x},t)$ such that

$$-\nabla\cdot(\boldsymbol{\sigma}_e - p\boldsymbol{I}) = \rho\boldsymbol{f} \quad \text{in } \Omega_t, \tag{2.38a}$$

$$\boldsymbol{k}^{-1}\boldsymbol{z} + \nabla p = \rho^f\boldsymbol{f} \quad \text{in } \Omega_t, \tag{2.38b}$$

$$\nabla\cdot(\boldsymbol{v}^s + \boldsymbol{z}) = g \quad \text{in } \Omega_t, \tag{2.38c}$$

$$\boldsymbol{\chi}(\boldsymbol{X},t)|_{\boldsymbol{X}=\boldsymbol{\chi}^{-1}(\boldsymbol{x},t)} = \boldsymbol{X} + \boldsymbol{u}_D \quad \text{on } \Gamma_D, \tag{2.38d}$$

$$(\boldsymbol{\sigma}_e - p\boldsymbol{I})\boldsymbol{n} = \boldsymbol{t}_N \quad \text{on } \Gamma_N, \tag{2.38e}$$

$$\boldsymbol{z}\cdot\boldsymbol{n} = q_D \quad \text{on } \Gamma_F, \tag{2.38f}$$

$$p = p_D \quad \text{on } \Gamma_P, \tag{2.38g}$$

$$\boldsymbol{\chi}(0) = \boldsymbol{X}, \quad \text{in } \Omega_0. \tag{2.38h}$$

This is the large deformation model we will consider from here onwards.

## 2.8 Linear poroelasticity

To allow us to perform rigorous analysis of the proposed finite element scheme presented in Chapter 4, we will now assume small deformations to yield a linear model of poroelasticity. This model is often referred to as the 'Biot model' in the geomechanics community and contains some additional terms. We will introduce the full Biot model here for use with a 2D cantilever bracket problem later tested in section 5.5, and to highlight that any subsequent theory developed in later chapters can be extended to the full Biot model. The governing equations of the Biot model, with displacement $\boldsymbol{u}$, fluid flux $\boldsymbol{z}$, and pressure $p$ as primary variables are summarised below:

$$-\nabla \cdot \boldsymbol{\sigma} = \boldsymbol{f} \quad \text{in } \Omega, \tag{2.39a}$$

$$\boldsymbol{k}^{-1}\boldsymbol{z} + \nabla p = \boldsymbol{b} \quad \text{in } \Omega, \tag{2.39b}$$

$$\nabla \cdot \boldsymbol{z} + \frac{\partial}{\partial t}(\alpha \nabla \cdot \boldsymbol{u} + c_0 p) = g \quad \text{in } \Omega, \tag{2.39c}$$

$$\boldsymbol{u} = \boldsymbol{u}_D \quad \text{on } \Gamma_D, \tag{2.39d}$$

$$\boldsymbol{\sigma n} = \boldsymbol{t}_N \quad \text{on } \Gamma_N, \tag{2.39e}$$

$$p = p_D \quad \text{on } \Gamma_P, \tag{2.39f}$$

$$\boldsymbol{z} \cdot \boldsymbol{n} = q_D \quad \text{on } \Gamma_F, \tag{2.39g}$$

$$\boldsymbol{u}(0) = \boldsymbol{u}^0, \quad p(0) = p^0, \quad \text{in } \Omega. \tag{2.39h}$$

Here $\boldsymbol{\sigma}$ is the total stress tensor given by $\boldsymbol{\sigma} = \lambda \text{tr}(\boldsymbol{\epsilon}(\boldsymbol{u}))\boldsymbol{I} + 2\mu_s \boldsymbol{\epsilon}(\boldsymbol{u}) - \alpha p \boldsymbol{I}$, with the linear strain tensor defined as $\boldsymbol{\epsilon}(\boldsymbol{u}) = \frac{1}{2}\left(\nabla \boldsymbol{u} + (\nabla \boldsymbol{u})^T\right)$, $g$ is the fluid source term, $\boldsymbol{f}$ is the body force on the mixture, and $\boldsymbol{b}$ is the body force on the fluid. Here $\Omega$ is a bounded domain in $\mathbb{R}^2$ or $\mathbb{R}^3$, and for the purpose of defining boundary conditions, $\partial\Omega = \Gamma_D \cup \Gamma_N$ for displacement and stress boundary conditions and $\partial\Omega = \Gamma_P \cup \Gamma_F$ for pressure and flux boundary conditions, with outward pointing unit normal $\boldsymbol{n}$.

The momentum and mass conservation equations are coupled through the Biot-Willis constant, $\alpha \in (0, 1]$, and the non-negative constrained specific storage coefficient $c_0 \geq 0$. The increment $\eta$ of fluid volume per unit volume of porous mixture (soil in Biot (1941)) may be written as: $\eta = \alpha \nabla \cdot \boldsymbol{u} + c_0 p$. From this one can observe that $c_0 p$ measures the amount of fluid that can be injected into a fixed material volume under pressure, and $\alpha \nabla \cdot \boldsymbol{u}$ represents the additional amount of fluid content that can be squeezed out due to the local change in volume. (Lipnikov, 2002; Phillips, 2005; Showalter, 2000). The parameters are summarised in Table 2.2.

| **Parameter** | |
|---|---|
| Lamé's first parameter | $\lambda$, |
| Lamé's second parameter (shear modulus) | $\mu_s$, |
| Dynamic permeability tensor | $\boldsymbol{k}$, |
| Biot-Willis constant | $\alpha$, |
| Constrained specific storage coefficient | $c_0$. |

Table 2.2: Poroelasticity parameters.

A derivation and more detailed explanation of these equations can be found in Phillips and Wheeler (2007a) and Showalter (2000). In this work we will mainly consider a simplification of the full Biot model (2.39), by setting $\alpha = 1$ and $c_0 = 0$. This yields a fully incompressible poroelastic model that retains all the numerical difficulties associated with approximating the original system of equations (2.39), see Remark 1. The linear fully saturated and incompressible poroelastic model is given by:

$$-(\lambda + \mu_s)\nabla(\nabla \cdot \boldsymbol{u}) - \mu_s \nabla^2 \boldsymbol{u} + \nabla p = \boldsymbol{f} \quad \text{in } \Omega, \tag{2.40a}$$

$$\boldsymbol{k}^{-1}\boldsymbol{z} + \nabla p = \boldsymbol{b} \quad \text{in } \Omega, \tag{2.40b}$$

$$\nabla \cdot (\boldsymbol{u}_t + \boldsymbol{z}) = g \quad \text{in } \Omega, \tag{2.40c}$$

$$\boldsymbol{u} = \boldsymbol{u}_D \quad \text{on } \Gamma_D, \tag{2.40d}$$

$$\boldsymbol{\sigma n} = \boldsymbol{t}_N \quad \text{on } \Gamma_N, \tag{2.40e}$$

$$p = p_D \quad \text{on } \Gamma_P, \tag{2.40f}$$

$$\boldsymbol{z} \cdot \boldsymbol{n} = q_D \quad \text{on } \Gamma_F, \tag{2.40g}$$

$$\boldsymbol{u}(0) = \boldsymbol{u}^0 \quad \text{in } \Omega, \tag{2.40h}$$

where $\boldsymbol{u}_t$ denotes $\frac{\partial \boldsymbol{u}}{\partial t}$. This model is the small deformation version of the simpli-fied and reformulated large deformation poroelasticity model (2.38), and will be the small deformation model considered from here onwards.

**Remark 1.** *The extension of the theoretical results presented in Chapter 4 from (2.40) to the full Biot equations (2.39), with* $\alpha \in \mathbb{R}_{>0}$ *and* $c_0 \in \mathbb{R}_{>0}$ *is straightforward. In the analysis in Chapter 4, the constant* $\alpha$ *would just get absorbed by a general constant* $C$*. When* $c_0 > 0$*, an additional pressure term is introduced into the mass conservation equation. Since this term is coercive, it only improves the stability of the system.*

# Chapter 3

# Finite element method

## 3.1 Introduction

A large proportion of the mathematical models in science and engineering take the form of differential equations. Only in the simplest cases, or under strong assumptions, is it possible to find exact analytical solutions to the equations in the model. Numerical methods are an established means of solving differential equations that are of practical interest in a variety of applied problems. Finite difference, finite volume and finite element methods are the most widely used of these methods. The basic idea is to replace the infinite-dimensional problem by a finite-dimensional approximation, which is, generally speaking, easier to compute. Finite element methods are based on weakening the restrictions on the solution space in the continuous setting, and searching for the approximate solution in the subspace which spans basis functions supported on small regions inside the domain. These methods are well-suited to solving problems on complex domains, and are therefore widely used in practical applications. In this work we consider only finite element methods (FEMs) for solving partial differential equations. This chapter comprises an overview of several theoretical and practical aspects of classical FEMs. The theory and notation presented here are essential in developing the techniques that form the core of this thesis. Most of the work presented in this chapter is based on work already presented in Arthurs (2012); Asner (2013); Bernabeu (2011); Brenner and Scott (2008); Brezzi and Fortin (1991). We conclude this chapter by discussing numerical methods currently available to solve the poroelastic equations.

## 3.2 Norms and spaces

Let $\Omega$ be a bounded domain in $\mathbb{R}^2$ or $\mathbb{R}^3$, and $\partial\Omega$ be the associated boundary. The space of square integrable functions is then given by

$$L^2(\Omega) = \left\{ u : \int_\Omega |u(x)|^2 \mathrm{d}x < \infty \right\},$$

with norm

$$\|u\|_{0,\Omega} = \left\{ \int_\Omega |u(x)|^2 \mathrm{d}x \right\}^{1/2}.$$

This space is equipped with the inner product

$$(u, v)^{1/2} = \int_\Omega u(x)v(x)\mathrm{d}x,$$

such that $\|u\|_{0,\Omega} = (u, v)^{1/2}$. Throughout this thesis we shall frequently refer to the Sobolev spaces $H^1(\Omega)$ and $H^2(\Omega)$. The definitions of these are as follows:

$$H^1(\Omega) = \left\{ u \in L^2(\Omega) : \frac{\partial u}{\partial x_j} \in L^2(\Omega),\ j = 1, \ldots, n, \right\},$$

$$H^2(\Omega) = \left\{ u \in L^2(\Omega) : \frac{\partial u}{\partial x_j} \in L^2(\Omega),\ j = 1, \ldots, n, \right.$$
$$\left. \frac{\partial^2 u}{\partial x_i \partial x_j} \in L^2(\Omega),\ i, j = 1, \ldots, n \right\}.$$

The corresponding norms are defined as

$$\|u\|_{1,\Omega} = \left\{ \|u\|_{0,\Omega}^2 + \sum_{j=1}^{n} \left\| \frac{\partial u}{\partial x_j} \right\|_{0,\Omega}^2 \right\}^{1/2},$$

$$\|u\|_{2,\Omega} = \left\{ \|u\|_{0,\Omega}^2 + \sum_{j=1}^{n} \left\| \frac{\partial u}{\partial x_j} \right\|_{0,\Omega}^2 + \sum_{i,j=1}^{n} \left\| \frac{\partial^2 u}{\partial x_i \partial x_j} \right\|_{0,\Omega}^2 \right\}^{1/2}.$$

We also define the divergence space

$$H_{div}(\Omega) = \left\{ \boldsymbol{v} \in L^2(\Omega) : \nabla \cdot \boldsymbol{v} \in L^2(\Omega) \right\}.$$

The set of functions of $L^2(\partial\Omega)$ which are traces of functions of $H^1(\Omega)$ onto the boundary, constitutes a subspace of $L^2(\partial\Omega)$ denoted by $H^{1/2}(\partial\Omega)$.

We will also briefly use linear and bounded functionals. For a continuous linear functional, $L : X \to \mathbb{R}$, the dual norm is defined as:

$$||L||_{X'} := \sup_{0 \neq v \in X} \frac{L(v)}{||v||_X},$$

where $X$ denotes a normed space, for example $H^1$ resulting in the norm for the dual space $H^{-1}$, see section 1.7 in Brenner and Scott (2008) for details. Similary the dual spaces $H^{-1/2}$ and $H_{div}^{-1}$ can be defined. We define the following norms for continuous in time functions $u$ such that the norm $L^2(0,T;X)$ satisfies

$$||u||_{L^2(X)} = \left( \int_0^T ||u(\cdot, s)||_X^2 \, ds \right)^{1/2},$$

and the norm $L^\infty(0,T;X)$ satisfies

$$||u||_{L^\infty(X)} = \sup \left\{ ||u(\cdot, s)||_X : s \in [0,T] \right\},$$

where $X$ is any given function space over $\Omega$. We partition $[0,T]$ into $N$ evenly spaced non-overlapping regions $(t_{n-1}, t_n]$, $n = 1, 2, \ldots, N$. For any sufficiently smooth function $u(x,t)$ we define $u^n(x) = u(x, t_n)$. Let the discrete approximation for all time to be the piecewise constant in time functions $v(x,t) = v^n(x)$ for $t \in (t_{n-1}, t_n]$. For such piecewise constant in time functions, $v$, we define the norms

$$||v||_{L^2(X)} = \left( \sum_{n=1}^{N} \Delta t ||v^n||_X^2 \right)^{1/2},$$

and

$$||v||_{L^\infty(X)} = \max \left\{ ||v^n||_X, n = 1, 2, ..., N \right\}.$$

## 3.3 Model problem

It is instructive to begin at a simple level and proceed by incrementally adding to the complexity of the equations we are discretising when explaining the use of the FEM, so we begin by considering the classical heat equation: given $T > 0$, for $t \in [0, T]$ find $u(x,t)$ such that

$$\frac{\partial u}{\partial t} - \nabla \cdot \nabla u = 0 \quad \text{in } \Omega, \tag{3.1a}$$

$$\boldsymbol{n} \cdot \nabla u = g_N \quad \text{on } \Gamma_N, \tag{3.1b}$$

$$u = g_D \quad \text{on } \Gamma_D, \tag{3.1c}$$

$$u(x, 0) = u^0(x) \quad \text{in } \Omega. \tag{3.1d}$$

Here $\Omega$ is a bounded domain in $\mathbb{R}^2$ or $\mathbb{R}^3$, with boundary $\partial\Omega = \Gamma_N \cup \Gamma_D$, that has an outward pointing unit normal $\boldsymbol{n}$. The initial condition is given by $u^0(x)$. In the case where $g_N = 0$, system (3.1) can describe the evolution of heat in an object with geometry described by $\Omega$, where we have perfect thermal insulation on $\Gamma_N$ and fixed temperature distributions given by the function $g_D$ defined on the boundary due to some part of the environment with fixed temperature contacting the object along $\Gamma_D$.

### 3.3.1 Weak formulation

The strong form of (3.1) requires $u$ to be at least twice differentiable. To weaken the regularity restrictions we multiply equation (3.1a) by an arbitrary function $v$, called a test function, and integrate over $\Omega$:

$$\left(\frac{\partial u}{\partial t}, v\right) - (\nabla \cdot \nabla u, v) = 0.$$

Applying the divergence theorem, this equation can be rewritten as

$$\begin{aligned}\left(\frac{\partial u}{\partial t}, v\right) &- (\nabla u \cdot \boldsymbol{n}, v)_{\partial\Omega} + (\nabla u, \nabla v) \\ &= \left(\frac{\partial u}{\partial t}, v\right) - (\nabla u \cdot \boldsymbol{n}, v)_{\Gamma_D} - (g_N, v)_{\Gamma_N} + (\nabla u, \nabla v) = 0.\end{aligned}$$

Here $(\cdot,\cdot)_{\Gamma_N}$ and $(\cdot,\cdot)_{\Gamma_D}$ denote the inner product taken over $\Gamma_N$ and $\Gamma_D$, respectively. Taking note of the Dirichlet condition (3.1c), and letting $v = 0$ on $\Gamma_D$, we arrive at the following equation:

$$\left(\frac{\partial u}{\partial t}, v\right) + (\nabla u, \nabla v) = (g_N, v)_{\Gamma_N}.$$

Note that in this equation the second derivatives of $u$ need not exist. With that in mind, both the solution and the test functions can come from the space $H^1(\Omega)$, as long as they satisfy the appropriate Dirichlet boundary conditions. For convenience we will use the notation $X_D = \{v \in H^1(\Omega) | v = u_D \text{ on } \Gamma_D\}$ and $X_0 = \{v \in H^1(\Omega) | v = 0 \text{ on } \Gamma_D\}$. The weak formulation of (3.1a) is as follows: Find $u \in X_D$ such that

$$\left(\frac{\partial u}{\partial t}, v\right) + (\nabla u, \nabla v) = (g_N, v)_{\Gamma_N} \quad \forall v \in X_0. \tag{3.2}$$

### 3.3.2 Time discretisation

We also need to choose a method of treating the time derivative. In this work, we do so using backward Euler difference quotients, and so we make the approximation $u_t(x, t + \Delta t) \approx \frac{u(x,t+\Delta t)-u(x,t)}{\Delta t}$ for some constant time step $\Delta t$. We write $u(x)^n$ for the the temporally-semidiscrete approximation to $u(x, n\Delta t)$, and our numerical scheme will yield approximations at times $t = 0, \Delta t, 2\Delta t, ..., T$. Inserting this difference quotient and assuming that $\Delta T$ divides $T$, equation (3.3) becomes: for $n = 1, 2, ..., \frac{T}{\Delta t}$, find $u^n \in X_D$ such that

$$(u^n, v) + \Delta t (\nabla u^n, \nabla v) = \Delta t (g_N, v)_{\Gamma_N} + \left(u^{n-1}, v\right) \quad \forall v \in X_0. \tag{3.3}$$

### 3.3.3 Spatial finite element discretisation

In order to solve this problem numerically, we must make it finite dimensional by discretising it suitably. The finite element approximation space is constructed as follows: first, the problem domain is partitioned into small element domains, and second, the element is defined by prescribing for each element domain a set of nodes and nodal values, and defining suitable basis functions on these, for example, as piecewise-linear basis functions.

Element domains are normally shaped as triangles or squares in $\mathbb{R}^2$, tetrahedra or hexahedra in $\mathbb{R}^3$. All the nodes, edges and faces of element domains constitute the problem mesh. Defining a set of local basis functions completes the finite element space. For a rigorous definition of finite elements, and a description of different types of elements we refer to Brenner and Scott (2008).

Let $\mathcal{T}^h$ be a partition of $\Omega$ into non-overlapping elements $K$. We denote by $h$ the size of the largest element in $\mathcal{T}^h$. On the given partition $\mathcal{T}^h$ we then define the following finite element spaces, to solve the model problem:

$$X_h = \left\{ u \in C^0(\Omega) : u|_K \in P_1(K); \forall K \in \mathcal{T}^h \right\},$$

$$X_{hD} = \left\{ u \in C^0(\Omega) : u|_K \in P_1(K); u = u_D \text{ on } \Gamma_D; \forall K \in \mathcal{T}^h \right\},$$

$$X_{h0} = \left\{ u \in C^0(\Omega) : u|_K \in P_1(K); u = 0 \text{ on } \Gamma_D; \forall K \in \mathcal{T}^h \right\},$$

where $P_1(K)$ is the space of linear polynomials on $K$, and $C^0(\Omega)$ is the space of continuous functions on $\Omega$. The discretised problem, for each time step, is to find $u_h^n \in X_{hD}$, for $n = 1, 2, ..., \frac{T}{\Delta t}$ such that

$$(u_h^n, v_h) + \Delta t \left( \nabla u_h^n, \nabla v_h \right) = \Delta t \left( g_N, v_h \right)_{\Gamma_N} + \left( u_h^{n-1}, v_h \right) \quad \forall v_h \in X_{h0}. \tag{3.4}$$

We now choose the Lagrangian basis $\{\phi_1, \phi_2, ..., \phi_m\}$ of $X_h$ defined by the nodal values at the nodes $\{\boldsymbol{x}_1, \boldsymbol{x}_2, ..., \boldsymbol{x}_m\}$, namely

$$\phi_i(\boldsymbol{x}_j) = \delta_{i,j} = \begin{cases} 1, & i = j \\ 0, & i \neq j \end{cases},$$

We observe that a basis of $X_{h0}$ can be constructed by removing $\phi_i$ with $\boldsymbol{x}_i \in \Gamma_D$ from the basis of $X_h$. Let us assume that the indices of such basis functions are $1, ..., m$, and therefore $X_{h0} = \text{span}\{\phi_1, ..., \phi_m\}$. The finite-dimensional weak problem (3.4) is equivalent to: Find $u_h^n \in X_{hD}$ such that

$$(u_h^n, \phi_i) + \Delta t \left( \nabla u_h^n, \nabla \phi_i \right) = \Delta t \left( g_N, \phi_i \right)_{\Gamma_N} + \left( u_h^{n-1}, \phi_i \right) \quad \forall i = 1, ..., m. \tag{3.5}$$

Any function from $X_h$ can be presented in the form of a basis expansion. Let this basis expansion for $u_h^n$ be

$$u_h^n(\boldsymbol{x}) = \sum_{i=1}^{m} u_i^n \phi_i(\boldsymbol{x}),$$

with $u_i^n = u_h^n(\boldsymbol{x}_i)$. We define the vector of nodal values to be $\boldsymbol{u}^n = [u_1^n, ..., u_m^n]^T$. Substituting this expression into (3.5), we finally obtain a linear system which we can solve for $\boldsymbol{u}^n$:

$$(\boldsymbol{M} + \Delta t \boldsymbol{A})\boldsymbol{u}^n = \boldsymbol{M}\boldsymbol{u}^{n-1} + \Delta t \boldsymbol{g}, \tag{3.6}$$

where we have defined the following matrices and vectors:

$$\boldsymbol{A} = [\boldsymbol{a}_{ij}], \quad \boldsymbol{a}_{ij} = \int_{\Omega} \nabla\phi_i \cdot \nabla\phi_j \, \mathrm{d}x,$$

$$\boldsymbol{M} = [\boldsymbol{m}_{ij}], \quad \boldsymbol{m}_{ij} = \int_{\Omega} \phi_i \cdot \phi_j \, \mathrm{d}x,$$

$$\boldsymbol{g} = [\boldsymbol{g}_i], \quad \boldsymbol{g}_i = \int_{\Gamma_N} g_N \cdot \phi_i \, \mathrm{d}s,$$

The linear system of equations (3.6) is sparse, symmetric and positive-definite. This makes it ideal for sparse elimination methods, such as frontal solvers (Irons, 1970) that exploit the sparsity in the matrix to improve performance. Alternatively iterative methods such as the popular conjugate gradient method could be applied. We refer to Chapter 2 in (Elman et al., 2005) for a detailed discussion.

## 3.4 Mixed methods

Before considering the discretisation of the poroelasticity equations in Chapter 4 we first consider the problems of Darcy and Stokes flow. This is because many of the difficulties in solving the three-field poroelasticity problem are present when coupling the Stokes equations (elasticity of the porous mixture) with the Darcy equations (fluid flow through pores), with a modified incompressibility constraint that combines the divergence of the displacement velocity and the fluid flux. We begin with a general formulation of both the Darcy and Stokes flow equations:

$$\boldsymbol{A}(\boldsymbol{u}) + \nabla p = \boldsymbol{f} \quad \text{in } \Omega, \tag{3.7a}$$

$$\nabla \cdot \boldsymbol{u} = 0 \quad \text{in } \Omega, \tag{3.7b}$$

where $\boldsymbol{u}$ denotes the velocity vector, $p$ the pressure, $\boldsymbol{f} \in [L^2(\Omega)]^d$, with $d = 2, 3$, and $\boldsymbol{A}$ represents the two cases:

- $\boldsymbol{A}(\boldsymbol{u}) = \boldsymbol{k}^{-1}\boldsymbol{u}$, corresponding to Darcy's equation.
- $\boldsymbol{A}(\boldsymbol{u}) = -2\mu_f \nabla \cdot \boldsymbol{\epsilon}(\boldsymbol{u})$, corresponding to Stokes equation.

For simplicity we assume Dirichlet conditions on the boundary, that is, $\boldsymbol{u} = 0$ on $\partial\Omega$ for Stokes and $\boldsymbol{u} \cdot \boldsymbol{n} = 0$ on $\partial\Omega$ for Darcy. Mixed methods refer to the discretisation of different variables using different finite elements. In order to formulate our finite element method we first need the weak formulation of problem (3.7). To do this we introduce the spaces

$$W^D = \{\boldsymbol{v} \in H_{div}(\Omega) : \boldsymbol{v} \cdot \boldsymbol{n} = 0 \text{ on } \partial\Omega\},$$

$$W^S = \left\{\boldsymbol{v} \in [H^1(\Omega)]^d : \boldsymbol{v} = 0 \text{ on } \Gamma_D\right\},$$

and

$$L_0^2 = \left\{q \in L_2(\Omega) : \int_\Omega q \, \mathrm{d}x = 0\right\}.$$

We denote the product space $\mathcal{W}^X = W^X \times L_0^2$, where $X$ is chosen to be $D$ for the Darcy equations or $S$ for the Stokes equations. We also define the following

norm on $\mathcal{W}^X$:

$$\|(\boldsymbol{u},p)\|^2_{\mathcal{W}^X} = \|\boldsymbol{u}\|^2_{l,\Omega} + \|\nabla\cdot\boldsymbol{u}\|^2_{0,\Omega} + \|p\|^2_{0,\Omega},$$

with $l = 0$ for Darcy and $l = 1$ for Stokes. Let $a(\boldsymbol{u},\boldsymbol{v})$ be the bilinear form corresponding to the weak formulation of $A(\boldsymbol{u})$:

$$a(\boldsymbol{u},\boldsymbol{v}) \;=\; \left\{\begin{array}{ll} (\boldsymbol{k}^{-1}\boldsymbol{u},\boldsymbol{v}) & \text{if Darcy's equation} \\ \int_\Omega 2\mu(\epsilon(\boldsymbol{u}):\epsilon(\boldsymbol{v})) + \lambda(\nabla\cdot\boldsymbol{u})(\nabla\cdot\boldsymbol{v})\,\mathrm{d}x & \text{if Stokes equation} \end{array}\right\}.$$

Now consider the combined bilinear form

$$B[(\boldsymbol{u},p),(\boldsymbol{v},q)] = a(\boldsymbol{u},\boldsymbol{v}) - (p,\nabla\cdot\boldsymbol{v}) + (q,\nabla\cdot\boldsymbol{u}).$$

The continuous weak formulation of (3.7) is now to find $(\boldsymbol{u},p)\in\mathcal{W}^X$ such that

$$B[(\boldsymbol{u},p),(\boldsymbol{v},q)] = (\boldsymbol{f},\boldsymbol{v}) \;\; \forall(\boldsymbol{v},q)\in\mathcal{W}^X.$$

For a given finite element subspace $\mathcal{W}^X_h \in \mathcal{W}^X$, we are left with the finite dimensional problem: find $(\boldsymbol{u}_h,p_h)\in\mathcal{W}^X_h$ such that:

$$B_h[(\boldsymbol{u}_h,p_h),(\boldsymbol{v}_h,q_h)] = (\boldsymbol{f},\boldsymbol{v}_h) \;\; \forall(\boldsymbol{v}_h,q_h)\in\mathcal{W}^X_h.$$

To ensure stability and convergence of the discretisation, the discrete subspace (mixed element) has to be chosen such that the following discrete inf-sup condition, (Babuška, 1971), is fulfilled:

$$\gamma\|(\boldsymbol{u}_h,p_h)\|_{\mathcal{W}^X_h} \le \sup_{(\boldsymbol{v}_h,q_h)\in\mathcal{W}^X_h} \frac{B_h[(\boldsymbol{u}_h,p_h),(\boldsymbol{v}_h,q_h)]}{\|(\boldsymbol{v}_h,q_h)\|_{\mathcal{W}^X_h}} \;\; \forall(\boldsymbol{u}_h,p_h)\in\mathcal{W}^X_h, \qquad (3.8)$$

where $\gamma > 0$ is a constant independent of any mesh parameters. Establishing this condition ensures wellposedness of the discretisation so that the linear system arising from the fully-discrete method is non-singular and can be solved using standard methods. It is not trivial to prove (3.8) for different combinations of finite elements. This task has resulted in its own research field within Numerical Analysis, and countless papers have been published on this topic. In table 3.1 we have documented some popular standard finite element pairs for solving the Stokes and Darcy equations, and outlined whether these satisfy (3.8), thereby yielding a stable and optimally converging method, or not. Note that many other possible discretisations exist.

| Mixed element | Stokes | Darcy |
|---|---|---|
| $P1 - P1$ | ✗ | ✗ |
| $P2 - P1$ | ✓ | ✗ |
| $P1 - P1 + stab$ | ✓ | ✓ |
| $P1 - P0$ | ✗ | ✗ |
| $RT - P0$ | ✗ | ✓ |
| $P1 - P0 + stab$ | ✓ | ✓ |

Table 3.1: Possible finite element combinations for Stokes and Darcy flow, showing whether a particular choice of elements is stable and optimally converging or not.

The naive choice of piecewise linear finite elements for both the velocities and the pressure, denoted by $(P1 - P1)$, or piecewise linear finite elements for the velocities and piecewise constants for the pressure, $(P1 - P0)$, result in an ill posed discretisation (Burman and Hansbo, 2007). Intuitively, this is because the velocity space is not rich enough to constrain the pressures, thus resulting in spurious pressure oscillations. A detailed explanation of this along with some worked examples can be found in Elman et al. (2005). The Taylor-Hood element, $(P2 - P1)$ - piecewise quadratic for the velocities and piecewise linear for the pressure, is a commonly used element for the Stokes equations. However for the Darcy equations this element does not convergence at the right order and fails to converge for the divergence of the velocities (Burman and Hansbo, 2007). The Raviart-Thomas element, $(RT - P0)$, first proposed in Raviart and Thomas (1977) is a divergence free element, often used to solve the Darcy equations. However this element is not able to control $H^1$ velocities, and therefore can not be used to solve the Stokes equations. When the finite element discretisation is based on a discrete subspace that does not satisfy the discrete inf-sup condition (3.8), a procedure aiming at stabilising the discrete system may be accomplished. The philosophy of stabilised methods is to strengthen formulations by adding an extra term, often to the mass conservation equation, so that discrete approximations, which would otherwise be unstable, become stable and convergent (Masud and Hughes, 2002). Numerous stabilisation techniques exist. To stabilise the equal order piecewise linear pair, a polynomial pressure projection has been proposed in Bochev and Dohrmann (2006) that results in a stable element for both the Stokes and Darcy equations, $(P1 - P1 + stab)$. Also, a pressure jump stabilisation, $(P1 - P0 + stab)$, that uses a piecewise constant pressure approximation and is stable and optimally converging for both the Stokes and Darcy equation has been analysed in Burman and Hansbo (2007). This is the stabilisation we will modify to solve the poroelastic equations.

## 3.5 Poroelastic finite element discretisations

### 3.5.1 Linear discretisations

The linear poroelastic equations are often solved in a reduced displacement and pressure formulation, from which the fluid flux can then be recovered (Murad and Loula, 1994; White and Borja, 2008). In Murad and Loula (1994) the stability and convergence of this reduced displacement pressure $(\boldsymbol{u}/p)$ formulation has been analysed. They were also able to show error bounds for inf-sup stable combinations of finite element spaces (e.g. Taylor-Hood elements). In this work we will keep the fluid flux variable resulting in a three-field, displacement, fluid flux, and pressure formulation. Keeping the fluid flux as a primary variable has the following advantages:

i It allows for greater accuracy in the fluid velocity field. This can be of interest whenever a poroelastic model is coupled with an advection diffusion equation, e.g. to account for gas exchange, thermal effects, contaminant transport or the transport of nutrients or drugs within a porous tissue (Khaled and Vafai, 2003).

ii Physically meaningful boundary conditions can be applied at the interface when modelling the interaction between a fluid and a poroelastic structure (Badia et al., 2009).

iii It allows for an easy extension of the fluid model from a Darcy to a Brinkman flow model, for which there are numerous applications in modelling biological tissues (Khaled and Vafai, 2003).

iv It reduces the order of the spatial derivative of the pressure, allowing for a discontinuous pressure approximation without any additional penalty terms.

v It avoids the calculation of the fluid flux in post-processing.

Error estimates have been proven in Phillips and Wheeler (2007a,b) for solving the three-field formulation problem using continuous piecewise linear approximations for displacements and mixed low-order Raviart Thomas elements for the fluid flux and pressure variables. However this method was found to be susceptible to spurious pressure oscillations (Phillips and Wheeler, 2009). To overcome these pressure oscillations, Li and Li (2012) analysed a discontinuous three-field method with moderate success, and Yi (2013) analysed a non-conforming three-field method. However no implementation of these methods in 3D has yet been presented. We hypothesize that this is due to the complexity of these non-standard elements used, making it very difficult to include them in existing finite element codes.

In addition to these monolithic approaches there has been considerable work on operator splitting (iterative) approaches where the poroelastic equations are separated into a fluid problem and elasticity problem. Each of these subsystems is then solved in a staggered fashion, and the solution is passed between the solvers. (Feng and He, 2010; Kim et al., 2011). For example the elasticity problem is solved, and the resulting deformation passed to the fluid solver for an improved solution of the fluid flux and pressure. The pressure is then passed back to the elasticity solver for an improved estimate of the deformation. This is repeated until convergence is achieved. The degree of coupling of the problem affects the stability and accuracy of the numerical solution (Wheeler and Gai, 2007). Although these methods are often able to take advantage of existing elasticity and fluid finite element software, and result in solving a smaller system of equations, these schemes are often only conditionally stable, and very small time steps may be required. The advantage of a monolithic approach is that the linear solver must solve simultaneously for the fluid variables and deformation variables, which ensures that a solution is always achieved, and any size time step can be used. Not having to deal with additional convergence tolerances and restrictions on the time step can significantly simplify the use of the method and improve the computational performance of problems that tightly couple the pore pressure with the deformation.

### 3.5.2 Discretisations valid in large deformations

We will now give a brief overview of different approaches for solving the poroelastic equations valid in large deformations. There has been some work on operator splitting (iterative) approaches (Chapelle et al., 2010). Again, such approaches are often only conditionally stable. Some notable quasi-static incompressible large deformation monolithic approaches include a mixed-penalty formulation, and a mixed solid velocity-pressure formulation, both outlined in Almeida and Spilker (1998), the solid velocity-pressure formulation is similar to the commonly used reduced ($\boldsymbol{u}/p$) formulation (Ateshian et al., 2010). These two-field formulations require a stable mixed element pair such as the popular Taylor-Hood element to satisfy the LBB inf-sup stability requirement. The key difficulty, however, that these elements cannot escape is that jumps in material coefficients may introduce large solution gradients across the interface, requiring severe mesh refinement. This is because a continuous pressure element is used, which is unable to reliably capture jumps in the pressure solution (White and Borja, 2008). In Levenston et al. (1998) a three-field (displacement, fluid flux, pressure) formulation has been outlined, however this method uses a low-order mixed finite element approximation without any stabilisation and therefore is not inf-sup stable. A dynamic three-field finite element using a continuous pressure approximation has been implemented in Vuong et al. (2015).

# Chapter 4

# Analysis of a stabilised finite element method for linear poroelasticity

The contents of this chapter closely follows the theoretical sections presented in the joint publication: L. Berger, R. Bordas, D. Kay, and S. Tavener; Stabilized low-order finite element approximation for linear three-field poroelasticity *SIAM Journal on Scientific Computing* 2015. D. Kay had the initial idea of applying a pressure jump stabilisation to three-field poroelasticity. L. Berger developed all the proofs, with guidance from D. Kay and R. Bordas, and wrote the original draft of the paper. S. Tavener assisted in simplifying the proofs, and improving the quality of the writing and the structure of the paper, along with the other authors.

## 4.1 Introduction

In this chapter we develop a stabilised, low-order, mixed finite element method for poroelastic models of biological tissues and restrict our attention to the fully saturated, incompressible, small deformation case. Our mixed scheme uses the lowest possible approximation order: piecewise constant approximation for the pressure and piecewise linear continuous elements for the displacement and fluid flux.

To ensure stability, a mixed finite element method must satisfy the Ladyzhenskaya-Babuska-Brezzi (LBB) condition. In this work we use a local pressure jump stabilisation method pioneered by Burman and Hansbo (2007) for the study of Stokes and Darcy flows that are coupled via an interface. This approach provides the natural $H^1$ stability for the displacements and $H_{div}$ stability for the fluid flux. In this Chapter we prove the stability of the mixed finite element method for poroelasticity using results and steps taken from Burman and Hansbo (2007). We also show that the naive approach of using the stabilisation of the pressure, as is done for the Darcy and Stokes equations in Burman and Hansbo (2007), results in an approximation that does not converge at an optimal rate. Stabilisation using the time derivative of pressure in the stabilisation term is shown to be crucial for stability and optimal convergence with refinement and counterexamples are provided in Section 6.5.

In section 4.2 we formulate the model and its continuous weak formulation and construct a fully-discrete approximation. In section 4.3 we will introduce some norms and inequalities. We prove existence and uniqueness of solutions to this discrete model at each time step in section 4.4, provide an energy estimate over time in section 4.5, and derive an optimal order a-priori error estimate in section 4.6.

## 4.2 The poroelastic model

### 4.2.1 Governing equations

Following Phillips and Wheeler (2007a) and Showalter (2000), we recall the governing equations (2.40) for a fully saturated, incompressible poroelastic model

$$-(\lambda+\mu)\nabla(\nabla\cdot\boldsymbol{u})-\mu\nabla^2\boldsymbol{u}+\nabla p=\boldsymbol{f} \quad \text{in } \Omega, \tag{4.1a}$$

$$\boldsymbol{k}^{-1}\boldsymbol{z}+\nabla p=\boldsymbol{b} \quad \text{in } \Omega, \tag{4.1b}$$

$$\nabla\cdot(\boldsymbol{u}_t+\boldsymbol{z})=g \quad \text{in } \Omega, \tag{4.1c}$$

$$\boldsymbol{u}=\boldsymbol{u}_D \quad \text{on } \Gamma_D, \tag{4.1d}$$

$$\boldsymbol{\sigma}\boldsymbol{n}=\boldsymbol{t}_N \quad \text{on } \Gamma_N, \tag{4.1e}$$

$$\boldsymbol{z}\cdot\boldsymbol{n}=q_D \quad \text{on } \Gamma_F, \tag{4.1f}$$

$$p=p_D \quad \text{on } \Gamma_P, \tag{4.1g}$$

$$\boldsymbol{u}(\cdot,0)=\boldsymbol{u}^0 \quad \text{in } \Omega. \tag{4.1h}$$

**Remark 4.2.1.** *Since the above resulting system of equations is linear, for ease of presentation, we will assume all Dirichlet boundary conditions are homogeneous, ie.,* $\boldsymbol{u}_D=\boldsymbol{0}, q_D=0, p_D=0$.

### 4.2.2 Weak formulation

We define the following spaces for displacement, fluid flux and pressure respectively,

$$\begin{aligned}
\boldsymbol{W}^E(\Omega) &= \{\boldsymbol{u}\in(H^1(\Omega))^d : \boldsymbol{u}=\boldsymbol{0} \text{ on } \Gamma_D\},\\
\boldsymbol{W}^D(\Omega) &= \{\boldsymbol{z}\in H_{div}(\Omega) : \boldsymbol{z}\cdot\boldsymbol{n}=0 \text{ on } \Gamma_F\},\\
\mathcal{L}(\Omega) &= \left\{\begin{array}{ll} L^2(\Omega) & \text{if } \Gamma_N\cup\Gamma_P\neq\emptyset \\ L^2_0(\Omega) & \text{if } \Gamma_N\cup\Gamma_P=\emptyset, \end{array}\right\},
\end{aligned}$$

where $L^2_0(\Omega) = \left\{q \in L^2(\Omega) : \int_\Omega q \, \mathrm{d}x = 0\right\}$, which we combine to construct the mixed solution space

$$\mathcal{W}^X = \left\{\boldsymbol{W}^E(\Omega) \times \boldsymbol{W}^D(\Omega) \times \mathcal{L}(\Omega)\right\}.$$

We define the bilinear form

$$a(\boldsymbol{u}, \boldsymbol{v}) = \int_\Omega 2\mu(\epsilon(\boldsymbol{u}) : \epsilon(\boldsymbol{v})) + \lambda(\nabla \cdot \boldsymbol{u})(\nabla \cdot \boldsymbol{v}) \, \mathrm{d}x,$$

for $\boldsymbol{u}, \boldsymbol{v} \in \boldsymbol{W}^E(\Omega)$. This bilinear form is continuous such that

$$a(\boldsymbol{u}, \boldsymbol{v}) \leq C_c \|\boldsymbol{u}\|_{1,\Omega} \|\boldsymbol{v}\|_{1,\Omega} \; \forall \boldsymbol{u}, \boldsymbol{v} \in (H^1(\Omega))^d. \tag{4.2}$$

Using Korn's inequality (Brenner and Scott, 2008; Ciarlet, 1978), and $\int_\Omega \lambda(\nabla \cdot \boldsymbol{v})(\nabla \cdot \boldsymbol{v}) \geq 0$ we have

$$\|\boldsymbol{v}\|^2_{a,\Omega} = a(\boldsymbol{v}, \boldsymbol{v}) \geq 2\mu \|\epsilon(\boldsymbol{v})\|^2_{0,\Omega} \geq C_k \|\boldsymbol{v}_h\|^2_{1,\Omega} \;\; \forall \boldsymbol{v} \in \boldsymbol{W}^E(\Omega_t). \tag{4.3}$$

Since $\boldsymbol{k}$ is assumed to be a symmetric and strictly positive definite tensor, there exists eigenfunctions $\lambda_{min}, \lambda_{max} > 0$ such that $\forall \boldsymbol{x} \in \Omega, \; \lambda_{min} \|\boldsymbol{\eta}\|^2_{0,\Omega} \leq \boldsymbol{\eta}^t \boldsymbol{k}(\boldsymbol{x}) \boldsymbol{\eta} \leq \lambda_{max} \|\boldsymbol{\eta}\|^2_{0,\Omega} \;\; \forall \boldsymbol{\eta} \in \mathbb{R}^d$, and

$$\lambda^{-1}_{min} \|\boldsymbol{w}\|^2_{0,\Omega} \geq (\boldsymbol{k}^{\boldsymbol{-1}} \boldsymbol{w}, \boldsymbol{w}) \geq \lambda^{-1}_{max} \|\boldsymbol{w}\|^2_{0,\Omega} \; \forall \boldsymbol{w} \in \boldsymbol{W}^D(\Omega_t). \tag{4.4}$$

The continuous weak problem is: Find $\boldsymbol{u}(x,t) \in \boldsymbol{W}^E(\Omega)$, $\boldsymbol{z}(x,t) \in \boldsymbol{W}^D(\Omega)$, and $p(x,t) \in \mathcal{L}(\Omega)$ for any time $t \in (0,T]$ such that:

$$a(\boldsymbol{u},\boldsymbol{v}) - (p, \nabla\cdot\boldsymbol{v}) = (\boldsymbol{f},\boldsymbol{v}) + (\boldsymbol{t}_N,\boldsymbol{v})_{\Gamma_N}\ \forall \boldsymbol{v} \in \boldsymbol{W}^E(\Omega_t), \tag{4.5a}$$

$$(\boldsymbol{k^{-1}}\boldsymbol{z},\boldsymbol{w}) - (p,\nabla\cdot\boldsymbol{w}) = (\boldsymbol{b},\boldsymbol{w})\ \forall \boldsymbol{w} \in \boldsymbol{W}^D(\Omega_t), \tag{4.5b}$$

$$(\nabla\cdot\boldsymbol{u}_t, q) + (\nabla\cdot\boldsymbol{z}, q) = (g,q)\ \forall q \in \mathcal{L}(\Omega_t). \tag{4.5c}$$

We will assume the following regularity requirements on the data,

$$\begin{aligned} &\boldsymbol{f} \in C^1((0,T];(H^{-1}(\Omega))^d), \qquad &&\boldsymbol{t}_N \in C^1((0,T];H^{-1/2}(\Gamma_N)),\\ &\boldsymbol{b} \in C^1((0,T];H^{-1}_{div}(\Omega)), &&g \in C^0((0,T];(L^2(\Omega))^d). \end{aligned} \tag{4.6}$$

For the initial conditions we require that $\boldsymbol{u}^0 \in (H^1(\Omega))^d$. The well-posedness of the continuous two-field formulation has been proven by Showalter (2000). Lipnikov (2002) proves well-posedness for the continuous three-field formulation (6.2). In this work we also establish the well-posedness of (6.2) as a result of the energy estimates proven in section 4.5, see remark 4.5.1.

### 4.2.3 Fully-discrete model

We define the following finite element spaces,

$$\begin{aligned} \boldsymbol{W}_h^E &= \left\{\boldsymbol{u}_h \in C^0(\Omega) : \boldsymbol{u}_h|_K \in P_1(K)\ \forall K \in \mathcal{T}^h, \boldsymbol{u}_h = \boldsymbol{0} \text{ on } \Gamma_D\right\},\\ \boldsymbol{W}_h^D &= \left\{\boldsymbol{z}_h \in C^0(\Omega) : \boldsymbol{z}_h|_K \in P_1(K)\ \forall K \in \mathcal{T}^h, \boldsymbol{z}_h\cdot\boldsymbol{n} = 0 \text{ on } \Gamma_F\right\},\\ Q_h &= \begin{cases} \left\{p_h : p_h|_K \in P_0(K)\ \forall K \in \mathcal{T}^h\right\} & \text{if } \Gamma_N \cup \Gamma_P \neq \emptyset\\ \left\{p_h : p_h|_K \in P_0(K), \int_\Omega p_h = 0\ \forall K \in \mathcal{T}^h\right\} & \text{if } \Gamma_N \cup \Gamma_P = \emptyset \end{cases}, \end{aligned}$$

where $P_0(K)$ and $P_1(K)$ are respectively the spaces of constant and linear polynomials on $K$. We partition $[0,T]$ into $N$ evenly spaced non-overlapping regions $(t_{n-1}, t_n]$, $n = 1, 2, \ldots, N$ , where $t_n - t_{n-1} = \Delta t$. For any sufficiently smooth function $v(x,t)$ we define $v^n(x) = v(x, t_n)$ and the discrete time derivative by $v^n_{\Delta t} = \frac{v^n - v^{n-1}}{\Delta t}$.

The fully-discrete weak problem is: For $n = 1, 2, \ldots, N$, find $\boldsymbol{u}^n_h \in \boldsymbol{W}^E_h$, $\boldsymbol{z}^n_h \in \boldsymbol{W}^D_h$ and $p^n_h \in Q_h$ such that

$$a(\boldsymbol{u}^n_h, \boldsymbol{v}_h) - (p^n_h, \nabla \cdot \boldsymbol{v}_h) = (\boldsymbol{f}^n, \boldsymbol{v}_h) + (\boldsymbol{t}_N, \boldsymbol{v}_h)_{\Gamma_N} \; \forall \boldsymbol{v}_h \in \boldsymbol{W}^E_h, \tag{4.7a}$$

$$(\boldsymbol{k}^{-1} \boldsymbol{z}^n_h, \boldsymbol{w}_h) - (p^n_h, \nabla \cdot \boldsymbol{w}_h) = (\boldsymbol{b}^n, \boldsymbol{w}_h) \; \forall \boldsymbol{w}_h \in \boldsymbol{W}^D_h, \tag{4.7b}$$

$$(\nabla \cdot \boldsymbol{u}^n_{\Delta t,h}, q_h) + (\nabla \cdot \boldsymbol{z}^n_h, q_h) + J\left(p^n_{\Delta t,h}, q_h\right) = (g^n, q_h) \; \forall q_h \in Q_h. \tag{4.7c}$$

The stabilisation term is

$$J(p,q) = \delta \sum_K \int_{\partial K \backslash \partial \Omega} h_{\partial K} [p][q] \, \mathrm{d}s. \tag{4.8}$$

Here $\delta$ is a stabilisation parameter that is independent of $h$ and $\Delta t$. Here $h_{\partial K}$ denotes the size (diameter) of an element edge in 2D or face in 3D, and $[\cdot]$ is the jump across an edge or face (taken on the interior edges only). We will see in the numerical results, Chapter 5 that the convergence is not sensitive to $\delta$. The set of all elements is denoted by $K$, $h_{\partial K}$ denotes the size of an element edge in 2D or face in 3D, and $[\cdot]$ is the jump across an edge. The jump in pressure $[p_h]$ across an element or face $E$ adjoining elements $T$ and $S$ is defined such that

$$(p_h|_T - p_h|_S)\boldsymbol{n}_{E,T} = (p_h|_S - p_h|_T)\boldsymbol{n}_{E,S}.$$

Here $\boldsymbol{n}_{E,T}$ is the outward normal from element $T$, with respect to edge $E$, $\boldsymbol{n}_{E,S}$ is the corresponding inward facing normal, and $p_h|_T$ and $p_h|_S$ denote the pressure in element $T$ and $S$, respectively.

We also assume

$$a(\boldsymbol{u}_h^0, \boldsymbol{v}_h) = a(\boldsymbol{u}^0, \boldsymbol{v}_h) \ \forall \boldsymbol{v}_h \in \boldsymbol{W}_h^E, \tag{4.9a}$$

$$J(p_h^0, q_h) = J(p^0, q_h) \ \forall q_h \in Q_h, \tag{4.9b}$$

where $p^0 \in \mathcal{L}(\Omega)$.

## 4.3 Norms and inequalities

In this section we will introduce some norms and inequalities required for the remainder of this chapter. Throughout this work, we will let $C$ denote a generic positive constant, whose value may change from instance to instance, but is independent of any mesh parameters.

### 4.3.1 Useful inequalities

Detailed derivations of the following four inequalities can be found in Brenner and Scott (2008). If $f, g \in L^2(\Omega)$ then by the **Cauchy-Schwarz** inequality we have

$$\int_\Omega |f(x)g(x)|dx \leq \|f\|_{0,\Omega}\|g\|_{0,\Omega}.$$

From the **triangle inequality** we have

$$\|f + g\|_{0,\Omega} \leq \|f\|_{0,\Omega} + \|g\|_{0,\Omega}.$$

For any real numbers $a$ and $b$, by **Young's inequality**,

$$ab \leq \frac{\epsilon}{2}a^2 + \frac{1}{2\epsilon}b^2 \ \ \forall \epsilon > 0.$$

This inequality is sometimes referred to as the arithmetic-geometric mean inequality.

Next, assuming $\int_{\Gamma_D} \mathrm{d}s \neq 0$ and $C_p > 0$, the **Poincaré inequality**, also known as Poincaré-Friedrich's inequality is given by

$$\|\boldsymbol{u}\|_{0,\Omega} \leq C_p \|\nabla \boldsymbol{u}\|_{0,\Omega} \ \ \forall \boldsymbol{u} \in \boldsymbol{W}^E(\Omega).$$

### 4.3.2 Properties of the J-norm

The stabilisation term gives rise to the semi-norm

$$|q|_{J,\Omega} = J(q,q)^{1/2}.$$

Using the scaling argument, also used in Burman and Hansbo (2007),

$$\left\|h^{1/2} p_h\right\|_{0,\partial K} \leq c_z \|p_h\|_{0,K} \ \forall p_h \in Q_h. \tag{4.10}$$

Cauchy-Schwarz and the triangle inequality the following bounds for the stabilisation term hold.

$$|p_h|_{J,\Omega} \leq C\|p_h\|_{0,\Omega} \text{ and } J(p_h, q_h) \leq |p_h|_{J,\Omega}|q_h|_{J,\Omega}, \ \forall p_h, q_h \in Q_h. \tag{4.11}$$

Furthermore, for any $q \in H^1(\Omega)$,

$$J(p,q) = 0, \ \forall p \in \mathcal{L}(\Omega), \tag{4.12}$$

which forms the corner stone of the method's error estimate and was originally proposed in Silvester and Kechkar (1990). Also see Lemma 1.23 in Di Pietro and

Ern (2011).

### 4.3.3 Approximation results

We now give some approximation results that will be useful later. Let $\pi_h^1 : H^1(\Omega) \to \boldsymbol{W}_h^E$ and $\pi_h^0 : L^2(\Omega) \to Q_h$ be Clément projections (interpolation operators), see Ciarlet (1978).

**Lemma 4.3.1.** *For all $v \in (H^2(\Omega))^d$ and $q \in H^1(\Omega)$ the interpolation operators satisfy: For $s = 0, 1$*

$$||v - \pi_h^1 v||_{s,\Omega} \leq Ch^{2-s}\|v\|_{2,\Omega}, \tag{4.13}$$

$$\left\|q - \pi_h^0 q\right\|_{0,\Omega} \leq Ch\|q\|_{1,\Omega}, \tag{4.14}$$

$$|q - \pi_h^0 q|_{J,\Omega} \leq Ch\|q\|_{1,\Omega}. \tag{4.15}$$

*Proof.* The first two results are standard Brenner and Scott (2008). The final result is obtained by using the element error estimate provided in Verfürth (1998) and then summing over all elements. □

Due to the surjectivity of the divergence operator, for every $p \in L^2(\Omega)$ there exists a function $\boldsymbol{v}_p \in (H_0^1(\Omega))^d$ such that $\nabla\cdot\boldsymbol{v}_p = -p$ and $\|\boldsymbol{v}_p\|_{1,\Omega} \leq c\|p\|_{0,\Omega}$. This last inequality can be shown to hold by considering the famous inf-sup condition related to the continous Stokes problem (Brenner and Scott, 2008; Brezzi and Fortin, 1991). We assume that the projection, $\pi_h^1\boldsymbol{v}_p$, is stable such that

$$\left\|\pi_h^1\boldsymbol{v}_p\right\|_{1,\Omega} \leq \hat{c}\|p\|_{0,\Omega}. \tag{4.16}$$

Furthermore, for any element $K \in \mathcal{T}^h$

$$||\boldsymbol{v}_p - \pi_h^1\boldsymbol{v}_p||_{L^2(K)} \leq Ch||\boldsymbol{v}_p||_{H^1(\omega_K)}, \tag{4.17}$$

where $\omega_K$ is a domain made of the elements in $\mathcal{T}^h$ neighbouring $K$. For more details about the properties of this projection we refer to section 4.8 in Brenner and Scott (2008). This projection will allow us to obtain stability of the pressure and avoid spurious pressure oscillations. The discrepancy between the projection and its continuous counterpart will eventually be made up by the stabilisation term, shown in section 4.4. Combining the above with the trace inequality, see lemma 3.1 in Verfürth (1998),

$$\big\|(\boldsymbol{v}_p - \pi_h^1 \boldsymbol{v}_p)\cdot \boldsymbol{n}\big\|_{0,\partial K}^2 \leq C\big\|\boldsymbol{v}_p - \pi_h^1 \boldsymbol{v}_p\big\|_{0,K}(h^{-1}\big\|\boldsymbol{v}_p - \pi_h^1 \boldsymbol{v}_p\big\|_{0,K} + \big\|\boldsymbol{v}_p - \pi_h^1 \boldsymbol{v}_p\big\|_{1,K}), \tag{4.18}$$

we obtain

$$\big\|(\boldsymbol{v}_p - \pi_h^1 \boldsymbol{v}_p)\cdot \boldsymbol{n})\big\|_{0,\partial K}^2 \leq Ch||\boldsymbol{v}_p||^2_{H^1(\omega_K)}. \tag{4.19}$$

Taking into account $\|\boldsymbol{v}_p\|_{1,\Omega} \leq c\|p\|_{0,\Omega}$, we may write

$$\sum_K \int_{\partial K} h^{-1}|(\boldsymbol{v}_p - \pi_h^1 \boldsymbol{v}_p)\cdot \boldsymbol{n}|^2 \, ds \leq c_t\|p\|_{0,\Omega}^2. \tag{4.20}$$

We also have the following approximation for the time-discretisation error: For all $v \in H^2(0,T;(L^2(\Omega))^d)$

$$\sum_{n=1}^{N} \Delta t \left\| v_{\Delta t}^n - \frac{\partial v}{\partial t}(t^n,\cdot)\right\|_{0,\Omega}^2 \leq \Delta t^2 \int_0^T \|v_{tt}\|_{0,\Omega}^2 \mathrm{d}s. \tag{4.21}$$

See (Brenner and Scott, 2008; Thomée, 2006) for details.

### 4.3.4 Triple-norms

We will now define some triple-norms that are designed to get the required cancellation of the divergence terms and will allow us to obtain control in Step 1 (4.27) in the proof of Theorem 4.4.1. For all $[v,w,q] \in \big[(H^1(\Omega))^d \times H_{div}(\Omega) \times L^2(\Omega)\big]$

we define the norm

$$|||[v,w,q]|||_A^2 = \|v\|_{1,\Omega}^2 + \Delta t^2 \|\nabla \cdot w\|_{0,\Omega}^2 + \Delta t \|w\|_{0,\Omega}^2 + \|q\|_{0,\Omega}^2 + |q|_{J,\Omega}^2. \quad (4.22)$$

The above triple-norm has also been chosen to satisfy the continuity property (4.25). For all $[v,w,q] \in \left[L^\infty(0,T;(H^1(\Omega))^d) \times L^2(0,T;H_{div}(\Omega)) \times L^2(0,T;L^2(\Omega))\right]$ we define the norm

$$|||[v,w,q]|||_B^2 = ||v||_{L^\infty(H^1)}^2 + ||w||_{L^2(L^2)}^2 + ||q||_{L^2(L^2)}^2. \quad (4.23)$$

## 4.4 Existence and uniqueness of solutions to the fully-discrete model

Well-posedness of the unstabilised fully-discretised system (4.7) (i.e., for $\delta = 0$), with the use of a low order Raviart-Thomas approximation for the fluid velocity is shown by Phillips and Wheeler (2007b) for $c_0 > 0$, and by Lipnikov (2002) for $c_0 \geq 0$. Although as the permeability tends to zero and the porous mixture becomes impermeable, the three-field linear poroelasticity tends to a mixed linear elasticity problem (Haga et al., 2012). Hence, in this case this element becomes unstable, as expected since the elasticity $P1-P0$ approximation is known to be unstable. Our method is stable for both the Darcy problem (as the elasticity coefficients tend to infinity) and the mixed linear elasticity problem (as the permeability tends to zero), and is therefore stable for all permeabilities and elasticity coefficients.

Combining the fully-discrete equations (4.7a), (4.7b) and (4.7c), after first multiplying (4.7b) and (4.7c) by $\Delta t$, gives the equivalent problem;

For $n = 1, 2, \ldots, n$, find $(\boldsymbol{u}_h, \boldsymbol{z}_h, p_h)$ such that

$$
\begin{aligned}
&B_h^n[(\boldsymbol{u}_h, \boldsymbol{z}_h, p_h), (\boldsymbol{v}_h, \boldsymbol{w}_h, q_h)] \\
&\quad = (\boldsymbol{f}^n, \boldsymbol{v}_h) + (\boldsymbol{t}_N, \boldsymbol{v}_h)_{\Gamma_N} + \Delta t(\boldsymbol{b}^n, \boldsymbol{w}_h) + \Delta t(g^n, q_h) \\
&\qquad + (\nabla \cdot \boldsymbol{u}_h^{n-1}, q_h) + J(p_h^{n-1}, q_h) \quad \forall (\boldsymbol{v}_h, \boldsymbol{w}_h, q_h) \in \mathcal{W}_h^X,
\end{aligned}
$$

where

$$
\begin{aligned}
&B_h^n[(\boldsymbol{u}_h, \boldsymbol{z}_h, p_h), (\boldsymbol{v}_h, \boldsymbol{w}_h, q_h)] \\
&\quad = a(\boldsymbol{u}_h^n, \boldsymbol{v}_h) + \Delta t(\boldsymbol{k}^{-1} \boldsymbol{z}_h^n, \boldsymbol{w}_h) - (p_h^n, \nabla \cdot \boldsymbol{v}_h) - \Delta t(p_h^n, \nabla \cdot \boldsymbol{w}_h) \\
&\qquad + (\nabla \cdot \boldsymbol{u}_h^n, q_h) + \Delta t(\nabla \cdot \boldsymbol{z}_h^n, q_h) + J(p_h^n, q_h). \qquad (4.24)
\end{aligned}
$$

The linear form satisfies the following continuity property

$$
|B_h^n[(\boldsymbol{u}_h, \boldsymbol{z}_h, p_h), (\boldsymbol{v}_h, \boldsymbol{w}_h, q_h)]| \le C \, |||(\boldsymbol{u}_h^n, \boldsymbol{z}_h^n, p_h^n)|||_A \, |||(\boldsymbol{v}_h, \boldsymbol{w}_h, q_h)|||_A \,. \qquad (4.25)
$$

We apply Babuška's theory (Babuška, 1971) to show well-posedness (existence and uniqueness) of this discretised system at a particular time step. This requires us to prove a discrete inf-sup type result (Theorem 4.4.1) for the combined bilinear form (4.24).

**Theorem 4.4.1.** *Let $\gamma > 0$ be a constant independent of any mesh parameters. Then the finite element formulation (4.7) satisfies the following discrete inf-sup condition*

$$
\gamma \, |||(\boldsymbol{u}_h^n, \boldsymbol{z}_h^n, p_h^n)|||_A \le \sup_{(\boldsymbol{v}_h, \boldsymbol{w}_h, q_h) \in \mathcal{V}_h^X} \frac{B_h^n[(\boldsymbol{u}_h, \boldsymbol{z}_h, p_h), (\boldsymbol{v}_h, \boldsymbol{w}_h, q_h)]}{|||(\boldsymbol{v}_h, \boldsymbol{w}_h, q_h)|||_A} \;\; \forall (\boldsymbol{u}_h, \boldsymbol{z}_h, p_h) \in \mathcal{W}_h^X. \qquad (4.26)
$$

*Hence, given a solution at the previous time step the linear system arising from the fully-discrete method for the subsequent time step is non-singular.*

The following proof follows ideas presented by Burman and Hansbo (2007).

*Proof.*

*Step 1, bounding $\|\boldsymbol{u}_h^n\|_{1,\Omega}$, $\Delta t^{1/2}\|\boldsymbol{z}_h^n\|_{0,\Omega}$, and $|p_h^n|_{J,\Omega}$.*

Choose $(\boldsymbol{v}_h, \boldsymbol{w}_h, q_h) = (\beta \boldsymbol{u}_h^n, \beta \boldsymbol{z}_h^n, \beta p_h^n)$, then using (4.3) and (4.4), we obtain,

$$\begin{aligned} B_h^n[(\boldsymbol{u}_h, \boldsymbol{z}_h, p_h), (\beta \boldsymbol{u}_h^n, \beta \boldsymbol{z}_h^n, \beta p_h^n)] &= a(\boldsymbol{u}_h^n, \beta \boldsymbol{u}_h^n) + \Delta t(\boldsymbol{k^{-1}} \boldsymbol{z}_h^n, \beta \boldsymbol{z}_h^n) + J(p_h^n, \beta p_h^n) \\ &\geq \beta C_k \|\boldsymbol{u}_h^n\|_{1,\Omega}^2 + \beta \lambda_{max}^{-1} \Delta t \|\boldsymbol{z}_h^n\|_{0,\Omega}^2 + \beta |p_h^n|_{J,\Omega}^2. \quad (4.27) \end{aligned}$$

By being able to choose $\beta$ arbitrarily large, this step will enable us to regain control in step 4, and thus using the stabilisation, $\beta |p_h^n|_{J,\Omega}^2$, we can control the pressure.

*Step 2, bounding $\|p_h^n\|_{0,\Omega}$.*

Choose $(\boldsymbol{v}_h, \boldsymbol{w}_h, q_h) = (\pi_h^1 \boldsymbol{v}_{p_h^n}, \boldsymbol{0}, 0)$ and add $0 = \|p_h^n\|_{0,\Omega}^2 + (p_h^n, \nabla \cdot \boldsymbol{v}_{p_h^n})$ to obtain

$$B_h^n[(\boldsymbol{u}_h, \boldsymbol{z}_h, p_h), (\pi_h^1 \boldsymbol{v}_{p_h^n}, \boldsymbol{0}, 0)] = a(\boldsymbol{u}_h^n, \pi_h^1 \boldsymbol{v}_{p_h^n}) + \|p_h^n\|_{0,\Omega}^2 + (p_h^n, \nabla \cdot (\boldsymbol{v}_{p_h^n} - \pi_h^1 \boldsymbol{v}_{p_h^n})). \quad (4.28)$$

Focusing on the third term in (4.28) only, we apply the divergence theorem and split the integral over local elements to get

$$\begin{aligned} (p_h^n, \nabla \cdot (\boldsymbol{v}_{p_h^n} - \pi_h^1 \boldsymbol{v}_{p_h^n})) &= \sum_K \int_{\partial K} p_h^n (\boldsymbol{v}_{p_h^n} - \pi_h^1 \boldsymbol{v}_{p_h^n}) \cdot \boldsymbol{n} \, ds \\ &= \sum_K \frac{1}{2} \int_{\partial K} [p_h^n] (\boldsymbol{v}_{p_h^n} - \pi_h^1 \boldsymbol{v}_{p_h^n}) \cdot \boldsymbol{n} \, ds. \end{aligned}$$

We thus have

$$B_h^n[(\boldsymbol{u}_h, \boldsymbol{z}_h, p_h), (\pi_h^1 \boldsymbol{v}_{p_h^n}, \boldsymbol{0}, 0)] = \|p_h^n\|_{0,\Omega}^2 + a(\boldsymbol{u}_h^n, \pi_h^1 \boldsymbol{v}_{p_h^n}) + \sum_K \frac{1}{2} \int_{\partial K} [p_h^n](\boldsymbol{v}_{p_h^n} - \pi_h^1 \boldsymbol{v}_{p_h^n}) \cdot \boldsymbol{n} \, ds.$$

Now first applying the Cauchy-Schwarz inequality and (4.2) on the right hand side to get

$$B_h^n[(\boldsymbol{u}_h, \boldsymbol{z}_h, p_h), (\pi_h^1 \boldsymbol{v}_{p_h^n}, \boldsymbol{0}, 0)] \geq \|p_h^n\|_{0,\Omega}^2 - C_c \|\boldsymbol{u}_h^n\|_{1,\Omega} \left\| \pi_h^1 \boldsymbol{v}_{p_h^n} \right\|_{1,\Omega} - \sum_K \frac{1}{2} \left( \int_{\partial K} \left( h^{1/2} [p_h^n] \right)^2 \, ds \right)^{1/2} \cdot \left( \int_{\partial K} \left( h^{-1/2} (\boldsymbol{v}_{p_h^n} - \pi_h^1 \boldsymbol{v}_{p_h^n}) \cdot \boldsymbol{n} \right)^2 \, ds \right)^{1/2}.$$

Now apply Young's inequality and (4.16) to obtain

$$B_h^n[(\boldsymbol{u}_h, \boldsymbol{z}_h, p_h), (\pi_h^1 \boldsymbol{v}_{p_h^n}, \boldsymbol{0}, 0)] \geq \|p_h^n\|_{0,\Omega}^2 - \frac{C_c^2}{2\epsilon} \|\boldsymbol{u}_h^n\|_{1,\Omega}^2 - \frac{\epsilon \hat{c}}{2} \|p_h^n\|_{0,\Omega}^2 - \frac{1}{2\epsilon\delta} J(p_h^n, p_h^n) - \frac{\epsilon}{2} \sum_K \int_{\partial K} h^{-1} |(\boldsymbol{v}_{p_h^n} - \pi_h^1 \boldsymbol{v}_{p_h^n}) \cdot \boldsymbol{n}|^2 \, ds.$$

Applying (4.20) we obtain

$$B_h^n[(\boldsymbol{u}_h, \boldsymbol{z}_h, p_h), (\pi_h^1 \boldsymbol{v}_{p_h^n}, \boldsymbol{0}, 0)] \geq -\frac{C_c^2}{2\epsilon} \|\boldsymbol{u}_h^n\|_{1,\Omega}^2 + \left( 1 - (\hat{c} + c_t) \frac{\epsilon}{2} \right) \|p_h^n\|_{0,\Omega}^2 - \frac{1}{2\epsilon\delta} |p_h^n|_{J,\Omega}^2. \quad (4.29)$$

*Step 3, bounding* $\Delta t \|\nabla \cdot \boldsymbol{z}_h^n\|_{0,\Omega}$.

Choosing $(\boldsymbol{v}_h, \boldsymbol{w}_h, q_h) = (\boldsymbol{0}, \boldsymbol{0}, \Delta t \nabla \cdot \boldsymbol{z}_h^n)$ yields

$$B_h^n[(\boldsymbol{u}_h, \boldsymbol{z}_h, p_h), (\boldsymbol{0}, \boldsymbol{0}, \Delta t \nabla \cdot \boldsymbol{z}_h^n)] = (\nabla \cdot \boldsymbol{u}_h^n, \Delta t \nabla \cdot \boldsymbol{z}_h^n) + \Delta t^2 \|\nabla \cdot \boldsymbol{z}_h^n\|_{0,\Omega}^2 + J(p_h^n, \Delta t \nabla \cdot \boldsymbol{z}_h^n).$$

We bound the first term using the Cauchy-Schwarz inequality followed by Young's inequality such that

$$(\nabla \cdot \boldsymbol{u}_h^n, \Delta t \nabla \cdot \boldsymbol{z}_h^n) \leq \frac{C_p}{2\epsilon}\|\boldsymbol{u}_h^n\|_{1,\Omega}^2 + \frac{\epsilon \Delta t^2}{2}\|\nabla \cdot \boldsymbol{z}_h^n\|_{0,\Omega}^2.$$

We can also bound the third term as before using the Cauchy-Schwarz inequality followed by Young's inequality such that

$$\begin{aligned} J(p_h^n, \Delta t \nabla \cdot \boldsymbol{z}_h^n) &\leq \frac{1}{2\epsilon} J(p_h^n, p_h^n) + \frac{\epsilon \Delta t^2}{2} J(\nabla \cdot \boldsymbol{z}_h^n, \nabla \cdot \boldsymbol{z}_h^n) \\ &= \frac{1}{2\epsilon} J(p_h^n, p_h^n) + \epsilon \delta \Delta t^2 \sum_K \int_{\partial K} |h^{1/2} \nabla \cdot \boldsymbol{z}_h^n|^2 \, ds \\ &\leq \frac{1}{2\epsilon} J(p_h^n, p_h^n) + \epsilon \delta c_z \Delta t^2 \|\nabla \cdot \boldsymbol{z}_h^n\|_{0,\Omega}^2. \end{aligned} \tag{4.30}$$

Here we have used the scaling argument (4.10) which relates line and surface integrals and assumes that $\nabla \cdot \boldsymbol{z}_h^n$ is element-wise constant, and (4.11). This yields

$$\begin{aligned} B_h^n[(\boldsymbol{u}_h, \boldsymbol{z}_h, p_h), (\mathbf{0}, \mathbf{0}, \Delta t \nabla \cdot \boldsymbol{z}_h^n)] \geq (1 - \epsilon \delta c_z - \frac{\epsilon}{2}) \Delta t^2 \|\nabla \cdot \boldsymbol{z}_h^n\|_{0,\Omega}^2 \\ - \frac{1}{2\epsilon}|p_h^n|_{J,\Omega}^2 - \frac{C_p}{2\epsilon}\|\boldsymbol{u}_h^n\|_{1,\Omega}^2. \end{aligned} \tag{4.31}$$

*Step 4, Combining steps 1-3.*

Finally we can combine (4.27), (4.29) and (4.31) to get control over all the norms by choosing $(\boldsymbol{v}_h, \boldsymbol{w}_h, q_h) = (\beta \boldsymbol{u}_h^n + \pi_h^1 \boldsymbol{v}_{p_h^n}, \beta \boldsymbol{z}_h^n, \beta p_h^n + \Delta t \nabla \cdot \boldsymbol{z}_h^n)$, which yields

$$
\begin{aligned}
&B_h^n[(\boldsymbol{u}_h, \boldsymbol{z}_h, p_h), (\beta \boldsymbol{u}_h^n + \pi_h^1 \boldsymbol{v}_{p_h^n}, \beta \boldsymbol{z}_h^n, \beta p_h^n + \Delta t \nabla \cdot \boldsymbol{z}_h^n)] \geq \\
&(\beta C_k - \frac{C_c^2 + C_p}{2\epsilon}) \|\boldsymbol{u}_h^n\|_{1,\Omega}^2 + \beta \lambda_{max}^{-1} \Delta t \|\boldsymbol{z}_h^n\|_{0,\Omega}^2 + \left(1 - \epsilon \delta c_z - \frac{\epsilon}{2}\right) \Delta t^2 \|\nabla \cdot \boldsymbol{z}_h^n\|_{0,\Omega}^2 \\
&\qquad + \left(1 - (\hat{c} + c_t) \frac{\epsilon}{2}\right) \|p_h^n\|_{0,\Omega}^2 + \left(\beta - \frac{1}{2\epsilon} - \frac{1}{2\epsilon\delta}\right) |p_h^n|_{J,\Omega}^2, \quad (4.32)
\end{aligned}
$$

where we can choose

$$
\beta \geq \max \left[ \frac{C_c^2 + C_p}{2\epsilon C_k} + \frac{1 - \bar{C}\epsilon}{C_k}, \lambda_{max} \left(1 - \bar{C}\epsilon\right), \frac{1}{2\epsilon} + \frac{1}{2\epsilon\delta} + 1 - \bar{C}\epsilon \right], \quad (4.33)
$$

with $\bar{C} = \max\left[\frac{\hat{c}+c_t}{2}, \delta c_z - \frac{1}{2}\right]$. This yields

$$
B_h^n[(\boldsymbol{u}_h, \boldsymbol{z}_h, p_h), (\beta \boldsymbol{u}_h^n + \pi_h^1 \boldsymbol{v}_{p_h^n}, \beta \boldsymbol{z}_h^n, \beta p_h^n + \nabla \cdot \boldsymbol{z}_h^n)] \geq (1 - \bar{C}\epsilon) \left|\!\left|\!\left| (\boldsymbol{u}_h^n, \boldsymbol{z}_h^n, p_h^n) \right|\!\right|\!\right|_A^2 .
$$

To complete the proof, we let $(\boldsymbol{v}_h, \boldsymbol{w}_h, q_h) = (\beta \boldsymbol{u}_h^n + \pi_h^1 \boldsymbol{v}_{p_h^n}, \beta \boldsymbol{z}_h^n, \beta p_h^n + \Delta t \nabla \cdot \boldsymbol{z}_h^n)$ and show that for $\epsilon$ sufficiently small there exists a constant $C$ such that $\left|\!\left|\!\left| (\boldsymbol{u}_h^n, \boldsymbol{z}_h^n, p_h^n) \right|\!\right|\!\right|_A \geq C \left|\!\left|\!\left| (\boldsymbol{v}_h, \boldsymbol{w}_h, q_h) \right|\!\right|\!\right|_A$. Using the triangle inequality and (4.16) we obtain

$$
\begin{aligned}
&\left|\!\left|\!\left| (\beta \boldsymbol{u}_h^n + \pi_h^1 \boldsymbol{v}_{p_h^n}, \beta \boldsymbol{z}_h^n, \beta p_h^n + \Delta t \nabla \cdot \boldsymbol{z}_h^n) \right|\!\right|\!\right|_A^2 \\
&\leq C \Big( \beta^2 \|\boldsymbol{u}_h^n\|_{1,\Omega}^2 + \left\|\pi_h^1 \boldsymbol{v}_{p_h^n}\right\|_{1,\Omega}^2 + \Delta t^2 (1+\beta)^2 \|\nabla \cdot \boldsymbol{z}_h^n\|_{0,\Omega}^2 + \beta^2 \Delta t \|\boldsymbol{z}_h^n\|_{0,\Omega}^2 \\
&\qquad + \beta^2 \|p_h^n\|_{0,\Omega}^2 + \beta^2 |p_h^n|_{J,\Omega}^2 + \Delta t^2 |\nabla \cdot \boldsymbol{z}_h^n|_{J,\Omega}^2 \Big) \\
&\leq C \left|\!\left|\!\left| (\boldsymbol{u}_h^n, \boldsymbol{z}_h^n, p_h^n) \right|\!\right|\!\right|_A^2 ,
\end{aligned}
$$

as desired. □

Due to the artifical parameters $\epsilon$ and $\beta$ it is difficult to pin down the effect of $\delta$ on the stability (coercivity) of the discretisation. However from looking at (4.33) we can see that $\beta$ needs to be chosen large enough such that $\beta \geq \frac{1}{2\epsilon\delta}$, to ensure that the final coercivity result holds. This suggests that as $\delta \to 0$ we would lose coercivity. Also note that if we were to employ the naive approach of using $J(p_h^n, q_h)$ as the stabilisation term in (4.7), we would require $\beta \geq \frac{1}{2\epsilon\delta\Delta t}$, resulting in a loss of stability as $\Delta t \to 0$. This is shown numerically in Figure 5.4.

## 4.5 Energy estimate for the fully-discrete model

In this section we construct two new combined bilinear forms, $B^n_{\Delta t,h}$ (Lemmas 4.5.1 and 4.5.2) and $\mathcal{B}^n_h$ (Lemmas 4.5.3 and 4.5.4). These bilinear forms are bounded below by Lemmas 4.5.1 and 4.5.3 respectively. Lemma 4.5.2 uses Lemma 4.5.1 to provide a bound on $\boldsymbol{u}_h, \boldsymbol{z}_h$ and $p_h$. Lemma 4.5.4 uses Lemma 4.5.3 to provide a bound on $\nabla \cdot \boldsymbol{z}_h$.

### 4.5.1 Bound on the displacement, fluid flux and pressure

Adding (4.7a), (4.7b) and (4.7c), and assuming $\boldsymbol{t}_N = 0$ on $\Gamma_N$, we get the following

$$B^n_{\Delta t,h}[(\boldsymbol{u}_h, \boldsymbol{z}_h, p_h), (\boldsymbol{v}_h, \boldsymbol{w}_h, q_h)] = (\boldsymbol{f}^n, \boldsymbol{v}_h) + (\boldsymbol{b}^n, \boldsymbol{w}_h) + (g^n, q_h)\ \forall (\boldsymbol{v}_h, \boldsymbol{w}_h, q_h) \in \mathcal{W}^X_h, \tag{4.34}$$

where

$$\begin{aligned} B^n_{\Delta t,h}[(\boldsymbol{u}_h, \boldsymbol{z}_h, p_h), (\boldsymbol{v}_h, \boldsymbol{w}_h, q_h)] &= a(\boldsymbol{u}^n_h, \boldsymbol{v}_h) + (\boldsymbol{k}^{-1}\boldsymbol{z}^n_h, \boldsymbol{w}_h) - (p^n_h, \nabla \cdot \boldsymbol{v}_h) \\ &- (p^n_h, \nabla \cdot \boldsymbol{w}_h) + (\nabla \cdot \boldsymbol{u}^n_{\Delta t,h}, q_h) + (\nabla \cdot \boldsymbol{z}^n_h, q_h) + J(p^n_{\Delta t,h}, q_h). \end{aligned} \tag{4.35}$$

**Lemma 4.5.1.** $(\boldsymbol{u}_h, \boldsymbol{z}_h, p_h)$ *satisfies*

$$\left(\sum_{n=1}^{N} \Delta t B^n_{\Delta t,h}[(\boldsymbol{u}_h, \boldsymbol{z}_h, p_h), (\boldsymbol{u}^n_{\Delta t,h} + \pi^1_h \boldsymbol{v}_{p^n_h}, \boldsymbol{z}^n_h, p^n_h)] + \left\|\boldsymbol{u}^0_h\right\|^2_{1,\Omega} + |p^0_h|^2_{J,\Omega} + \|\boldsymbol{u}_h\|^2_{L^2(H^1)} + \|p_h\|^2_{L^2(J)}\right)$$
$$\geq C\left(\left\|\boldsymbol{u}^N_h\right\|^2_{1,\Omega} + |p^N_h|^2_{J,\Omega} + \|\boldsymbol{z}_h\|^2_{L^2(L^2)} + \|p_h\|^2_{L^2(L^2)}\right).$$

*Proof.* For $n = 1, 2, \ldots, N$ we choose $(\boldsymbol{v}_h, \boldsymbol{w}_h, q_h) = (\boldsymbol{u}^n_{\Delta t,h} + \pi^1_h \boldsymbol{v}_{p^n_h}, \boldsymbol{z}^n_h, p^n_h)$ in (4.35), multiplying by $\Delta t$, and summing over all time steps, we get

$$\sum_{n=1}^{N} \Delta t B^n_{\Delta t,h}[(\boldsymbol{u}_h, \boldsymbol{z}_h, p_h), (\boldsymbol{u}^n_{\Delta t,h} + \pi^1_h \boldsymbol{v}_{p^n_h}, \boldsymbol{z}^n_h, p^n_h)]$$
$$= \sum_{n=1}^{N} \Delta t a(\boldsymbol{u}^n_h, \boldsymbol{u}^n_{\Delta t,h}) + \sum_{n=1}^{N} \Delta t J(p^n_{\Delta t,h}, p^n_h) + \sum_{n=1}^{N} \Delta t \boldsymbol{k^{-1}}(\boldsymbol{z}^n_h, \boldsymbol{z}^n_h)$$
$$+ \sum_{n=1}^{N} \Delta t a(\boldsymbol{u}^n_h, \pi^1_h \boldsymbol{v}_{p^n_h}) - \sum_{n=1}^{N} \Delta t (p^n_h, \nabla \cdot \pi^1_h \boldsymbol{v}_{p^n_h}). \quad (4.36)$$

We now bound each of the above terms on the right hand side of (4.36) individually before combining the results.

$$\sum_{n=1}^{N} \Delta t a(\boldsymbol{u}^n_h, \boldsymbol{u}^n_{\Delta t,h}) = \sum_{n=1}^{N} \Delta t \left(\frac{1}{\Delta t}\|\boldsymbol{u}^n_h\|^2_{a,\Omega} - \frac{1}{\Delta t} a(\boldsymbol{u}^n_h, \boldsymbol{u}^{n-1}_h)\right)$$
$$\geq \frac{C_k}{2}\left\|\boldsymbol{u}^N_h\right\|^2_{1,\Omega} - \frac{C_c}{2}\left\|\boldsymbol{u}^0_h\right\|^2_{1,\Omega}, \quad (4.37)$$

where we have used (4.2) and (4.3) in the last step. Using (4.29) we have

$$\sum_{n=1}^{N} \Delta t a(\boldsymbol{u}^n_h, \pi^1_h \boldsymbol{v}_p) - \sum_{n=1}^{N} \Delta t (p^n_h, \nabla \cdot \pi^1_h \boldsymbol{v}_p) \geq -\frac{C^2_c}{2\epsilon}\|\boldsymbol{u}_h\|^2_{L^2(H^1)}$$
$$+ \left(1 - \left(\hat{c} + \frac{c_t}{2}\right)\frac{\epsilon}{2}\right)\|p_h\|^2_{L^2(L^2)} - \frac{1}{4\epsilon\delta}\|p_h\|^2_{L^2(J)}. \quad (4.38)$$

Using (4.4),

$$\sum_{n=1}^{N} \Delta t(\boldsymbol{k}^{-1}(\boldsymbol{z}_h^n, \boldsymbol{z}_h^n))) \geq \lambda_{max}^{-1} \|\boldsymbol{z}_h\|_{L^2(L^2)}^2. \tag{4.39}$$

The intermediate steps for the next bound have been omitted because they are very similar to (4.37). Thus

$$\sum_{n=1}^{N} \Delta t J(p_{\Delta t,h}^n, p_h^n) \geq \frac{1}{2}|p_h^N|_{J,\Omega}^2 - \frac{1}{2}|p_h^0|_{J,\Omega}^2. \tag{4.40}$$

We can now combine these intermediate results (4.37), (4.38), (4.39) and (4.40) to obtain from (4.36)

$$\begin{aligned}\sum_{n=0}^{N} \Delta t B_{\Delta t,h}^n[(\boldsymbol{u}_h, \boldsymbol{z}_h, p_h), (\boldsymbol{u}_{\Delta t,h}^n + \pi_h^1 \boldsymbol{v}_p, \boldsymbol{z}_h^n, p_h^n)] + \frac{C_c}{2} \left\|\boldsymbol{u}_h^0\right\|_{1,\Omega}^2 \\ + \frac{C_c^2}{2\epsilon}\|\boldsymbol{u}_h\|_{L^2(H^1)}^2 + \frac{1}{4\epsilon\delta}\|p_h\|_{L^2(J)}^2 + \frac{1}{2}|p_h^0|_{J,\Omega}^2 \\ \geq \frac{C_k}{2}\left\|\boldsymbol{u}_h^N\right\|_{1,\Omega}^2 + \frac{1}{2}|p_h^N|_{J,\Omega}^2 + \lambda_{max}^{-1}\|\boldsymbol{z}_h\|_{L^2(L^2)}^2 + (1 - C\epsilon)\,\|p_h\|_{L^2(L^2)}^2. \end{aligned} \tag{4.41}$$

Finally, choosing $\epsilon$ sufficiently small completes the proof.

□

**Lemma 4.5.2.** $(\boldsymbol{u}_h, \boldsymbol{z}_h, p_h)$ *satisfies*

$$\left\|\boldsymbol{u}_h^N\right\|_{1,\Omega}^2 + |p_h^N|_{J,\Omega}^2 + \|\boldsymbol{z}_h\|_{L^2(L^2)}^2 + \|p_h\|_{L^2(L^2)}^2 \leq C(T).$$

*Proof.* For $n = 1, 2, \ldots, N$ we choose $(\boldsymbol{v}_h, \boldsymbol{w}_h, q_h) = (\boldsymbol{u}_{\Delta t,h}^n + \pi_h^1 \boldsymbol{v}_{p_h^n}, \boldsymbol{z}_h^n, p_h^n)$ in

(4.35), multiplying by $\Delta t$, and summing yields

$$\sum_{n=1}^{N} \Delta t B^n_{\Delta t,h}[(\boldsymbol{u}^n_h, \boldsymbol{z}^n_h, p^n_h), (\boldsymbol{u}^n_{\Delta t,h} + \pi^1_h \boldsymbol{v}_{p^n_h}, \boldsymbol{z}^n_h, p^n_h)] = \sum_{n=1}^{N} \Delta t(\boldsymbol{f}^n, \boldsymbol{u}^n_{\Delta t,h} + \pi^1_h \boldsymbol{v}_{p^n_h})$$
$$+ \sum_{n=1}^{N} \Delta t(\boldsymbol{b}^n, \boldsymbol{z}^n_h) + \sum_{n=1}^{N} \Delta t(g^n, p^n_h).$$

Let us note that,

$$\sum_{n=1}^{N} \Delta t(\boldsymbol{f}^n, \boldsymbol{u}^n_{\Delta t,h}) = \sum_{n=1}^{N} (\boldsymbol{f}^n, \boldsymbol{u}^n_h - \boldsymbol{u}^{n-1}_h)$$
$$= (\boldsymbol{f}^N, \boldsymbol{u}^N_h) - (\boldsymbol{f}^1, \boldsymbol{u}^0_h) - \sum_{n=1}^{N-1} (\boldsymbol{f}^{n+1} - \boldsymbol{f}^n, \boldsymbol{u}^n_h), \quad (4.42)$$

and further that

$$-\sum_{n=1}^{N-1} (\boldsymbol{f}^{n+1} - \boldsymbol{f}^n, \boldsymbol{u}^n_h) \leq C \sum_{n=1}^{N-1} \left\| \boldsymbol{f}^{n+1} - \boldsymbol{f}^n \right\|_{0,\Omega} \| \boldsymbol{u}^n_h \|_{0,\Omega}$$
$$\leq C \sum_{n=1}^{N-1} \left\{ \int_{t_n}^{t_{n+1}} \| \boldsymbol{f}_t \|_{0,\Omega} \right\}^{1/2} \| \boldsymbol{u}^n_h \|_{1,\Omega} \leq C \left( \frac{1}{2\epsilon} \| \boldsymbol{f}_t \|^2_{L^2(L^2)} + \frac{\epsilon}{2} \| \boldsymbol{u}_h \|^2_{L^2(L^2)} \right).$$

Now using the above, Lemma 4.5.1, the Cauchy-Schwarz and Young's inequalities, choosing $\epsilon$ sufficiently small, and noting (4.16), we arrive at

$$\left\| \boldsymbol{u}^N_h \right\|^2_{1,\Omega} + |p^N_h|^2_{J,\Omega} + \| \boldsymbol{z}_h \|^2_{L^2(L^2)} + \| p_h \|^2_{L^2(L^2)} \leq C \Bigg( \| \boldsymbol{u}_h \|^2_{L^2(H^1)} + \frac{1}{\delta} \| p_h \|^2_{L^2(J)} + \left\| \boldsymbol{f}^N \right\|^2_{0,\Omega}$$
$$+ \| \boldsymbol{f}_t \|^2_{L^2(L^2)} + \left\| \boldsymbol{u}^0_h \right\|^2_{0,\Omega} + |p^0_h|^2_{J,\Omega} + \left\| \boldsymbol{f}^1 \right\|^2_{L^2(L^2)} + \| \boldsymbol{f} \|^2_{L^2(L^2)} + \| \boldsymbol{b} \|^2_{L^2(L^2)} + \| g \|^2_{L^2(L^2)} \Bigg).$$

Using assumed regularity of the given data to bound the third term and upwards on the righthand side we obtain

$$\left\| \boldsymbol{u}^N_h \right\|^2_{1,\Omega} + |p^N_h|^2_{J,\Omega} + \| \boldsymbol{z}_h \|^2_{L^2(L^2)} + \| p_h \|^2_{L^2(L^2)} \leq C \left( 1 + \| \boldsymbol{u}_h \|^2_{L^2(H^1)} + \| p_h \|^2_{L^2(J)} \right).$$

Upon applying the Gronwall lemma to the above inequality we obtain the desired result.

□

### 4.5.2 Bound on the divergence of the fluid flux

In order to bound the divergence of the fluid flux we now define the bilinear form $\mathcal{B}_h^n$. We first show how we derive $\mathcal{B}_h^n$ from the fully-discrete weak form (4.7), for which we know that a solution $(\boldsymbol{u}_h, \boldsymbol{z}_h, p_h)$ exists for test functions $(\boldsymbol{v}_h, \boldsymbol{w}_h, q_h) \in \mathcal{V}_h^X$. Adding (4.7a) and (4.7b), assuming $\boldsymbol{t}_N = 0$ on $\Gamma_N$, and summing we have

$$\sum_{n=1}^{N} a(\boldsymbol{u}_h^n, \boldsymbol{v}_h) + \sum_{n=1}^{N} (\boldsymbol{k}^{-1}\boldsymbol{z}_h^n, \boldsymbol{w}_h) - \sum_{n=1}^{N} (p_h^n, \nabla \cdot \boldsymbol{v}_h) - \sum_{n=1}^{N} (p_h^n, \nabla \cdot \boldsymbol{w}_h)$$
$$= \sum_{n=1}^{N} (\boldsymbol{f}^n, \boldsymbol{v}_h) + \sum_{n=1}^{N} (\boldsymbol{b}^n, \boldsymbol{w}_h) \quad \forall (\boldsymbol{v}_h, \boldsymbol{w}_h, q_h) \in \mathcal{V}_h^X. \quad (4.43)$$

For the purposes of this proof we now introduce initial conditions for the fluid flux and the pressure, $\boldsymbol{z}^0 \in H_{div}(\Omega)$ and $p^0 \in \mathcal{L}(\Omega_t)$ respectively. We also define their projections into their respective finite element spaces by $\boldsymbol{z}_h^0 := \pi_h^0 \boldsymbol{z}^0$ and $p_h^0 := \pi_h^0 p^0$.

Adding (4.7a) and (4.7b), and summing from 0 to $N-1$, we have

$$\sum_{n=1}^{N} a(\boldsymbol{u}_h^{n-1}, \boldsymbol{v}_h) + \sum_{n=1}^{N} (\boldsymbol{k}^{-1}\boldsymbol{z}_h^{n-1}, \boldsymbol{w}_h) - \sum_{n=1}^{N} (p_h^{n-1}, \nabla \cdot \boldsymbol{v}_h) - \sum_{n=1}^{N} (p_h^{n-1}, \nabla \cdot \boldsymbol{w}_h)$$
$$= \sum_{n=1}^{N} (\boldsymbol{f}^{n-1}, \boldsymbol{v}_h) + \sum_{n=1}^{N} (\boldsymbol{b}^{n-1}, \boldsymbol{w}_h) \quad \forall (\boldsymbol{v}_h, \boldsymbol{w}_h, q_h) \in \mathcal{V}_h^X. \quad (4.44)$$

Taking (4.7c), multiplying by $\Delta t$, and summing we have

$$\sum_{n=1}^{N}\Delta t(\nabla\cdot\boldsymbol{u}^n_{\Delta t,h},q_h)+\sum_{n=1}^{N}\Delta t(\nabla\cdot\boldsymbol{z}^n_h,q_h)+\sum_{n=1}^{N}\Delta tJ(p^n_{\Delta t,h},q_h)$$
$$=\sum_{n=1}^{N}\Delta t(g^n,q_h)\quad\forall(\boldsymbol{v}_h,\boldsymbol{w}_h,q_h)\in\mathcal{V}^X_h. \quad (4.45)$$

Now adding (4.43) and (4.45), and subtracting (4.44) we get

$$\sum_{n=1}^{N}\Delta t\mathcal{B}^n_h[(\boldsymbol{u}_h,\boldsymbol{z}_h,p_h),(\boldsymbol{v}_h,\boldsymbol{w}_h,q_h)]$$
$$=\sum_{n=1}^{N}\Delta t(\boldsymbol{f}^n_{\Delta t},\boldsymbol{v}_h)+\sum_{n=1}^{N}\Delta t(\boldsymbol{b}^n_{\Delta t},\boldsymbol{w}_h)+\sum_{n=1}^{N}\Delta t(g^n,q_h)\ \forall(\boldsymbol{v}_h,\boldsymbol{w}_h,q_h)\in\mathcal{V}^X_h,$$

where

$$\mathcal{B}^n_h[(\boldsymbol{u}_h,\boldsymbol{z}_h,p_h),(\boldsymbol{v}_h,\boldsymbol{w}_h,q_h)]=a(\boldsymbol{u}^n_{\Delta t,h},\boldsymbol{v}_h)+(\boldsymbol{k^{-1}}\boldsymbol{z}^n_{\Delta t,h},\boldsymbol{w}_h)$$
$$-(p^n_{\Delta t,h},\nabla\cdot\boldsymbol{v}_h)-(p^n_{\Delta t,h},\nabla\cdot\boldsymbol{w}_h)+(\nabla\cdot\boldsymbol{u}^n_{\Delta t,h},q_h)+(\nabla\cdot\boldsymbol{z}^n_h,q_h)+J(p^n_{\Delta t,h},q_h). \quad (4.46)$$

With these preliminaries, we may now bound $\mathcal{B}^n_h$ from below.

**Lemma 4.5.3.** *For all $\beta>\beta^\star>0$, $(\boldsymbol{u}_h,\boldsymbol{z}_h,p_h)$ satisfies*

$$\sum_{n=1}^{N}\Delta t\ \mathcal{B}^n_h[(\boldsymbol{u}_h,\boldsymbol{z}_h,p_h),(\beta\boldsymbol{u}^n_{\Delta t,h}+\pi^1_hv_p,\beta\boldsymbol{z}^n_h,\beta p^n_{\Delta t,h}+\nabla\cdot\boldsymbol{z}^n_h)]+\left\|\boldsymbol{z}^0_h\right\|^2_{0,\Omega}\geq$$
$$C\left(\|\boldsymbol{u}_{\Delta t,h}\|^2_{L^2(H^1)}+\left\|\boldsymbol{z}^N_h\right\|^2_{0,\Omega}+\|p_{\Delta t,h}\|^2_{L^2(L^2)}+\|p_{\Delta t,h}\|^2_{L^2(J)}+\|\nabla\cdot\boldsymbol{z}_h\|^2_{L^2(L^2)}\right).$$

*where $\beta^\star$ takes the value of $\beta$ previously chosen in (4.33).*

*Proof.* For $n=1,2,\dots,N$ we choose $(\boldsymbol{v}_h,\boldsymbol{w}_h,q_h)=(\beta\boldsymbol{u}^n_{\Delta t,h}+\pi^1_h\boldsymbol{v}_{p^n_h},\beta\boldsymbol{z}^n_h,\beta p^n_{\Delta t,h}+$

$\nabla \cdot \boldsymbol{z}_h^n)$ in (4.46)

$$
\begin{aligned}
&\sum_{n=1}^{N} \Delta t \mathcal{B}_h^n[(\boldsymbol{u}_h, \boldsymbol{z}_h, p_h), (\beta \boldsymbol{u}_{\Delta t,h}^n + \pi_h^1 \boldsymbol{v}_p, \beta \boldsymbol{z}_h^n, \beta p_{\Delta t,h}^n + \nabla \cdot \boldsymbol{z}_h^n)] \\
&\quad = \sum_{n=1}^{N} \Delta t a(\boldsymbol{u}_{\Delta t,h}^n, \beta \boldsymbol{u}_{\Delta t,h}^n) + \sum_{n=1}^{N} \Delta t \boldsymbol{k}^{-1}(\boldsymbol{z}_{\Delta t,h}^n, \beta \boldsymbol{z}_h^n) + \sum_{n=1}^{N} \Delta t (\nabla \cdot \boldsymbol{z}_h^n, \nabla \cdot \boldsymbol{z}_h^n) \\
&\quad + \sum_{n=1}^{N} \Delta t (\nabla \cdot \boldsymbol{u}_{\Delta t,h}^n, \nabla \cdot \boldsymbol{z}_h^n) + \sum_{n=1}^{N} \Delta t J(p_{\Delta t,h}^n, \nabla \cdot \boldsymbol{z}_h^n) + \sum_{n=1}^{N} \Delta t J(p_{\Delta t,h}^n, \beta p_{\Delta t,h}^n) \\
&\qquad + \sum_{n=1}^{N} \Delta t a(\boldsymbol{u}_{\Delta t,h}^n, \pi_h^1 \boldsymbol{v}_p) - \sum_{n=1}^{N} \Delta t (p_{\Delta t,h}^n, \nabla \cdot \pi_h^1 \boldsymbol{v}_p). \quad (4.47)
\end{aligned}
$$

For all $\epsilon > 0$ using (4.3), (4.4), the Cauchy-Schwarz, Young's and Poincaré inequalities, (4.11) and (4.10) on $\nabla \cdot \boldsymbol{z}_h^n$, and an approach similar to step 2 in the proof of Theorem 4.4.1 for the final two terms on the righthand side, we obtain

$$
\begin{aligned}
&\sum_{n=1}^{N} \Delta t \mathcal{B}_h^n[(\boldsymbol{u}_h, \boldsymbol{z}_h, p_h), (\beta \boldsymbol{u}_{\Delta t,h}^n + \pi_h^1 \boldsymbol{v}_p, \beta \boldsymbol{z}_h^n, \beta p_{\Delta t,h}^n + \nabla \cdot \boldsymbol{z}_h^n)] \\
&\geq \left(\beta C_k - \frac{C_p + C_c^2}{2\epsilon}\right) \|\boldsymbol{u}_{\Delta t,h}\|_{L^2(H^1)}^2 + \frac{\beta \lambda_{max}^{-1}}{2} \left\|\boldsymbol{z}_h^N\right\|_{0,\Omega}^2 + \left(\beta - \frac{1}{2\epsilon} - \frac{1}{2\epsilon\delta}\right) \|p_{\Delta t,h}\|_{L^2(J)}^2 \\
&\quad + \left(1 - \epsilon\delta c_z - \frac{\epsilon}{2}\right) \|\nabla \cdot \boldsymbol{z}_h\|_{L^2(L^2)}^2 - \frac{\beta \lambda_{min}^{-1}}{2} \left\|\boldsymbol{z}_h^0\right\|_{0,\Omega}^2 + (1 - C\epsilon) \|p_{\Delta t,h}\|_{L^2(L^2)}^2.
\end{aligned}
\quad (4.48)
$$

Finally choosing $\epsilon$ sufficiently small and $\beta \geq \max\left[\frac{C_p + C_c^2}{2C_k\epsilon}, \frac{1}{2\epsilon} + \frac{1}{2\epsilon\delta}\right]$ completes the proof.

□

The following Lemma shows the divergence control of the fluid flux.

**Lemma 4.5.4.** $\boldsymbol{z}_h$ *obtained from (4.46) satisfies*

$$
\|\nabla \cdot \boldsymbol{z}_h\|_{L^2(L^2)}^2 \leq C(T).
$$

*Proof.* For $n = 1, 2, \dots, N$ we choose $(\boldsymbol{v}_h, \boldsymbol{w}_h, q_h) = (\beta \boldsymbol{u}^n_{\Delta t,h} + \pi^1_h \boldsymbol{v}_{p^n_h}, \beta \boldsymbol{z}^n_h, \beta p^n_{\Delta t,h} + \nabla \cdot \boldsymbol{z}^n_h)$ in (4.46) yielding

$$\begin{aligned}
&\sum_{n=1}^{N} \Delta t \mathcal{B}^n_h[(\boldsymbol{u}^n_h, \boldsymbol{z}^n_h, p^n_h), (\beta \boldsymbol{u}^n_{\Delta t,h} + \pi^1_h \boldsymbol{v}_{p^n_h}, \boldsymbol{z}^n_h, \beta p^n_{\Delta t,h} + \nabla \cdot \boldsymbol{z}^n_h)] \\
&\quad = \sum_{n=1}^{N} \Delta t (\boldsymbol{f}^n_{\Delta t}, \beta \boldsymbol{u}^n_{\Delta t,h} + \pi^1_h \boldsymbol{v}_{p^n_h}) + \sum_{n=1}^{N} \Delta t (\boldsymbol{b}^n_{\Delta t}, \beta \boldsymbol{z}^n_h) \\
&\quad + \sum_{n=1}^{N} \Delta t (g^n, \beta p^n_{\Delta t,h} + \nabla \cdot \boldsymbol{z}^n_h).
\end{aligned}$$

Using Lemma 4.5.3, the Cauchy-Schwarz and Young's inequalities, and (4.16), along with ideas already presented in the proof of Lemma 4.5.2

$$\begin{aligned}
&\|\boldsymbol{u}_{\Delta t,h}\|^2_{L^2(H^1)} + \|p_{\Delta t,h}\|^2_{L^2(L^2)} + \|p_{\Delta t,h}\|^2_{L^2(J)} + \left\|\boldsymbol{z}^N_h\right\|^2_{0,\Omega} + \|\nabla \cdot \boldsymbol{z}_h\|^2_{L^2(L^2)} \\
&\quad \le C \left( \|\boldsymbol{f}_t\|^2_{L^2(L^2)} + \|\boldsymbol{b}_t\|^2_{L^2(L^2)} + \|p_h\|^2_{L^2(L^2)} + \|\boldsymbol{z}_h\|^2_{L^2(L^2)} + \|g\|^2_{L^2(L^2)} \right).
\end{aligned}$$

Finally, using Lemma 4.5.2 to bound $\|p_h\|_{L^2(L^2)}$, applying a Gronwall lemma, and using regularity, we obtain the desired result.

□

### 4.5.3 The energy estimate

**Theorem 4.5.5.** *The solution to the fully-discrete problem (4.7) satisfies the energy estimate*

$$\|\boldsymbol{u}_h\|^2_{L^\infty(H^1)} + \|p_h\|^2_{L^\infty(J)} + \|\boldsymbol{z}_h\|^2_{L^2(L^2)} + \|p_h\|^2_{L^2(L^2)} + \|\nabla \cdot \boldsymbol{z}_h\|^2_{L^2(L^2)} \le C.$$

*Proof.* The proof follows from combining Lemma 4.5.2 and Lemma 4.5.4, and noting that these lemmas hold for all time steps $n = 0, 1, ..., N$. This then gives the desired discrete in time $L^\infty$ bounds. □

**Remark 4.5.1.** *Having proven Theorem 4.5.5, it is now a standard calculation to show that the discrete Galerkin approximation converges weakly, as $\Delta t, h \to 0$, to the continuous problem with respect to continuous versions of the norms of the energy estimate in Theorem 4.5.5. This in turn shows that the continuous variational problem is well-posed. Due to the linearity of the variational form and noting that $|\boldsymbol{v}|_{J,\Omega} \to 0$ as $h \to 0$, these calculations are straight forward and closely follow the existence and uniqueness proofs presented in Ženíšek (1984) and Barucq et al. (2005) for the linear two-field Biot problem and a nonlinear Biot problem, respectively.*

## 4.6 A-priori error analysis

Lemma 4.6.1 provides a Galerkin orthogonality result obtained by comparing continuous and discrete weak forms, which is the corner stone of the error analysis. Lemma 4.6.2 bounds the auxiliary errors for displacement, flux and pressure in the appropriate norms and Lemma 4.6.3 bounds the auxiliary error for the divergence of the flux. Since Lemmas 4.6.2 and 4.6.3 bound the auxiliary errors at the same order as the projection errors, combining projection and auxiliary errors in Theorem 4.6.4 provides an optimal error estimate.

We define the finite element error functions

$$\boldsymbol{e}_{\boldsymbol{u}} = \boldsymbol{u} - \boldsymbol{u}_h, \quad \boldsymbol{e}_{\boldsymbol{z}} = \boldsymbol{z} - \boldsymbol{z}_h, \quad e_p = p - p_h.$$

We introduce the following projection errors:

$$\boldsymbol{\eta}_{\boldsymbol{u}} = \boldsymbol{u} - \pi_h^1 \boldsymbol{u}, \quad \boldsymbol{\eta}_{\boldsymbol{z}} = \boldsymbol{z} - \pi_h^1 \boldsymbol{z}, \quad \eta_p = p - \pi_h^0 p,$$

where we have assumed $\boldsymbol{z}(\cdot, t_n) \in (H^1(\Omega))^d$. Auxiliary errors:

$$\boldsymbol{\theta}_{\boldsymbol{u}}^n(\cdot) = \pi_h^1 \boldsymbol{u}(\cdot, t_n) - \boldsymbol{u}_h^n(\cdot),\ \boldsymbol{\theta}_{\boldsymbol{z}}^n(\cdot) = \pi_h^1 \boldsymbol{z}(\cdot, t_n) - \boldsymbol{z}_h^n(\cdot),\ \theta_p^n(\cdot) = \pi_h^0 p(\cdot, t_n) - p_h^n(\cdot), \tag{4.49}$$

and time-discretisation errors:

$$\boldsymbol{\rho}_{\boldsymbol{u}}^n(\cdot) = \frac{\boldsymbol{u}(\cdot, t_n) - \boldsymbol{u}(\cdot, t_{n-1})}{\Delta t} - \frac{\partial \boldsymbol{u}(\cdot, t_n)}{\partial t},\ \rho_p^n = \frac{p(\cdot, t_n) - p(\cdot, t_{n-1})}{\Delta t} - \frac{\partial p(\cdot, t_n)}{\partial t}. \tag{4.50}$$

### 4.6.1 Galerkin orthogonality

We now give a Galerkin orthogonality type argument for analysing the difference between the fully-discrete approximation and the true solution. For this we introduce the continuous counterpart of the fully-discrete combined weak form (4.34) given by

$$B^n[(\boldsymbol{u}, \boldsymbol{z}, p), (\boldsymbol{v}, \boldsymbol{w}, q)] = (\boldsymbol{f}(\cdot, t_n), \boldsymbol{v}) + (\boldsymbol{b}(\cdot, t_n), \boldsymbol{w}) + (g(\cdot, t_n), q)\ \forall (\boldsymbol{v}, \boldsymbol{w}, q) \in \mathcal{V}^X, \tag{4.51}$$

where

$$\begin{aligned} B^n[(\boldsymbol{u}, \boldsymbol{z}, p), (\boldsymbol{v}, \boldsymbol{w}, q)] \ &= \ a(\boldsymbol{u}(\cdot, t_n), \boldsymbol{v}) + \boldsymbol{k}^{-1}(\boldsymbol{z}(\cdot, t_n), \boldsymbol{w}) - (p(\cdot, t_n), \nabla \cdot \boldsymbol{v}) \\ &\quad -(p(\cdot, t_n), \nabla \cdot \boldsymbol{w}) + (\nabla \cdot \boldsymbol{u}_t(\cdot, t_n), q) + (\nabla \cdot \boldsymbol{z}(\cdot, t_n), q). \end{aligned}$$

**Lemma 4.6.1.** *Assuming* $(\boldsymbol{u}(\cdot, t_n), \boldsymbol{z}(\cdot, t_n), p(\cdot, t_n)) \in (H^1(\Omega))^d \times H_{div}(\Omega) \times (H^1(\Omega) \cap \mathcal{L}(\Omega_t))$

$$B_{\Delta t,h}^n[(\boldsymbol{e}_{\boldsymbol{u}}, \boldsymbol{e}_{\boldsymbol{z}}, e_p), (\boldsymbol{v}_h, \boldsymbol{w}_h, q_h)] = (\nabla \cdot \boldsymbol{\rho}_{\boldsymbol{u}}^n, q_h) + J(\rho_p^n, q_h)\ \ \forall (\boldsymbol{v}_h, \boldsymbol{w}_h, q_h) \in \mathcal{V}_h^X.$$

*Proof.* Subtracting the discrete weak form (4.34) from the continuous weak form

(4.51), we obtain

$$B^n[(\boldsymbol{u},\boldsymbol{z},p),(\boldsymbol{v}_h,\boldsymbol{w}_h,q_h)]-B^n_{\Delta t,h}[(\boldsymbol{u}_h,\boldsymbol{z}_h,p_h),(\boldsymbol{v}_h,\boldsymbol{w}_h,q_h)]=0,\ \ \forall(\boldsymbol{v}_h,\boldsymbol{w}_h,q_h)\in\mathcal{V}_h^X.$$

Now add $J(p_t(\cdot,t_n),q)=0$ to the left hand side, see (4.12). Finally add $(\nabla\cdot(\boldsymbol{u}_{\Delta t}(\cdot,t_n)-\boldsymbol{u}_t(\cdot,t_n)),q)+J(p_{\Delta t}(\cdot,t_n)-p_t(\cdot,t_n),q)$ to the left and the righthand side to obtain the desired result. □

### 4.6.2 Auxiliary error estimates

**Lemma 4.6.2.**

$$|||[\boldsymbol{\theta}_{\boldsymbol{u}},\boldsymbol{\theta}_{\boldsymbol{z}},\theta_p]|||_B^2+\|\theta_p\|^2_{L^\infty(J)}\le C(\delta,T)(h^2+\Delta t^2). \tag{4.52}$$

*Proof.* Using Lemma 4.6.1 and choosing $\boldsymbol{v}_h^n=\boldsymbol{\theta}^n_{\Delta t,\boldsymbol{u}}+\pi_h^1\boldsymbol{v}_{p_h^n}$, $\boldsymbol{w}_h^n=\boldsymbol{\theta}^n_{\boldsymbol{z}}$, $q_h^n=\theta_p^n$, we get

$$\begin{aligned}&B^n_{\Delta t,h}[(\boldsymbol{\theta}^n_{\boldsymbol{u}}+\boldsymbol{\eta}^n_{\boldsymbol{u}},\boldsymbol{\theta}^n_{\boldsymbol{z}}+\boldsymbol{\eta}^n_{\boldsymbol{z}},\theta^n_p+\eta^n_p),(\boldsymbol{\theta}^n_{\Delta t,\boldsymbol{u}}+\pi_h^1\boldsymbol{v}_{p_h^n},\boldsymbol{\theta}^n_{\boldsymbol{z}},\theta^n_p)]\\&\quad=(\nabla\cdot\boldsymbol{\rho}^n_{\boldsymbol{u}},\theta^n_p)+J(\rho^n_p,\theta^n_p).\end{aligned}$$

Rearranging gives

$$\begin{aligned}&B^n_{\Delta t,h}[(\boldsymbol{\theta}^n_{\boldsymbol{u}},\boldsymbol{\theta}^n_{\boldsymbol{z}},\theta^n_p),(\boldsymbol{\theta}^n_{\Delta t,\boldsymbol{u}}+\pi_h^1\boldsymbol{v}_{p_h^n},\boldsymbol{\theta}^n_{\boldsymbol{z}},\theta^n_p)]\\&\quad=(\nabla\cdot\boldsymbol{\rho}^n_{\boldsymbol{u}},\theta^n_p)+J(\rho^n_p,\theta^n_p)-B^n_{\Delta t,h}[(\boldsymbol{\eta}^n_{\boldsymbol{u}},\boldsymbol{\eta}^n_{\boldsymbol{z}},\eta^n_p),(\boldsymbol{\theta}^n_{\Delta t,\boldsymbol{u}}+\pi_h^1\boldsymbol{v}_{p_h^n},\boldsymbol{\theta}^n_{\boldsymbol{z}},\theta^n_p)].\end{aligned}$$

Expanding the righthand side, noting that $(\eta^n_p,\nabla\cdot(\boldsymbol{\theta}^n_{\Delta t,\boldsymbol{u}}+\pi_h^1\boldsymbol{v}_{ph}))=0$, $(\eta^n_p,\nabla\cdot\boldsymbol{\theta}^n_{\boldsymbol{z}})=0$, multiplying both sides by $\Delta t$ and summing gives

$$\sum_{n=1}^N\Delta t B^n_{\Delta t,h}[(\boldsymbol{\theta}^n_{\boldsymbol{u}},\boldsymbol{\theta}^n_{\boldsymbol{z}},\theta^n_p),(\boldsymbol{\theta}^n_{\Delta t,\boldsymbol{u}}+\pi_h^1\boldsymbol{v}_{p_h^n},\boldsymbol{\theta}^n_{\boldsymbol{z}},\theta^n_p)]=\sum_{i=1}^7\Phi_i, \tag{4.53}$$

where

$$\Phi_1 = -\sum_{n=1}^{N} \Delta t a(\boldsymbol{\eta}_{\boldsymbol{u}}^n, \boldsymbol{\theta}_{\Delta t,\boldsymbol{u}}^n), \qquad \Phi_2 = -\sum_{n=1}^{N} \Delta t(\boldsymbol{k}^{-1}(\boldsymbol{\eta}_{\boldsymbol{z}}^n, \boldsymbol{\theta}_{\boldsymbol{z}}^n)),$$

$$\Phi_3 = -\sum_{n=1}^{N} \Delta t a(\boldsymbol{\eta}_{\boldsymbol{u}}^n, \pi_h^1 \boldsymbol{v}_p), \qquad \Phi_4 = -\sum_{n=1}^{N} \Delta t J(\eta_{\Delta t,p}^n, \theta_p^n),$$

$$\Phi_5 = \sum_{n=1}^{N} \Delta t(\nabla \cdot \boldsymbol{\rho}_{\boldsymbol{u}}^n, \theta_p^n), \qquad \Phi_6 = \sum_{n=1}^{N} \Delta t J(\rho_p^n, \theta_p^n),$$

$$\Phi_7 = -\sum_{n=1}^{N} \Delta t(\theta_p^n, \nabla \cdot (\boldsymbol{\eta}_{\Delta t,\boldsymbol{u}}^n + \boldsymbol{\eta}_{\boldsymbol{z}}^n)).$$

We now individually consider the terms on the right hand side of (4.53): To bound the first quantity, we use (4.21), Lemma 4.3.1, the triangle, Cauchy-Schwarz and Young's inequalities, $\boldsymbol{\theta}_{\boldsymbol{u}}^0 = \mathbf{0}$, and (4.2),

$$\begin{aligned}
\Phi_1 &= -\sum_{n=1}^{N} a(\boldsymbol{\eta}_{\boldsymbol{u}}^n, \boldsymbol{\theta}_{\boldsymbol{u}}^n - \boldsymbol{\theta}_{\boldsymbol{u}}^{n-1}) \\
&= -a(\boldsymbol{\eta}_{\boldsymbol{u}}^N, \boldsymbol{\theta}_{\boldsymbol{u}}^N) + \sum_{n=1}^{N} a(\boldsymbol{\eta}_{\boldsymbol{u}}^n - \boldsymbol{\eta}_{\boldsymbol{u}}^{n-1}, \boldsymbol{\theta}_{\boldsymbol{u}}^{n-1}) \\
&= -a(\boldsymbol{\eta}_{\boldsymbol{u}}^N, \boldsymbol{\theta}_{\boldsymbol{u}}^N) + \Delta t \sum_{n=1}^{N} a\left(\left(I - \pi_h^1\right)\left(\boldsymbol{\rho}_{\boldsymbol{u}}^n + \frac{\partial \boldsymbol{u}(\cdot, t_n)}{\partial t}\right), \boldsymbol{\theta}_{\boldsymbol{u}}^{n-1}\right) \\
&\leq \epsilon C \left\|\boldsymbol{\theta}_{\boldsymbol{u}}^N\right\|_{1,\Omega}^2 + \frac{Ch^2}{\epsilon}\left\|\boldsymbol{u}^N\right\|_{2,\Omega}^2 + \epsilon C\|\boldsymbol{\theta}_{\boldsymbol{u}}\|_{L^2(H^1)}^2 + \frac{Ch^2}{2\epsilon}\|\boldsymbol{u}_t\|_{L^2(H^2)}^2 \\
&\quad + \frac{C\Delta t^2}{2\epsilon}\|\boldsymbol{u}_{tt}\|_{L^2(H^1)}^2.
\end{aligned}$$

Next, using (4.4), Young's inequality, (4.16) and Lemma 4.3.1,

$$\Phi_2 \leq \frac{\epsilon}{2}\|\boldsymbol{\theta}_{\boldsymbol{z}}\|_{L^2(L^2)}^2 + \frac{\lambda_{min}^{-2} h^2}{2\epsilon}\|\boldsymbol{z}\|_{L^2(H^1)}^2.$$

Using (4.2), Young's inequality and Lemma 4.3.1,

$$\Phi_3 \leq \frac{\epsilon}{2}\left\|\pi_h^1 \boldsymbol{v}_{p_h^n}\right\|_{L^2(H^1)}^2 + \frac{C}{2\epsilon}\|\boldsymbol{\eta}_{\boldsymbol{u}}\|_{L^2(H^1)}^2 \leq \frac{\epsilon\hat{c}^2}{2}\|\theta_p\|_{L^2(L^2)}^2 + \frac{Ch^2}{2\epsilon}\|\boldsymbol{u}\|_{L^2(H^2)}^2.$$

The bound on $\Phi_4$ is obtained using a similar argument to the bound on $\Phi_1$,

$$\Phi_4 \leq \epsilon \|\theta_p\|^2_{L^2(J)} + \frac{\delta C h^2}{2\epsilon}\|p_t\|^2_{L^2(H^1)} + \frac{\delta C \Delta t^2}{2\epsilon}\|p_{tt}\|^2_{L^2(H^1)}.$$

Using the Cauchy-Schwarz and Young's inequalities and Lemma 4.3.1,

$$\Phi_5 \leq \frac{\epsilon}{2}\|\theta_p\|^2_{L^2(L^2)} + \frac{\Delta t^2}{2\epsilon}\|\boldsymbol{u}_{tt}\|^2_{L^2(L^2)} \text{ and } \Phi_6 \leq \frac{\epsilon}{2}\|\theta_p\|^2_{L^2(J)} + \frac{\delta C \Delta t^2}{2\epsilon}\|p_{tt}\|^2_{L^2(L^2)}.$$

As can be seen from the bound on $\Phi_4$ and $\Phi_6$ we lose control of the auxillary error if $\delta$ is very large. This is reflected in the numerical experiments in Chapter 5, where simulations with a large $\delta$ carry a larger error. Further if we were to employ $\frac{1}{\Delta t}J(p_h^n, q_h)$ as the stabilisation term in (4.7), which would result in an inf-sup stable method and pass the proof of Theorem 4.4.1, we would now have $\Phi_4 \leq \epsilon\|\theta_p\|^2_{L^2(J)} + \frac{\delta C h^2}{2\epsilon \Delta t^2}\|p_t\|^2_{L^2(H^1)} + \frac{\delta C}{2\epsilon}\|p_{tt}\|^2_{L^2(H^1)}$. This would cause the error to increase as $\Delta t \to 0$, shown numerically in Figure 5.4. Thus the choice of stabilisation, $J(p^n_{\Delta t,h}, q_h)$, is key to creating a stable and converging method.

Finally, using the Cauchy-Schwarz and Young's inequalities, and a similar argument to the bound on $\Phi_1$,

$$\Phi_7 \leq \frac{3\epsilon}{2}\|\theta_p\|^2_{L^2(L^2)} + \frac{h^2}{2\epsilon}\|\boldsymbol{u}_t\|^2_{L^2(H^2)} + \frac{\Delta t^2}{2\epsilon}\|\boldsymbol{u}_{tt}\|^2_{L^2(H^1)} + \frac{h^2}{2\epsilon}\|\boldsymbol{z}\|^2_{L^2(H^2)}.$$

Combining these bounds with an application of coercivity Lemma 4.5.1 to (4.53), noting the assumed regularity of the continuous solution and choosing $\epsilon$ sufficiently small, gives

$$\|\boldsymbol{\theta}_{\boldsymbol{u}}^N\|^2_{1,\Omega} + |\theta_p^N|^2_{J,\Omega} + \|\boldsymbol{\theta}_{\boldsymbol{z}}\|^2_{L^2(L^2)} + \|\theta_p\|^2_{L^2(L^2)} \leq C(\delta)\left(\|\boldsymbol{\theta}_{\boldsymbol{u}}\|^2_{L^2(H^1)} + \|\theta_p\|^2_{L^2(J)} + h^2 + \Delta t^2\right). \tag{4.54}$$

An application of Gronwall's lemma gives

$$\|\boldsymbol{\theta}_{\boldsymbol{u}}^N\|_{1,\Omega}^2 + |\theta_p^N|_{J,\Omega}^2 + \|\boldsymbol{\theta}_{\boldsymbol{z}}\|_{L^2(L^2)}^2 + \|\theta_p\|_{L^2(L^2)}^2 \le C(\delta,T)\left(h^2+\Delta t^2\right).$$

Because the above holds for all time steps $n = 0, 1, ..., N$, we can get the desired $L^\infty$ bounds to complete the proof of the theorem. □

We now present an a-priori auxiliary error estimate of the fluid flux, in its natural $H_{div}$ norm.

**Lemma 4.6.3.** *Assuming* $\boldsymbol{u} \in H^2\left(0,T;(H^1(\Omega))^d\right) \cap H^1\left(0,T;(H^2(\Omega))^d\right)$, $\boldsymbol{z} \in L^2\left(0,T;(H^2(\Omega))^d\right)$ *and* $p \in H^2(0,T;J\cap\mathcal{L}(\Omega_t)) \cap H^1(0,T;H^1(\Omega))$, *then the finite element solution (4.7) satisfies the auxillary error estimate*

$$\|\nabla\cdot\boldsymbol{\theta}_{\boldsymbol{z}}\|_{L^2(L^2)}^2 \le C(\delta,T)(h^2+\Delta t^2). \tag{4.55}$$

*Proof.* Similarly to the approach taken in obtaining (4.46) we may easily obtain the following identity

$$\sum_{n=1}^{N}\Delta t\mathcal{B}_h^n[(\boldsymbol{\theta}_{\boldsymbol{u}}^n,\boldsymbol{\theta}_{\boldsymbol{z}}^n,\theta_p^n),(\beta\boldsymbol{\theta}_{\Delta t,\boldsymbol{u}}^n+\pi_h^1\boldsymbol{v}_{\theta_{\Delta t,p}^n},\beta\boldsymbol{\theta}_{\boldsymbol{z}}^n,\beta\theta_{\Delta t,p}^n+\nabla\cdot\boldsymbol{\theta}_{\boldsymbol{z}}^n)] = \sum_{i=1}^{6}\Psi_i, \tag{4.56}$$

where

$$\Psi_1 = -\sum_{n=1}^{N}\Delta t a(\boldsymbol{\eta}_{\Delta t,\boldsymbol{u}}^n,\beta\boldsymbol{\theta}_{\Delta t,\boldsymbol{u}}^n+\pi_h^1\boldsymbol{v}_{\theta_{\Delta t,p}^n}),$$
$$\Psi_2 = -\sum_{n=1}^{N}\Delta t(\nabla\cdot(\boldsymbol{\eta}_{\Delta t,\boldsymbol{u}}^n+\boldsymbol{\eta}_{\boldsymbol{z}}^n),\nabla\cdot\boldsymbol{\theta}_{\boldsymbol{z}}^n+\beta\theta_{\Delta t,p}^n),$$
$$\Psi_3 = \sum_{n=1}^{N}\Delta tJ(\eta_{\Delta t,p}^n,\beta\theta_{\Delta t,p}^n+\nabla\cdot\boldsymbol{\theta}_{\boldsymbol{z}}^n), \quad \Psi_4 = -\sum_{n=1}^{N}\Delta t(\boldsymbol{k}^{-1}(\boldsymbol{\eta}_{\Delta t,\boldsymbol{z}}^n,\beta\boldsymbol{\theta}_{\boldsymbol{z}}^n)),$$
$$\Psi_5 = \sum_{n=1}^{N}\Delta tJ(\rho_p^n,\beta\theta_{\Delta t,p}^n+\nabla\cdot\boldsymbol{\theta}_{\boldsymbol{z}}^n), \quad \Psi_6 = \sum_{n=1}^{N}\Delta t(\nabla\cdot\boldsymbol{\rho}_{\boldsymbol{u}}^n,\beta\theta_{\Delta t,p}^n+\nabla\cdot\boldsymbol{\theta}_{\boldsymbol{z}}^n).$$

We now bound the terms on the right hand side of (4.56) using machinery developed during the previous proof:

$$
\begin{aligned}
\Psi_1 &\leq \frac{C\epsilon}{2}\|\boldsymbol{\theta}_{\Delta t,\boldsymbol{u}}\|^2_{L^2(H^1)} + \frac{\hat{c}^2\epsilon}{2}\|\theta_{\Delta t,p}\|^2_{L^2(L^2)} + \frac{Ch^2}{2\epsilon}\|\boldsymbol{u}_t\|^2_{L^2(H^2)} \\
&\quad + \frac{C}{2\epsilon}\Delta t^2\|\boldsymbol{u}_{tt}\|^2_{L^2(H^1)}, \qquad (4.57)\\
\Psi_2 &\leq \epsilon\|\nabla\cdot\boldsymbol{\theta}_{\boldsymbol{z}}\|^2_{L^2(L^2)} + \epsilon\|\theta_{\Delta t,p}\|^2_{L^2(L^2)} + \frac{Ch^2}{2\epsilon}\left(\|\boldsymbol{u}_t\|^2_{L^2(H^2)} + \|\boldsymbol{z}\|^2_{L^2(H^2)}\right) \\
&\quad + \frac{C}{2\epsilon}\Delta t^2\|\boldsymbol{u}_{tt}\|^2_{L^2(H^1)}, \qquad (4.58)\\
\Psi_3 &\leq \epsilon C\|\nabla\cdot\boldsymbol{\theta}_{\boldsymbol{z}}\|^2_{L^2(L^2)} + \epsilon\left\|\theta^n_{\Delta t,p}\right\|^2_{L^2(J)} + \frac{\delta Ch^2}{2\epsilon}\|p_t\|^2_{L^2(H^1)} \\
&\quad + \frac{\delta C}{2\epsilon}\Delta t^2\|p_{tt}\|^2_{L^2(J)}, \qquad (4.59)\\
\Psi_4 &\leq \epsilon\|\boldsymbol{\theta}_{\boldsymbol{z}}\|^2_{L^2(L^2)} + \frac{Ch^2}{2\epsilon}\|\boldsymbol{z}_t\|^2_{L^2(H^1)} + \frac{C}{2\epsilon}\Delta t^2\|\boldsymbol{z}_{tt}\|^2_{L^2(L^2)}, \qquad (4.60)\\
\Psi_5 &\leq \epsilon\|\theta_{\Delta t,p}\|^2_{L^2(J)} + \epsilon C\|\nabla\cdot\boldsymbol{\theta}_{\boldsymbol{z}}\|^2_{L^2(L^2)} + \frac{C\Delta t^2}{2\epsilon}\|p_{tt}\|^2_{L^2(J)}, \qquad (4.61)\\
\Psi_6 &\leq \epsilon\|\theta_{\Delta t,p}\|^2_{L^2(L^2)} + \epsilon\|\nabla\cdot\boldsymbol{\theta}_{\boldsymbol{z}}\|^2_{L^2(L^2)} + \frac{C}{2\epsilon}\Delta t^2\|\boldsymbol{u}_{tt}\|^2_{L^2(H^1)}. \qquad (4.62)
\end{aligned}
$$

We can now combine the individual bounds (4.57), (4.58), (4.59), (4.60), (4.61), and (4.62), with the coercivity result Lemma 4.5.3, choose $\beta$ sufficiently large, use the assumption $\boldsymbol{\theta}^0_{\boldsymbol{z}} = \mathbf{0}$, the assumed regularity of $\boldsymbol{u}, \boldsymbol{z}$ and $p$, and choose $\epsilon$ sufficiently small to obtain

$$
\left\|\boldsymbol{\theta}^N_{\boldsymbol{z}}\right\|^2_{0,\Omega} + \|\nabla\cdot\boldsymbol{\theta}_{\boldsymbol{z}}\|^2_{L^2(L^2)} \leq C\|\boldsymbol{\theta}_{\boldsymbol{z}}\|^2_{L^2(L^2)} + C(\delta)(h^2 + \Delta t^2).
$$

Applying Gronwall's lemma, we get the desired result. □

### 4.6.3 The a-priori error estimate

Combining the previous lemmas we have the following.

**Theorem 4.6.4.** *Assuming* $\boldsymbol{u} \in H^2\left(0,T;(L^2(\Omega))^d\right) \cap H^1\left(0,T;(H^2(\Omega))^d\right)$, $\boldsymbol{z} \in L^2\left(0,T;(H^1(\Omega))^d\right)$ *and* $p \in H^2\left(0,T;H^1(\Omega)\cap\mathcal{L}(\Omega_t)\right)$, *then the finite element* *solution (4.7) satisfies the error estimate*

$$|||\boldsymbol{e_u}, \boldsymbol{e_z}, e_p|||_B^2 \leq C(h^2 + \Delta t^2).$$

*Assuming* $\boldsymbol{u} \in H^2\left(0,T;(H^1(\Omega))^d\right) \cap H^1\left(0,T;(H^2(\Omega))^d\right)$, $\boldsymbol{z} \in L^2\left(0,T;(H^2(\Omega))^d\right)$ *and* $p \in H^2(0,T;J \cap \mathcal{L}(\Omega_t)) \cap H^1(0,T;H^1(\Omega))$, *then the finite element solution (4.7) satisfies the error estimate*

$$|||\boldsymbol{e_u}, \boldsymbol{e_z}, e_p|||_B^2 + \|\nabla \cdot \boldsymbol{e_z}\|_{L^2(L^2)}^2 \leq C(h^2 + \Delta t^2).$$

*Proof.* We first write the errors as $\boldsymbol{e}_{\boldsymbol{u}}^n = \boldsymbol{\eta}_{\boldsymbol{u}}^n + \boldsymbol{\theta}_{\boldsymbol{u}}^n$, and similarly for the other variables. Using lemma 4.3.1 we can bound the projection errors, and using lemma 4.6.2 and lemma 4.6.3 we can bound the auxillary errors to give the desired result. □

## 4.7 Conclusion

The local pressure jump stabilisation method (Burman and Hansbo, 2007) is commonly used to solve the Stokes or Darcy equations using piecewise linear approximations for the velocities, and piecewise constant approximations for the pressure variable. The main contribution of this chapter has been to extend these ideas to three-field poroelasticity. We have presented a stability result for the discretised equations that guarantees the existence of a unique solution at each time step, and derived an energy estimate which can be used to prove weak convergence of the solution to the discretised system to the solution to the continuous problem as the mesh parameters tend to zero. We also derived an optimal error estimate which includes an error for the fluid flux in its natural $H_{div}$ norm.

# Chapter 5

# Numerical results for the stabilised finite element method

The contents of this chapter closely follows the numerical results section presented in the joint publication: L. Berger, R. Bordas, D. Kay, and S. Tavener; Stabilized low-order finite element approximation for linear three-field poroelasticity *SIAM Journal on Scientific Computing* 2015. The numerical tests were designed by L. Berger, with guidance from D. Kay and R. Bordas, and were implemented by L. Berger. S. Tavener assisted in improving the quality of the writing along with the other authors.

## 5.1 Introduction

In this chapter we detail the implementation of the finite element method presented in the previous chapter (section 5.2), followed by numerical experiments that illustrate the convergence of the method and its ability to overcome pressure oscillations. We present convergence studies for both two- and three-dimensional test problems which illustrate the predicted convergence rates for the fully-discrete finite element method. We then apply our method to the popular 2D cantilever bracket problem and demonstrate that our stabilisation technique overcomes the spurious pressure oscillations that have been experienced by other methods. Finally, a 3D unconfined compression problem is presented that highlights the added mass effect of the method for different choices of the stabilisation parameter $\delta$.

## 5.2 Implementation

For the implementation we used the C++ library libMesh (Kirk et al., 2006), and the multi-frontal direct solver mumps (Amestoy et al., 2000) to solve the resulting linear system. To solve the full Biot model problem (2.39), we need to solve the following linear system at each time step:

$$\begin{bmatrix} \boldsymbol{A} & 0 & \alpha \boldsymbol{B}^T \\ 0 & \Delta t \boldsymbol{M} & \Delta t \boldsymbol{B}^T \\ \alpha \boldsymbol{B} & \Delta t \boldsymbol{B} & -c_0 \boldsymbol{Q} - \boldsymbol{J} \end{bmatrix} \begin{bmatrix} \boldsymbol{u}^n \\ \boldsymbol{z}^n \\ \boldsymbol{p}^n \end{bmatrix} = \begin{bmatrix} \boldsymbol{r} \\ \Delta t \boldsymbol{s} \\ \boldsymbol{B}\boldsymbol{u}^{n-1} - c_0 \boldsymbol{Q}\boldsymbol{p}^{n-1} - \boldsymbol{J}\boldsymbol{p}^{n-1} - \Delta t\, \boldsymbol{g} \end{bmatrix}, \tag{5.1}$$

where we have defined the following matrices and vectors:

$$\boldsymbol{A} = [\boldsymbol{a}_{ij}], \;\; \boldsymbol{a}_{ij} = \int_\Omega 2\mu_s \nabla \boldsymbol{\phi}_i : \nabla \boldsymbol{\phi}_j + \lambda (\nabla \cdot \boldsymbol{\phi}_i)(\nabla \cdot \boldsymbol{\phi}_j),$$

$$\boldsymbol{M} = [\boldsymbol{m}_{ij}], \;\; \boldsymbol{m}_{ij} = \int_\Omega k^{-1} \boldsymbol{\phi}_i \cdot \boldsymbol{\phi}_j,$$

$$\boldsymbol{B} = [\boldsymbol{b}_{ij}], \;\; \boldsymbol{b}_{ij} = -\int_\Omega \psi_i \nabla \cdot \boldsymbol{\phi}_j,$$

$$\boldsymbol{Q} = [\boldsymbol{q}_{ij}], \;\; \boldsymbol{q}_{ij} = \int_\Omega \psi_i \cdot \psi_j,$$

$$\boldsymbol{J} = [\boldsymbol{j}_{ij}], \;\; \boldsymbol{j}_{ij} = \delta \sum_K \int_{\partial k \setminus \partial \Omega} h_{\partial K} [\psi_i][\psi_j] \,\mathrm{d}s,$$

$$\boldsymbol{r} = [\boldsymbol{r}_i], \;\; \boldsymbol{r}_i = \int_\Omega \boldsymbol{f}_i \cdot \phi_i + \int_{\Gamma_N} \boldsymbol{t}_{Ni} \cdot \phi_i,$$

$$\boldsymbol{s} = [\boldsymbol{s}_i], \quad \boldsymbol{s}_i = \int_\Omega \boldsymbol{b}_i \cdot \phi_i - \int_{\Gamma_P} p_D \phi_i \cdot \boldsymbol{n},$$

$$\boldsymbol{g} = [\boldsymbol{g}_i], \quad \boldsymbol{g}_i = \int_\Omega g \psi_i.$$

Here $\boldsymbol{\phi}_i$ are vector valued linear basis functions such that the displacement vector can be written as $\boldsymbol{u}^n = \sum_{i=1}^{n_u} \boldsymbol{u}_i^n \boldsymbol{\phi}_i$, with $\sum_{i=1}^{n_u} \boldsymbol{u}_i^n \boldsymbol{\phi}_i \in \boldsymbol{W}_h^E$. Similarly for the fluid flux vector we have $\boldsymbol{z}^n = \sum_{i=1}^{n_z} \boldsymbol{z}_i^n \boldsymbol{\phi}_i$, with $\sum_{i=1}^{n_z} \boldsymbol{z}_i^n \boldsymbol{\phi}_i \in \boldsymbol{W}_h^D$. The scalar valued constant basis functions $\psi_i$ are used to approximate the pressure, such that $\boldsymbol{p}^n = \sum_{i=1}^{n_p} p_i^n \psi_i$, with $\sum_{i=1}^{n_p} p_i^n \psi_i \in Q_h$.

### 5.2.1 Algorithm to assemble the stabilisation matrix

Let $K \in \mathcal{T}_h$ be an element and $\mathcal{D}(K)$ be the pressure degree of freedom associated with element $K$. We define $\mathcal{A}(K)$ to be the set of elements $L \in \mathcal{T}_h$ neighbouring $K$.

```
for every K ∈ T_h do
  for every L ∈ A(K) do
    Calculate h_∂K
    i ← D(K)
    j ← D(L)
    J_ii ← J_ii + (δh²_∂K in 2D, δh³_∂K in 3D)
    J_ij ← J_ij − (δh²_∂K in 2D, δh³_∂K in 3D)
  end for
end for
```

Figure 5.1: Stabilisation matrix $\boldsymbol{J}$ assembly

## 5.3 2D test problem

Choosing $\lambda = \mu = \alpha = 1$, $c_0 = 0$ and $\boldsymbol{k} = \boldsymbol{I}$ in (4.1) we solve the problem

$$-2\nabla(\nabla \cdot \boldsymbol{u}) - \nabla^2 \boldsymbol{u} + \nabla p = \boldsymbol{f} \quad \text{in } \Omega, \tag{5.2a}$$

$$\boldsymbol{z} + \nabla p = 0 \quad \text{in } \Omega, \tag{5.2b}$$

$$\nabla \cdot (\boldsymbol{u}_t + \boldsymbol{z}) = g \quad \text{in } \Omega, \tag{5.2c}$$

$$\boldsymbol{u}(t) = \boldsymbol{u}_D \quad \text{on } \Gamma_d, \tag{5.2d}$$

$$\boldsymbol{z}(t) \cdot \boldsymbol{n} = q_D \quad \text{on } \Gamma_f, \tag{5.2e}$$

$$\boldsymbol{u}(0) = 0, \;\; p(0) = 0 \quad \text{in } \Omega. \tag{5.2f}$$

The domain, $\Omega$, is the unit square and the source terms and boundary conditions are chosen so that the true solution is

$$\boldsymbol{u} = \begin{pmatrix} -\frac{1}{4\pi}\cos(2\pi x)\sin(2\pi y)\sin(2\pi t) \\ -\frac{1}{4\pi}\sin(2\pi x)\cos(2\pi y)\sin(2\pi t) \end{pmatrix}, \;\; \boldsymbol{z} = \begin{pmatrix} -2\pi\cos(2\pi x)\sin(2\pi y)\sin(2\pi t) \\ -2\pi\sin(2\pi x)\cos(2\pi y)\sin(2\pi t) \end{pmatrix},$$

and $p = \sin(2\pi x)\sin(2\pi y)\sin(2\pi t)$, with $t \in [0, 0.25]$.

### 5.3.1 Choice of $\delta$

The most appropriate choice of stabilisation parameter $\delta$ is not known a priori. Small values of $\delta$ can result in spurious pressure solutions, as shown in Figure 5.2a for $\delta = 0.1$. Larger values of the stabilisation parameter produce smooth pressure solutions, as shown in Figure 5.2b for a value of $\delta = 1$. The value of $\delta$ required to produce a stable solution depends on the geometry and material parameters of the particular problem under investigation, but is independent of any mesh parameters. However, care should be taken that $\delta$ does not get chosen to be excessively large. Due to the global nature of the stabilisation, this can cause loss of incompressibility as shown in Figure 5.10. In the extreme case of pressure jump stabilisation (penalization), $\delta \to \infty$, the pressure will tend to a constant solution. This loss in accuracy has already been highlighted in the error analysis performed in section 4.6.2. To circumvent this issue a local stabilisation method has been developed for the closely related Stokes problem where stabilisation is performed on individual macroelements within the mesh, avoiding coupling throughout the whole domain (Kechkar and Silvester, 1992; Silvester and Kechkar, 1990). Using this local stabilisation approach, even in the extreme case of $\delta \to \infty$, the pressure will only tend to a constant solution on each individual macroelement but remain discontinuous between macroelements. It is therefore very robust with respect to the magnitude of the stabilisation parameter and prevents loss of incompressibility (Kay and Silvester, 1999; Kechkar and Silvester, 1992).

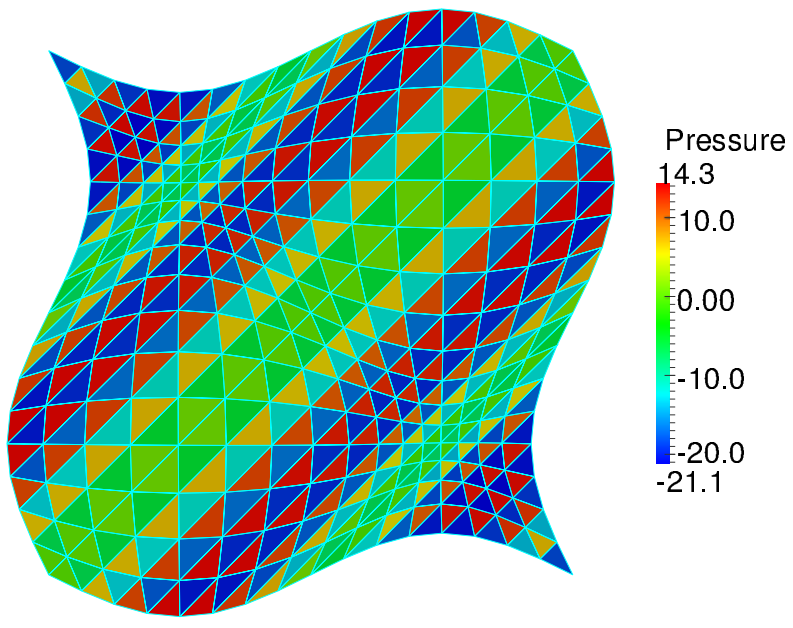

(a)

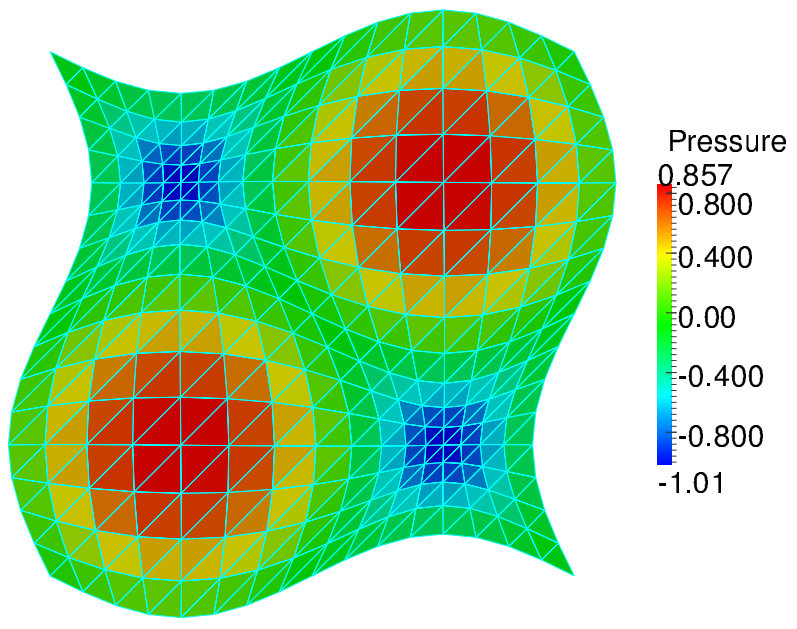

(b)

Figure 5.2: (a) Unstable pressure field, caused by not choosing the stabilisation parameter $\delta$ large enough, with $\delta = 0.1$, at $t = 0.25$. (b) Stable pressure field, with $\delta = 1$ at $t = 0.25$.

### 5.3.2 2D convergence study

The convergence of the method with discretisation parameters is illustrated in Figure 5.3a – 5.3e for $\delta = 1, 10, 100$. The convergence rates observed in the appropriate norms agree with the theoretically derived error estimates.

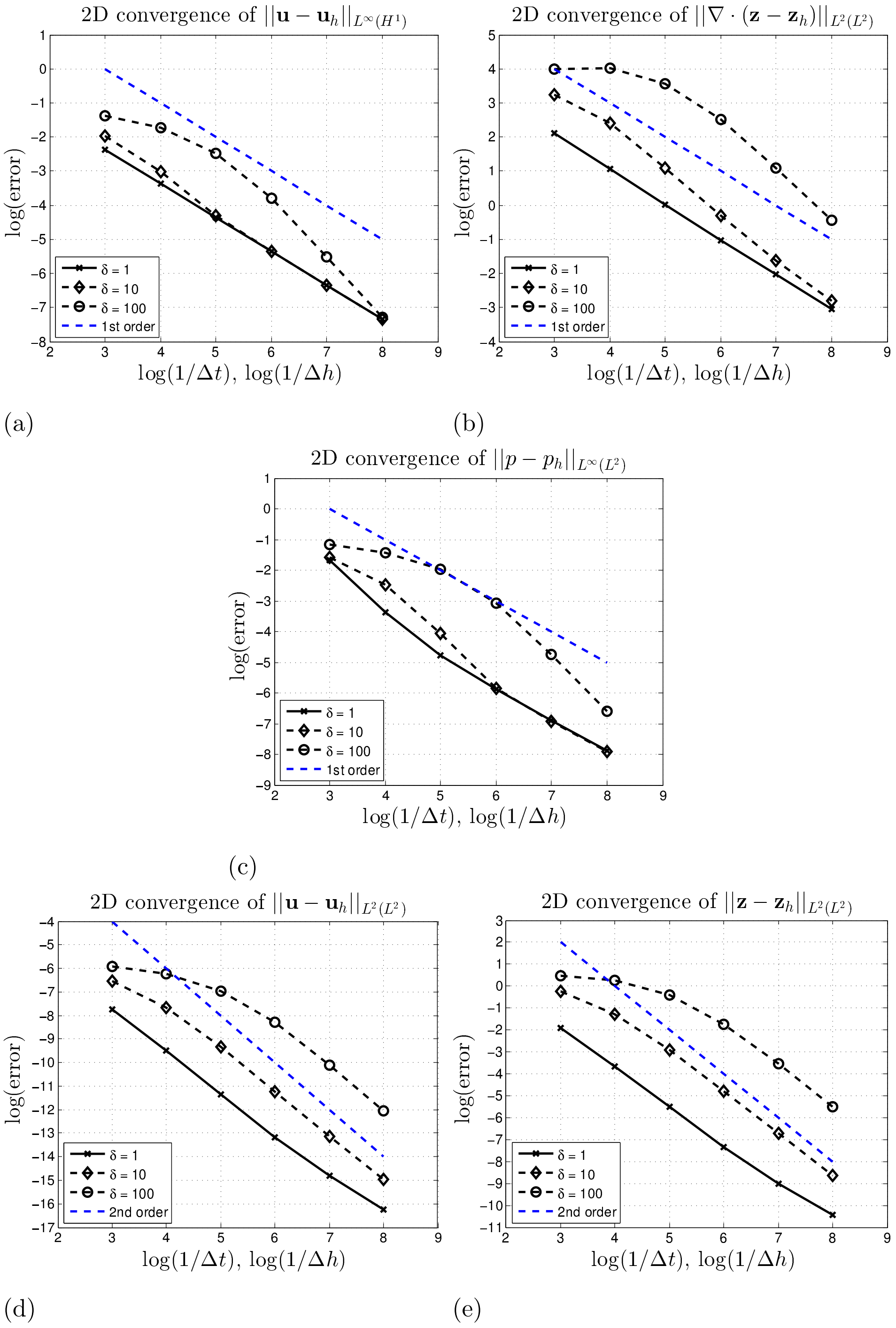

Figure 5.3: Convergence of the displacement, fluid flux, and pressure errors in their respective norms of the simplified poroelastic 2D test problem with different (stable) values for the stabilisation parameter $\delta$.

### 5.3.3 Alternative stabilisation techniques

In Figure 5.4 we illustrate the convergence of the pressure error for three possible stabilisation forms with decreasing time step. The test problem is the same as in the previous 2D convergence study, with $\delta = 1$, and $t \in [0, 0.025]$. As demonstrated in section 5.3, the stabilisation $J(p_{\Delta t,h}, q_h)$ yields a stable solution (Figure 5.5a) and converges, for all sizes of $\Delta t$, to the spatial error (Figure 5.5b), as expected. The more naive approach, inserting the stabilisation $J(p_h, q_h)$, results in the solution becoming unstable, and introducing an oscillating pressure mode into the approximation, see Figure 5.5c and Figure 5.5d. This is because the stabilisation becomes relatively small as $\Delta t$ decreases. Also note that the final refinement step is not possible when using $J(p_{\Delta t,h}, q_h)$ in Figure 5.4 because the numerical solver that solves the resulting linear system breaks due to the relative large pressure mode present in the solution, see Figure 5.6a and Figure 5.6b, showing the pressure solution and error after the first time step, for the last possible refinement level of $\Delta t$. To overcome this issue one could chose to scale the stabilisation, and try $\frac{1}{\Delta t} J(p_h, q_h)$. Although this stabilisation now stays stable during refinement, it does not converge. Instead the error builds up with decreasing $\Delta t$ and the stabilisation starts to dominate the solution by preventing any jumps in pressure and causing extreme smoothing of the pressure, as seen in Figure 5.5e, where the pressure solution is now almost zero throughout.

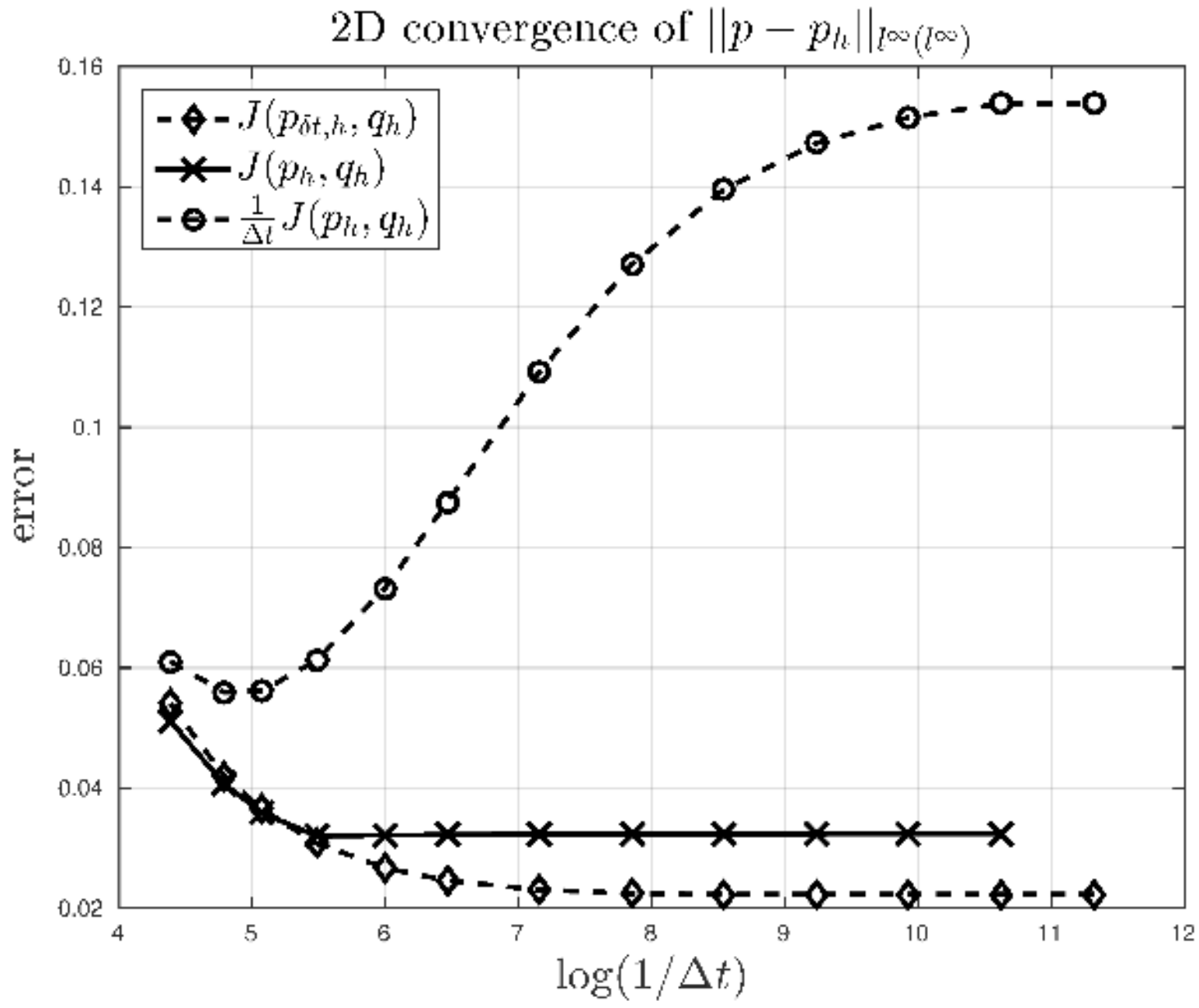

Figure 5.4: Convergence of the pressure error for three different stabilisation forms, with $\delta = 1$.

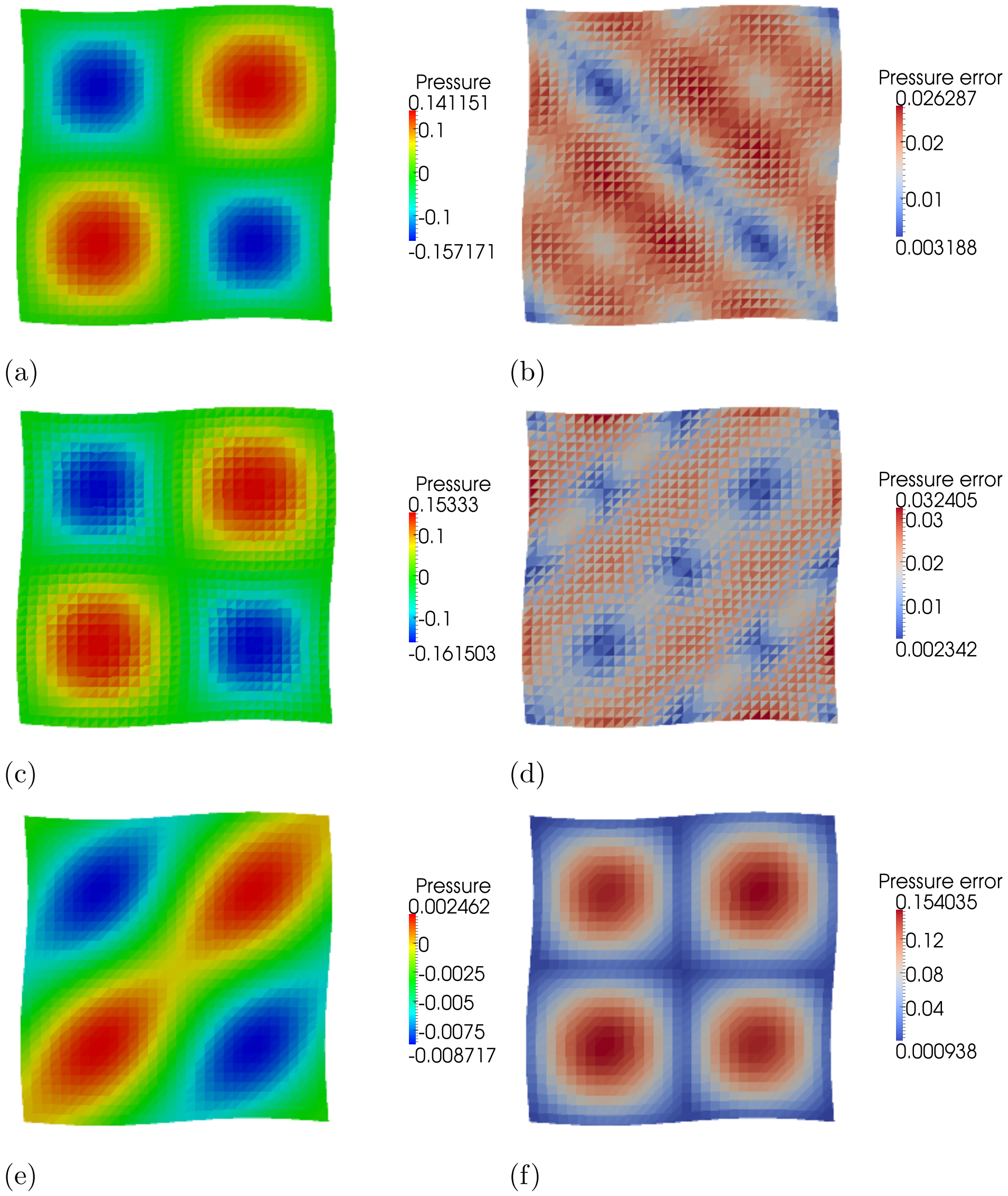

Figure 5.5: Pressure solution and pressure error after 1028 timesteps at $t = 0.025$ using the stabilisation $J(p_{\Delta t,h}, q_h)$, (a) and (b), $J(p_h, q_h)$, (c) and (d), and $\frac{1}{\Delta t}J(p_h, q_h)$, (e) and (f).

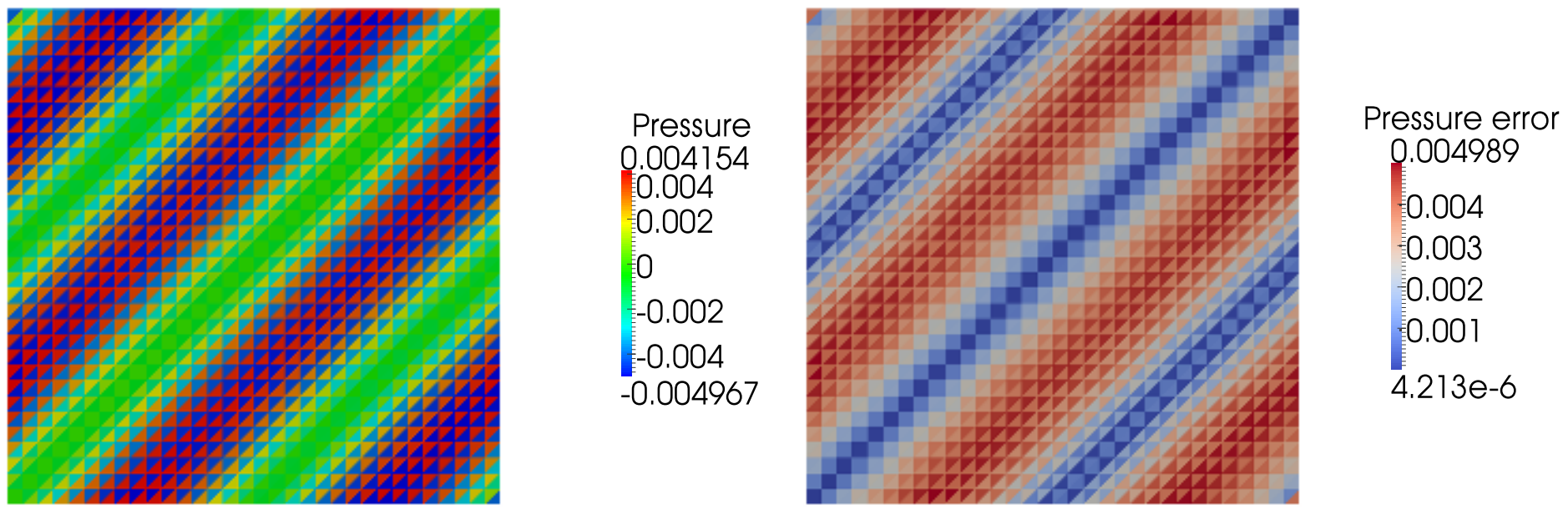

(a) (b)

Figure 5.6: Pressure solution (a) and pressure error (b) after the first time step at ($t = 0.025s/1028$) using the stabilisation $J(p_h, q_h)$.

## 5.4 3D test problem

Extending the test problem in Section 5.3 to the unit cube, we set

$$\boldsymbol{u} = \begin{pmatrix} -\frac{1}{6\pi}\cos(2\pi x)\sin(2\pi y)\sin(2\pi z)\sin(2\pi t) \\ -\frac{1}{6\pi}\sin(2\pi x)\cos(2\pi y)\sin(2\pi z)\sin(2\pi t) \\ -\frac{1}{6\pi}\sin(2\pi x)\sin(2\pi y)\cos(2\pi z)\sin(2\pi t) \end{pmatrix},$$

$$\boldsymbol{z} = \begin{pmatrix} -2\pi\cos(2\pi x)\sin(2\pi y)\sin(2\pi z)\sin(2\pi t) \\ -2\pi\sin(2\pi x)\cos(2\pi y)\sin(2\pi z)\sin(2\pi t) \\ -2\pi\sin(2\pi x)\sin(2\pi y)\cos(2\pi z)\sin(2\pi t) \end{pmatrix},$$

and

$$p = \sin(2\pi x)\sin(2\pi y)\sin(2\pi z)\sin(2\pi t).$$

The expected rates of convergence for each variable in the appropriate norm are illustrated in the numerical results presented in figure 5.7a – 5.7e for $\delta = 0.001, 0.01, 0.1$. The stabilisation factor $\delta$ may be chosen to be very much smaller for 3D problems as compared to 2D problems and the effect of the stabilisation term on the solution is negligible. This can be explained by the improved ratio of solid displacement and fluid flux nodes to pressure nodes in three dimensions, making the LBB condition easier to satisfy.

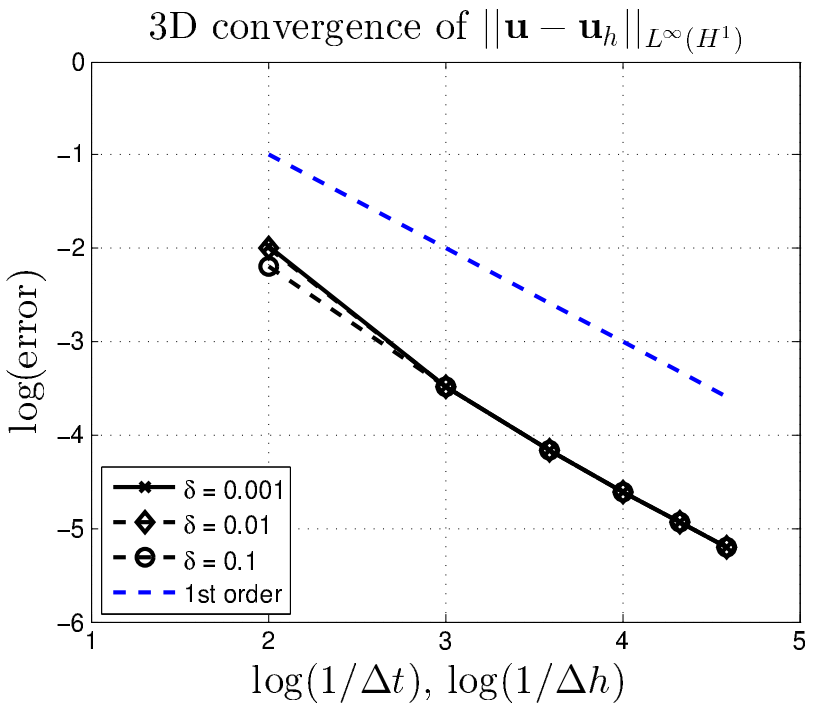

(a)

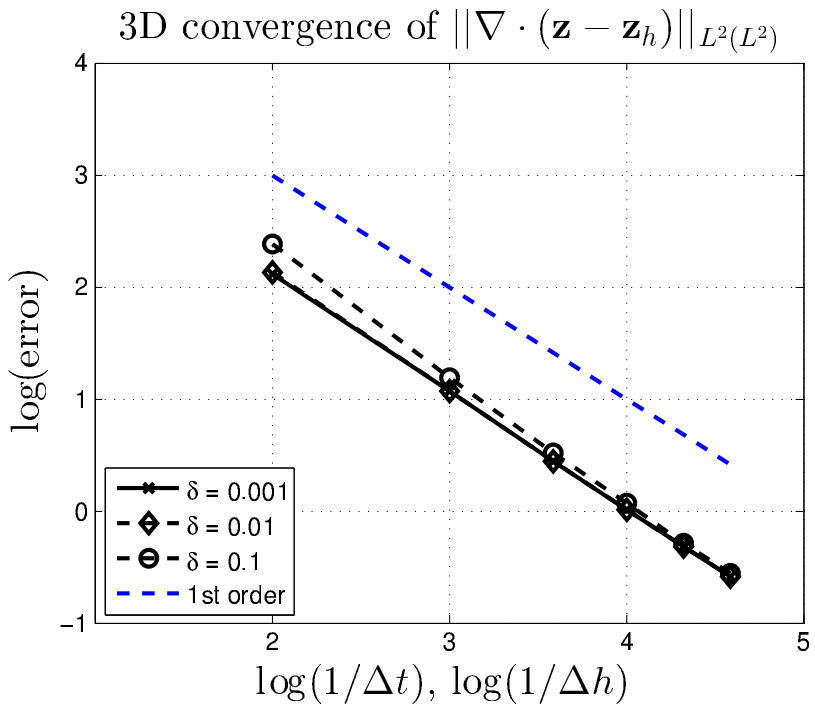

(b)

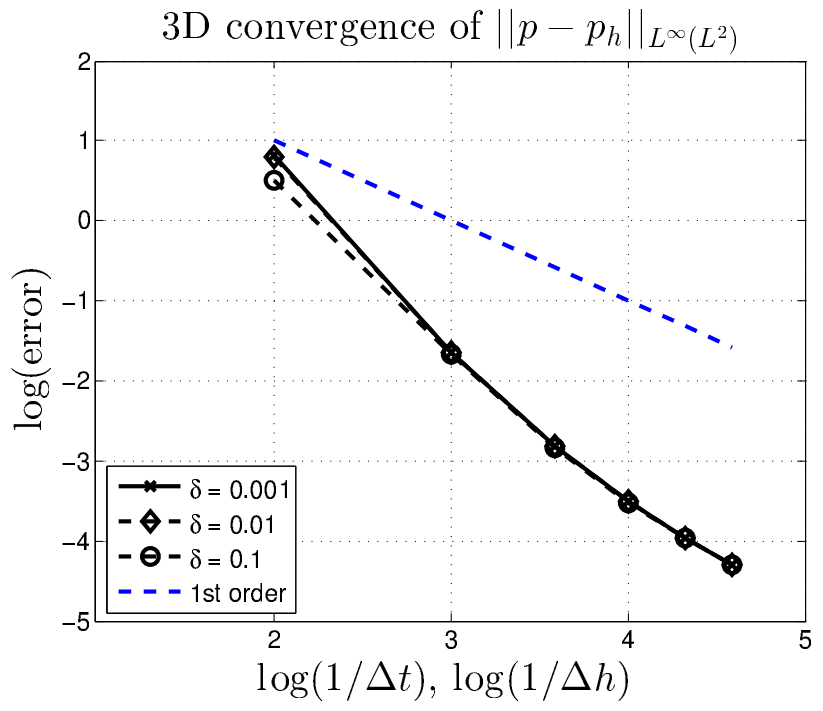

(c)

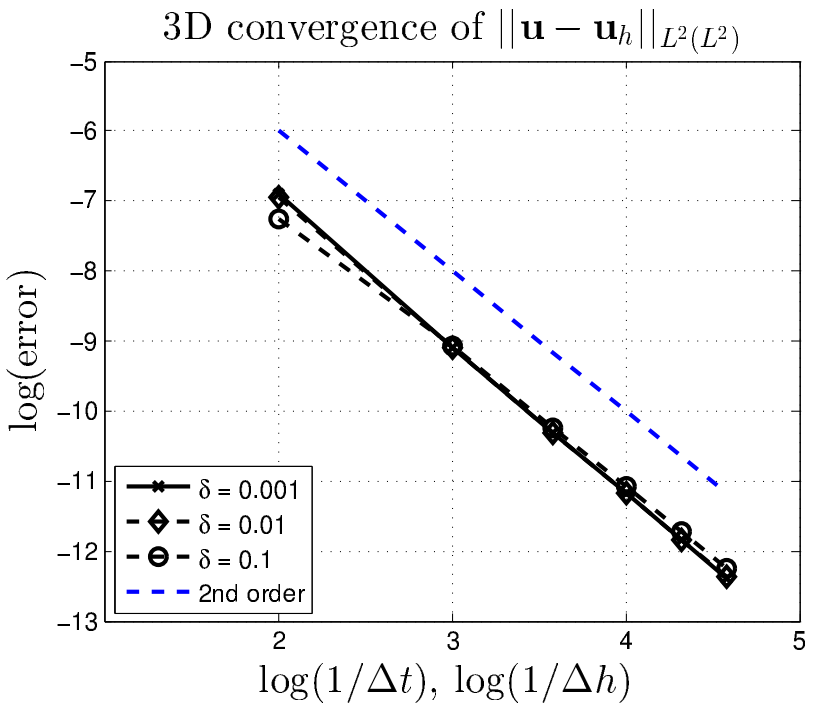

(d)

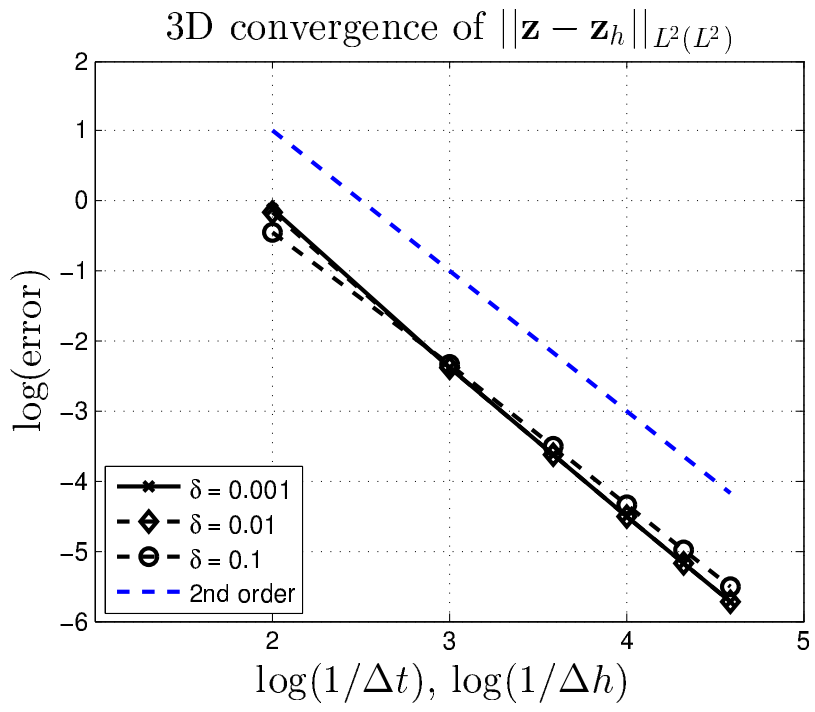

(e)

Figure 5.7: Convergence of the displacement, fluid flux, and pressure errors in their respective norms of the simplified poroelastic 3D test problem with different (stable) values for the stabilisation parameter $\delta$.

## 5.5 2D cantilever bracket problem

We consider the 2D cantilever bracket problem used in Phillips and Wheeler (2009) to illustrate the problem of spurious pressure oscillation. This problem was also used in Liu (2004) and Yi (2013) to demonstrate the ability of their method to overcome these spurious pressure oscillations. The cantilever bracket problem (shown in Figure 5.8a) is solved on a unit square $[0,1]^2$. No-flow flux boundary conditions are applied along all sides, the deformation is fixed ($\boldsymbol{u} = 0$) along the left hand-side ($x = 0$), and a downward traction force, $\boldsymbol{t}_N \cdot \boldsymbol{n} = -1$, is applied along the top edge ($y = 1$). The right and bottom sides are traction-free. For this numerical experiment, we set $\Delta t = 0.001$, $h = 1/96$, $\delta = 5 \times 10^{-6}$. The material parameters $\lambda$ and $\mu$ are chosen such that Youngs's modulus, $E = 10^5$ and Poisson's ratio $\nu = 0.4$ and $\alpha = 0.93, c_0 = 0, \boldsymbol{k} = 1 \times 10^{-7}\boldsymbol{I}$, values shown in Phillips and Wheeler (2009) to typically cause locking. The proposed stabilised finite element method yields a smooth pressure solution without any oscillations as is shown in Figure 5.8b.

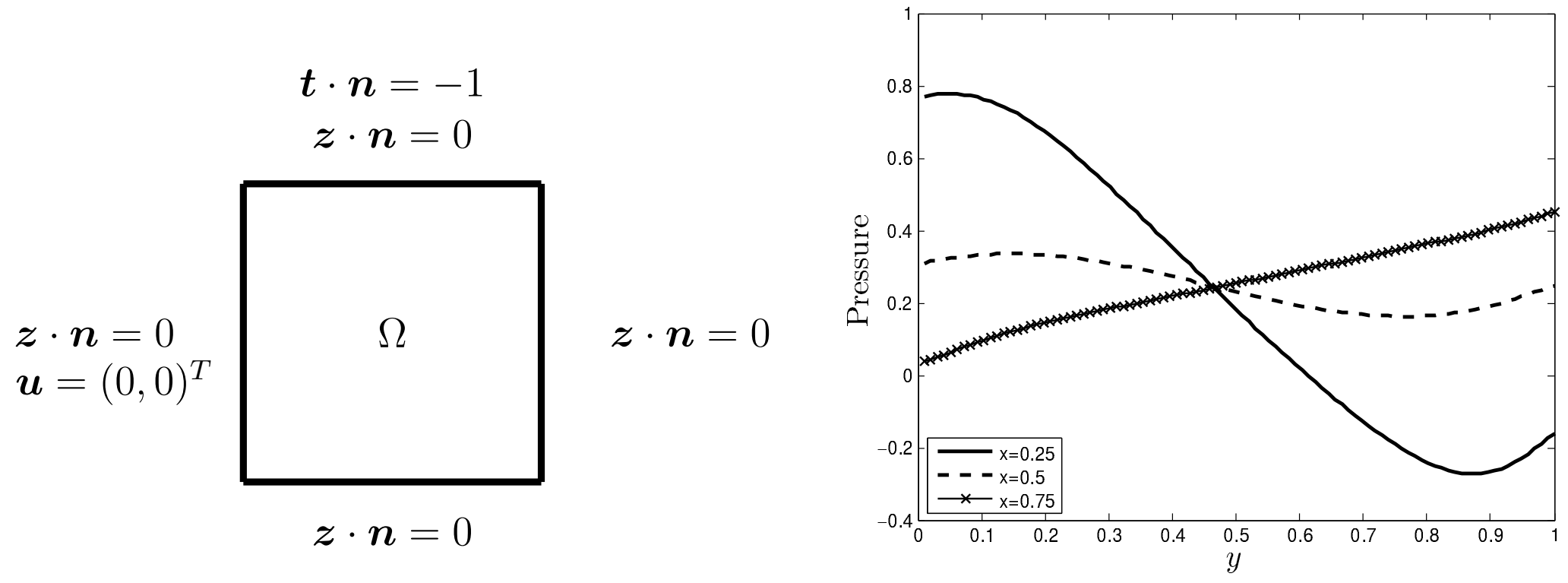

(a) (b)

Figure 5.8: (a) Boundary conditions for the cantilever bracket problem. (b) Pressure solution of the cantilever bracket problem at $t = 0.005$.

## 5.6 3D unconfined compression stress relaxation

In this test, a cylindrical specimen of porous tissue is exposed to a prescribed displacement in the axial direction while left free to expand radially, see Figure 5.9a. (Note that the two plates are not explicitly modelled in the simulation, but are realised through displacement boundary conditions.) After the initial loading, the displacement is held constant while the tissue relaxes in the radial direction due to interstitial fluid flow through the radial boundary. The outer radius and height of the cylinder is $1mm$, whereas the axial compression is $\epsilon_0 = 0.05mm$. The bottom of the tissue is constrained in the vertical direction. The fluid pressure is set to zero at the outer radial surface. The parameters used for the simulation can be found in Table 5.1. The material parameters $\mu_s$ and $\lambda$ can be related to the more familiar Young's modulus $E$ and the Poisson ratio $\nu$ by $\mu_s = \frac{E}{2(1+\nu)}$ and $\lambda = \frac{E\nu}{(1+\nu)(1-2\nu)}$. For the special case of a cylindrical geometry, Armstrong et al. (1984) found a closed-form analytical solution for the radial displacement $u$ on the porous medium, given by

$$\frac{u}{a}(a,t) = \epsilon_0 \left[ \nu + (1-2\nu)(1-\nu) \sum_{n=1}^{\infty} \frac{\exp\left(-\alpha_n^2 \frac{Mkt}{a^2}\right)}{\alpha_n^2 (1-\nu)^2 - (1-\nu)} \right], \tag{5.3}$$

where $\alpha_n$ are the solutions to the characteristic equation, given by $J_1(x) - (1-\nu)xJ_0(x)/(1-2\nu) = 0$, where $J_0$ and $J_1$ are Bessel functions, $\epsilon_0$ is the amplitude of the applied axial strain, $a$ is the radius of the cylinder, and $t_g$ is the characteristic time of diffusion (relaxation) given by $t_g = a^2/Mk$, where $M = \lambda + 2\mu$ is the P-wave modulus of the elastic solid skeleton, and $k$ is the permeability. Figure 5.9b shows the pressure solution after 5 seconds. The normalised radial displacement predicted by our implementation (Figure 5.10) using a value of $\delta = 0.001$ gives a root mean squared error of $6.7 \times 10^{-4}$ against the analytical solution provided by Armstrong et al. (1984), and yields a stable solution. The same test problem has also been used to verify other poroelastic software such as FEBio (Maas et al., 2012). The analytical solution available for this test problem describes the displacement of the outer radius which is directly dependent on the amount of mass in the system since the porous medium is assumed to be incompressible and fully saturated. It is therefore an ideal test problem for analyzing the effect that the added stabilisation term has on the conservation of mass. In Figure 5.10 we can see that for large values of $\delta$ the numerical solution loses mass faster and comes to a steady state that has less mass than the analytical solution. This is a clear limitation of the method and the stability parameter therefore needs to be chosen carefully. However, for 3D problems $\delta$ can be chosen to be very small so this effect is negligible, as can be seen in Figure 5.10 for a stable value of $\delta = 0.001$.

| Parameter | Description | Value |
|---|---|---|
| $k$ | Dynamic permeability | $10^{-1}\ \mathrm{m^3\,s\,kg^{-1}}$ |
| $\nu$ | Poisson ratio | 0.15 |
| $E$ | Young's modulus | $1000\ \mathrm{kg\,m^{-1}\,s^{-2}}$ |
| $\Delta t$ | Time step used in the simulation | $0.1\,\mathrm{s}$ |
| $T$ | Final time of the simulation | $10\,\mathrm{s}$ |

Table 5.1: Parameters used for the unconfined compression test problem.

## 5.7 Conclusion

We have presented numerical experiments in 2D and 3D that illustrate the convergence of the method, the effectiveness of the method in overcoming spurious pressure oscillations, and the added mass effect of the stabilisation term.

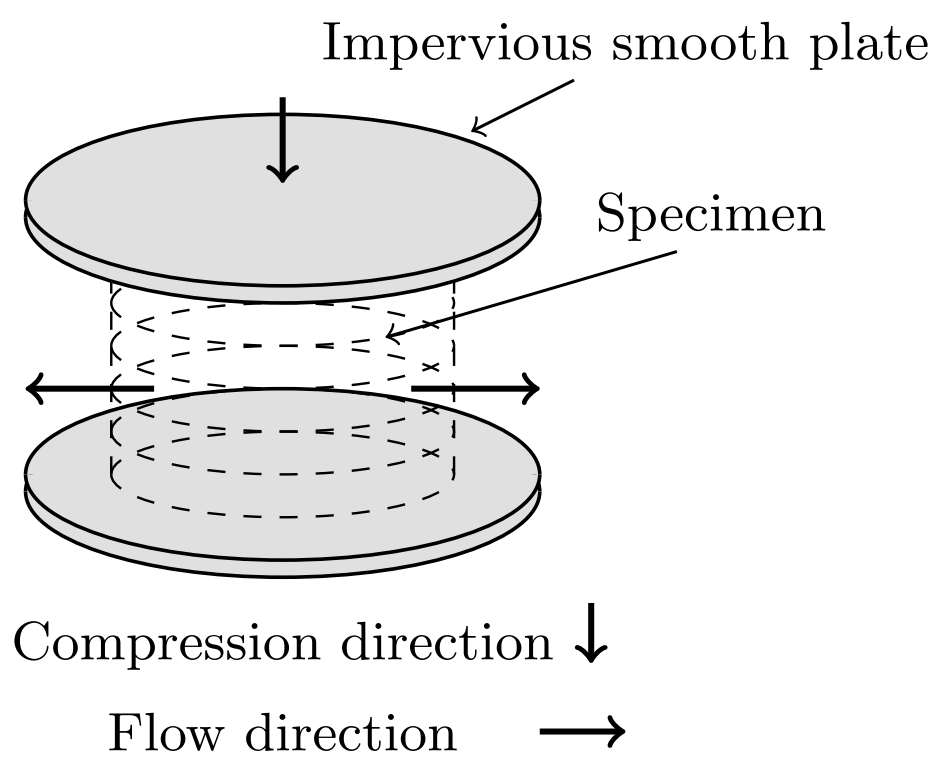

(a)

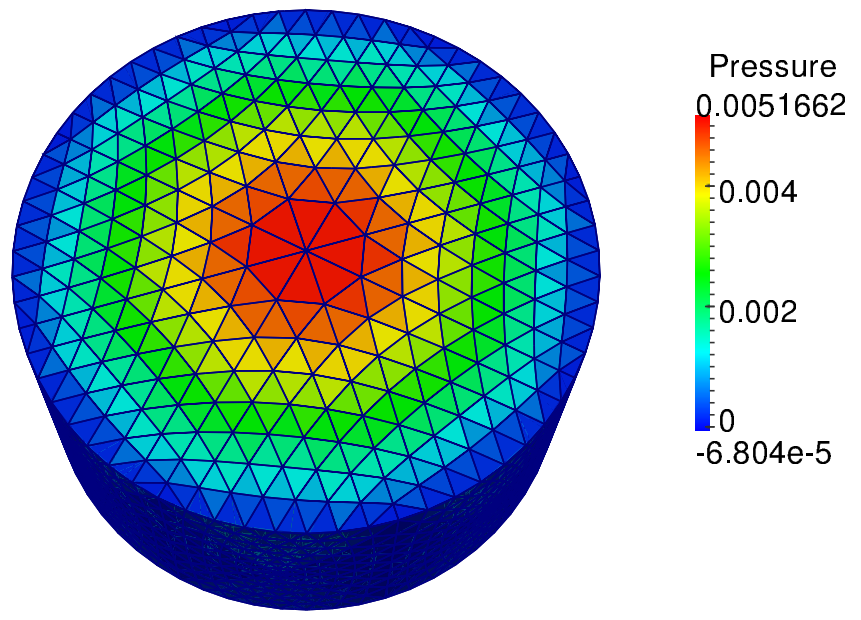

(b)

Figure 5.9: (a) Sketch of the test problem. The porous medium is being compressed between two smooth impervious plates. The frictionless plates permit the porous medium to expand in order to conserve volume and then to gradually relax as the fluid seeps out radially. (b) Pressure field solution at $t = 5s$, using a mesh with 28160 tetrahedra.

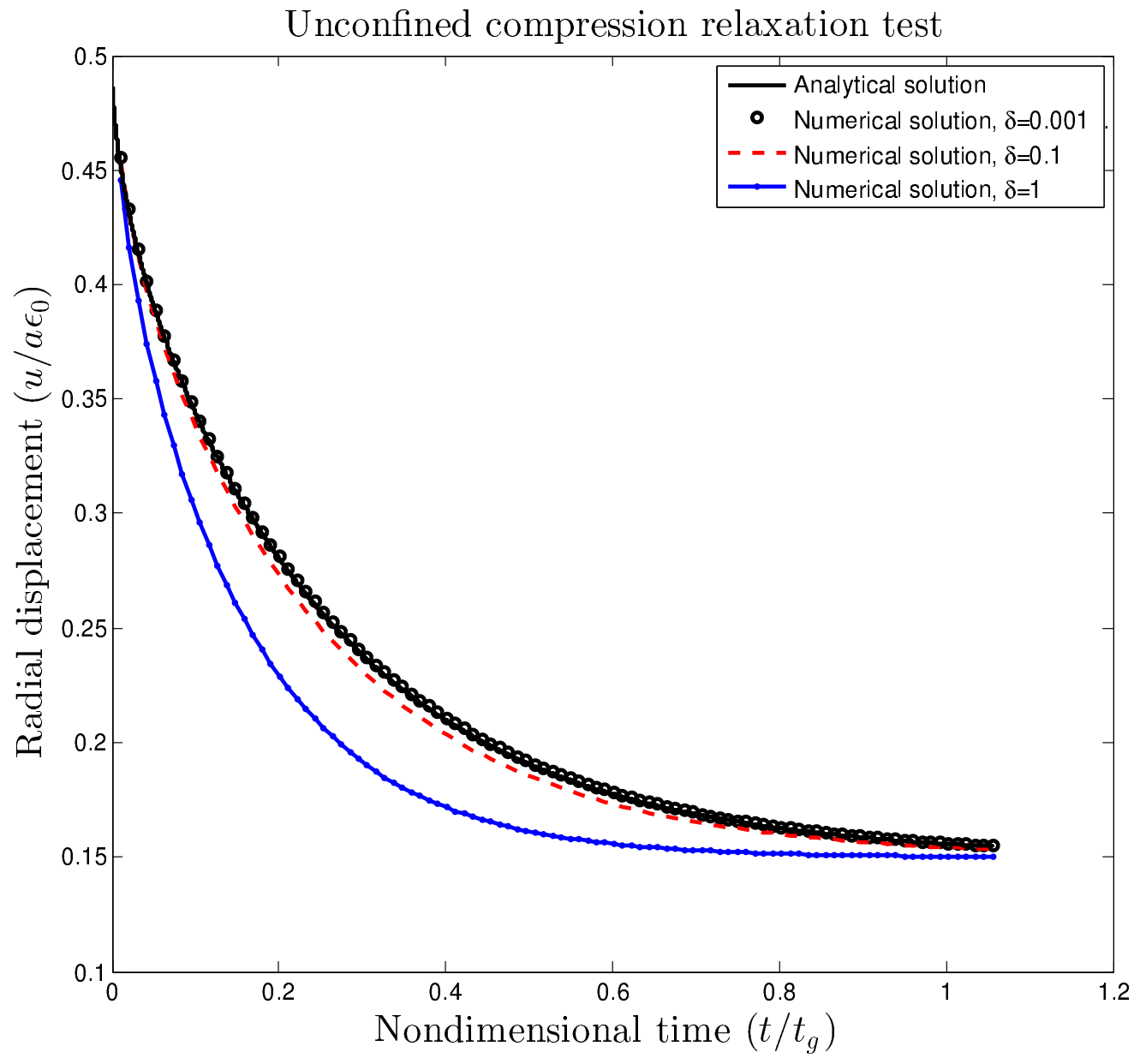

Figure 5.10: Normalised radial displacement versus normalised time calculated using the analytical solution, and using the proposed numerical method with different values of $\delta$. At $t = 0$ the radial expansion is half of the axial compression indicating the instantaneous incompressibility of the poroelastic tissue. The final amount of tissue recoil depends on the intrinsic Poisson ratio of the tissue skeleton.

# Chapter 6

# A stabilised finite element method for poroelasticity valid in large deformations

## 6.1 Introduction

In Chapter 4, we developed a stabilised, low-order, mixed finite-element method for the fully saturated, incompressible, poroelasticity equations, in the linear, small deformation case. In this Chapter we extend this work to the nonlinear, large deformation case.

In section 6.2, we recall the large deformation quasi-static incompressible poroelastic model. In section 6.3 we present the stabilised nonlinear finite-element method, and provide some implementation details in section 6.4. In section 6.5, we present a range of 3D numerical experiments to verify the accuracy of the method and illustrate its ability to reliably capture steep pressure gradients.

## 6.2 The model

Following Ateshian et al. (2010) and Almeida and Spilker (1998), we recall the governing equations (2.38) for a fully saturated, incompressible poroelastic model valid in large deformations. The problem is to find $\boldsymbol{\chi}(\boldsymbol{X}, t)$, $\boldsymbol{z}(\boldsymbol{x}, t)$ and $p(\boldsymbol{x}, t)$ such that

$$
\begin{aligned}
-\nabla \cdot (\boldsymbol{\sigma}_e - p\boldsymbol{I}) &= \rho \boldsymbol{f} && \text{in } \Omega_t, \\
\boldsymbol{k}^{-1}\boldsymbol{z} + \nabla p &= \rho^f \boldsymbol{f} && \text{in } \Omega_t, \\
\nabla \cdot (\boldsymbol{\chi}_t + \boldsymbol{z}) &= g && \text{in } \Omega_t, \\
\boldsymbol{\chi}(\boldsymbol{X}, t)|_{\boldsymbol{X}=\boldsymbol{\chi}^{-1}(\boldsymbol{x},t)} &= \boldsymbol{X} + \boldsymbol{u}_D && \text{on } \Gamma_D, \\
(\boldsymbol{\sigma}_e - p\boldsymbol{I})\boldsymbol{n} &= \boldsymbol{t}_N && \text{on } \Gamma_N, \\
\boldsymbol{z} \cdot \boldsymbol{n} &= q_D && \text{on } \Gamma_F, \\
p &= p_D && \text{on } \Gamma_P, \\
\boldsymbol{\chi}(\boldsymbol{X}, 0) &= \boldsymbol{X} && \text{in } \Omega_0.
\end{aligned}
\tag{6.1}
$$

**Remark 6.2.1.** *It is a straightforward extension to include the solid inertia $\boldsymbol{a}^s$ which can then be discretised using a Newmark scheme, see e.g. Chapelle et al. (2010), Li et al. (2004), Sauter and Wieners (2010).*

## 6.3 The stabilised finite element method

For ease of presentation, we will assume all Dirichlet boundary conditions are homogeneous, ie., $\boldsymbol{u}_D = \boldsymbol{0}, q_D = 0, p_D = 0$.

### 6.3.1 Weak formulation

We define the following spaces for deformed location, fluid flux and pressure respectively,

$$\begin{aligned}
\boldsymbol{W}^E(\Omega_t) &= \{\boldsymbol{v} \in (H^1(\Omega_t))^d : \boldsymbol{v} = \boldsymbol{0} \text{ on } \Gamma_D\}, \\
\boldsymbol{W}^D(\Omega_t) &= \{\boldsymbol{w} \in H_{div}(\Omega_t) : \boldsymbol{w} \cdot \boldsymbol{n} = 0 \text{ on } \Gamma_F\}, \\
\mathcal{L}(\Omega_t) &= \left\{ \begin{array}{ll} L^2(\Omega_t) & \text{if } \Gamma_N \cup \Gamma_P \neq \emptyset \\ L^2_0(\Omega_t) & \text{if } \Gamma_N \cup \Gamma_P = \emptyset, \end{array} \right\},
\end{aligned}$$

where $L^2_0(\Omega_t) = \left\{ q \in L^2(\Omega_t) : \int_{\Omega_t} q \, \mathrm{d}x = 0 \right\}$.

We make use of the identity $\nabla \cdot (\boldsymbol{\sigma}_e \boldsymbol{v}) = \nabla \cdot \boldsymbol{\sigma}_e \cdot \boldsymbol{v} + \boldsymbol{\sigma}_e : \nabla \boldsymbol{v}$, and the symmetry of $\boldsymbol{\sigma}_e$ to yield the following continuous weak problem. Find $\boldsymbol{\chi}(\boldsymbol{X}, t) \in \boldsymbol{W}^E(\Omega_t)$, $\boldsymbol{z}(\boldsymbol{x}, t) \in \boldsymbol{W}^D(\Omega_t)$ and $p(\boldsymbol{x}, t) \in \mathcal{L}(\Omega_t)$ for any time $t \in (0, T]$ such that

$$\begin{aligned}
\int_{\Omega_t} \left[ \boldsymbol{\sigma}_e : \nabla^S \boldsymbol{v} - p \nabla \cdot \boldsymbol{v} \right] \, \mathrm{d}\Omega_t &= \int_{\Omega_t} \rho \boldsymbol{f} \cdot \boldsymbol{v} \, \mathrm{d}\Omega_t \\
&\quad + \int_{\Gamma_N} \boldsymbol{t}_N \cdot \boldsymbol{v} \, \mathrm{d}\Gamma_N \quad \forall \boldsymbol{v} \in \boldsymbol{W}^E(\Omega_t), \\
\int_{\Omega_t} \left[ \boldsymbol{k}^{-1} \boldsymbol{z} \cdot \boldsymbol{w} - p \nabla \cdot \boldsymbol{w} \right] \, \mathrm{d}\Omega_t &= \int_{\Omega_t} \rho^f \boldsymbol{f} \cdot \boldsymbol{w} \, \mathrm{d}\Omega_t \quad \forall \boldsymbol{w} \in \boldsymbol{W}^D(\Omega_t), \\
\int_{\Omega_t} \left[ q \nabla \cdot \boldsymbol{\chi}_t + q \nabla \cdot \boldsymbol{z} \right] \, \mathrm{d}\Omega_t &= \int_{\Omega_t} g q \, \mathrm{d}\Omega_t \quad \forall q \in \mathcal{L}(\Omega_t).
\end{aligned} \tag{6.2}$$

Here $\nabla^S \boldsymbol{v} = \frac{1}{2} \left( \nabla \boldsymbol{v} + (\nabla \boldsymbol{v})^T \right)$ for some vector $\boldsymbol{v}$.

### 6.3.2 The fully-discrete model

Let $\mathcal{T}^h$ be a partition of $\Omega_t$ into non-overlapping elements $K$, where $h$ denotes the size of the largest element in $\mathcal{T}^h$. We define the following finite element spaces,

$$\begin{aligned}
\boldsymbol{W}_h^E &= \left\{\boldsymbol{v}_h \in C^0(\Omega_t) : \boldsymbol{v}_h|_K \in P_1(K)\ \forall K \in \mathcal{T}^h, \boldsymbol{v}_h = \boldsymbol{0} \text{ on } \Gamma_D\right\},\\
\boldsymbol{W}_h^D &= \left\{\boldsymbol{w}_h \in C^0(\Omega_t) : \boldsymbol{w}_h|_K \in P_1(K)\ \forall K \in \mathcal{T}^h, \boldsymbol{w}_h \cdot \boldsymbol{n} = 0 \text{ on } \Gamma_F\right\},\\
Q_h &= \begin{cases} \left\{q_h : q_h|_K \in P_0(K)\ \forall K \in \mathcal{T}^h\right\} & \text{if } \Gamma_N \cup \Gamma_p \neq \emptyset \\ \left\{q_h : q_h|_K \in P_0(K), \int_{\Omega_t} q_h = 0\ \forall K \in \mathcal{T}^h\right\} & \text{if } \Gamma_N \cup \Gamma_p = \emptyset \end{cases},
\end{aligned}$$

where $P_0(K)$ and $P_1(K)$ are respectively the spaces of constant and linear polynomials on $K$.

We define the combined solution space $\mathcal{U}_h(t) = \boldsymbol{W}_h^E \times \boldsymbol{W}_h^D \times Q_h$. The discretisation in time is given by partitioning $[0,T]$ into $N$ evenly spaced non-overlapping regions $(t_{n-1}, t_n]$, $n = 1, 2, \ldots, N$ , where $t_n - t_{n-1} = \Delta t$. For any sufficiently smooth function $v(x,t)$ we define $v^n(x) = v(x, t_n)$ and the discrete time derivative by $v^n_{\Delta t} := \frac{v^n - v^{n-1}}{\Delta t}$. The fully-discrete weak problem is: For $n = 1, \ldots, N$, find $\boldsymbol{\chi}_h^n \in \boldsymbol{W}_h^E(\Omega_{t_n})$, $\boldsymbol{z}_h^n \in \boldsymbol{W}_h^D(\Omega_{t_n})$ and $p_h^n \in Q_h(\Omega_{t_n})$ such that

$$\begin{aligned}
\int_{\Omega_{t_n}} \left[\boldsymbol{\sigma}^n_{e,h} : \nabla^S \boldsymbol{v}_h - p_h^n \nabla \cdot \boldsymbol{v}_h\right] \,\mathrm{d}\Omega_{t_n} &= \int_{\Omega_{t_n}} \rho \boldsymbol{f}^n \cdot \boldsymbol{v}_h \,\mathrm{d}\Omega_{t_n} \\
&\quad + \int_{\Gamma_N} \boldsymbol{t}_N^n \cdot \boldsymbol{v}_h \,\mathrm{d}\Gamma_N \quad \forall \boldsymbol{v}_h \in \boldsymbol{W}_h^E(\Omega_{t_n}),\\
\int_{\Omega_{t_n}} \left[\boldsymbol{k}^{-1} \boldsymbol{z}_h^n \cdot \boldsymbol{w}_h - p_h^n \nabla \cdot \boldsymbol{w}_h\right] \,\mathrm{d}\Omega_{t_n} &= \int_{\Omega_{t_n}} \rho^f \boldsymbol{f}^n \cdot \boldsymbol{w}_h \,\mathrm{d}\Omega_{t_n} \quad \forall \boldsymbol{w}_h \in \boldsymbol{W}_h^D(\Omega_{t_n}),\\
\int_{\Omega_{t_n}} \left[q_h \nabla \cdot \boldsymbol{\chi}^n_{\Delta t,h} + q_h \nabla \cdot \boldsymbol{z}_h^n\right] \,\mathrm{d}\Omega_{t_n} + J(p^n_{\Delta t,h}, q_h) &= \int_{\Omega_{t_n}} g^n q_h \,\mathrm{d}\Omega_{t_n} \\
&\quad \forall q_h \in Q_h(\Omega_{t_n}). \qquad (6.3)
\end{aligned}$$

### 6.3.3 Solution via Newton iteration at $t_n$

Since the system of equations (6.3) is highly nonlinear, its solution requires a scheme such as Newton's method. With Newton's method, an improved solution is obtained from a linear approximation of the nonlinear equation at an already computed solution. This first order Taylor expansion corresponds in finite element applications to the linearisation of the weak form, and can be obtained by the directional derivative, explained in section 6.3.4.

Let $\mathfrak{u}_h^n = \{\boldsymbol{\chi}_h^n, \boldsymbol{z}_h^n, p_h^n\} \in \mathcal{U}_h(t_n)$ denote the solution vector at a particular time step, $\xi\mathfrak{u}_h = \{\xi\boldsymbol{\chi}_h, \xi\boldsymbol{z}_h, \xi p_h\}$ denote the solution increment vector, and $\mathfrak{v}_h = \{\boldsymbol{v}_h, \boldsymbol{w}_h, q_h\} \in \mathcal{V}_h(t)$ where $\mathcal{V}_h(t) = \boldsymbol{W}_{h0}^E \times \boldsymbol{W}_{h0}^D \times Q_h$. The nonlinear system of equations (6.3) can be recast in the form: Find $\mathfrak{u}_h^n \in \mathcal{U}_h(t_n)$ such that

$$G^n(\mathfrak{u}_h^n, \mathfrak{v}_h) = 0 \ \forall \mathfrak{v}_h \in \mathcal{V}_h(t_n), \tag{6.4}$$

where

$$\begin{aligned} G^n(\mathfrak{u}_h^n, \mathfrak{v}_h) = \int_{\Omega_{t_n}} \big(&\boldsymbol{\sigma}_{e,h}^n : \nabla^S \boldsymbol{v}_h - p_h^n \nabla \cdot \boldsymbol{v}_h + \boldsymbol{k}^{-1} \boldsymbol{z}_h^n \cdot \boldsymbol{w}_h - p_h^n \nabla \cdot \boldsymbol{w}_h \\ &+ q_h \nabla \cdot (\boldsymbol{\chi}_{\Delta t,h}^n + \boldsymbol{z}_h^n) - \rho \boldsymbol{f}^n \cdot \boldsymbol{v}_h + \rho^f \boldsymbol{f}^n \cdot \boldsymbol{w}_h + g q_h \big) \ \mathrm{d}\Omega_{t_n} \\ &- \int_{\Gamma_N} \boldsymbol{t}_N^n \cdot \boldsymbol{v}_h \ \mathrm{d}\Gamma_N. \end{aligned} \tag{6.5}$$

Given an approximate solution $\bar{\mathfrak{u}}_h^n$, we approximate (6.4) by

$$G^n(\bar{\mathfrak{u}}_h^n, \mathfrak{v}_h) + DG^n(\bar{\mathfrak{u}}_h^n, \mathfrak{v}_h)[\xi\mathfrak{u}_h] = 0 \ \forall \mathfrak{v}_h \in \mathcal{V}_h(t_n),$$

and solve

$$DG^n(\bar{\mathfrak{u}}_h^n, \mathfrak{v}_h)[\xi\mathfrak{u}_h] = -G(\bar{\mathfrak{u}}_h^n, \mathfrak{v}_h) \ \forall \mathfrak{v}_h \in \mathcal{V}_h(t_n), \tag{6.6}$$

for the Newton step $\xi\mathfrak{u}_h$, where $DG$ is the directional derivative of $G$, at $\bar{\mathfrak{u}}_h^n$, in the direction $\xi\mathfrak{u}_h$.

### 6.3.4 Approximation of $DG^n$.

In biphasic tissue problems, it is common to approximate directional derivative of $G$ by assuming the nonlinear elasticity term is the dominant nonlinearity and ignoring the other nonlinearities (Ün and Spilker, 2006; White and Borja, 2008). Let

$$E^n((\boldsymbol{\chi}_h^n, p_h^n), \boldsymbol{v}_h) = \int_{\Omega_{t_n}} \boldsymbol{\sigma}_{e,h}^n : \nabla^S \boldsymbol{v}_h - p_h^n \nabla \cdot \boldsymbol{v}_h \, \mathrm{d}\Omega_{t_n}. \tag{6.7}$$

For Newton's method we require the directional derivative of $E^n((\boldsymbol{\chi}_{\Delta t,h}^n, p_h^n), \boldsymbol{v}_h)$ at a particular trial solution $(\overline{\boldsymbol{\chi}_{\Delta t,h}^n}, \overline{p_h^n})$ in the direction $\xi\boldsymbol{\chi}_h$, given by (see Wriggers (2008, section 3.5.3))

$$DE^n((\overline{\boldsymbol{\chi}_h^n}, \overline{p_h^n}), \boldsymbol{v}_h)[\xi\boldsymbol{\chi}_h] = \int_{\overline{\Omega}_{t_n}} \nabla^S \boldsymbol{v}_h : \overline{\boldsymbol{\Theta}_h^n} : \nabla^S \xi\boldsymbol{\chi}_h + \overline{\boldsymbol{\sigma}_{e,h}^n} : \big((\nabla\xi\boldsymbol{\chi}_h)^T \cdot \nabla \boldsymbol{v}_h\big) \, \mathrm{d}\Omega_{t_n}, \tag{6.8}$$

where $\overline{\boldsymbol{\Theta}_h^n}$ is a fourth-order tensor and $\overline{\boldsymbol{\sigma}_{e,h}^n}$ is the effective (elastic) stress tensor, both evaluated at a trial solution $\overline{\boldsymbol{\chi}_h^n}$. Further, any variable with a bar above it will correspond to it being evaluated at a trial solution. The fourth-order spatial tangent modulus tensor $\boldsymbol{\Theta}$ is described in A.1. For a detailed explanation and derivation see Bonet and Wood (1997); Wriggers (2008). The approximate linearisation of the nonlinear problem is thus given by

$$\begin{aligned} DG^n(\overline{\mathfrak{u}}_h^n, \mathfrak{v}_h)[\xi\mathfrak{u}_h] \approx \\ \int_{\overline{\Omega}_{t_n}} \Big[\nabla^S \boldsymbol{v}_h : \overline{\boldsymbol{\Theta}_h^n} : \nabla^S \xi\boldsymbol{\chi}_h + \overline{\boldsymbol{\sigma}_{e,h}} : \big((\nabla\xi\boldsymbol{\chi}_h)^T \cdot \nabla \boldsymbol{v}_h\big) - \xi p_h \nabla \cdot \boldsymbol{v}_h \\ + \bar{\boldsymbol{k}}^{-1} \xi \boldsymbol{z}_h \cdot \boldsymbol{w}_h - \xi p_h \nabla \cdot \boldsymbol{w}_h + q_h \nabla \cdot \left(\frac{\xi\boldsymbol{\chi}_h}{\Delta t} + \xi \boldsymbol{z}_h\right)\Big] \, \mathrm{d}\Omega_{t_n}, \end{aligned} \tag{6.9}$$

Using (6.5), (6.9) and equation (6.6) the Newton solve becomes: Find $\xi\boldsymbol{\chi}_h \in \boldsymbol{W}_h^E(\Omega_{t_n})$, $\xi\boldsymbol{z}_h \in \boldsymbol{W}_h^D(\Omega_{t_n})$ and $\xi p_h \in Q_h(\Omega_{t_n})$ such that

$$
\begin{aligned}
&\int_{\overline{\Omega}_{t_n}} \left[\nabla^S \boldsymbol{v}_h : \overline{\boldsymbol{\Theta}_h^n} : \nabla^S \xi\boldsymbol{\chi}_h + \overline{\boldsymbol{\sigma}_{e,h}^n} : \left((\nabla\xi\boldsymbol{\chi}_h)^T \cdot \nabla\boldsymbol{v}_h\right) - \xi p_h \nabla\cdot\boldsymbol{v}_h\right] \mathrm{d}\overline{\Omega}_{t_n} \\
&\quad = -\int_{\overline{\Omega}_{t_n}} \left[\overline{\boldsymbol{\sigma}_{e,h}^n} : \nabla^S\boldsymbol{v}_h - \overline{p_h^n}\nabla\cdot\boldsymbol{v}_h - \overline{\rho}\boldsymbol{f}^n\cdot\boldsymbol{v}_h\right] \mathrm{d}\overline{\Omega}_{t_n} \\
&\qquad - \int_{\overline{\Gamma}_N} \boldsymbol{t}_N^n \cdot \boldsymbol{v}_h \,\mathrm{d}\overline{\Gamma}_N \quad \forall \boldsymbol{v}_h \in \boldsymbol{W}_h^E(\Omega_{t_n}), \\
&\int_{\overline{\Omega}_{t_n}} \left[\bar{\boldsymbol{k}}^{-1}\xi\boldsymbol{z}_h\cdot\boldsymbol{w}_h - \xi p_h\nabla\cdot\boldsymbol{w}_h\right] \mathrm{d}\overline{\Omega}_{t_n} \\
&\quad = -\int_{\overline{\Omega}_{t_n}} \left[\bar{\boldsymbol{k}}^{-1}\overline{\boldsymbol{z}_h^n}\cdot\boldsymbol{w}_h - \overline{p_h^n}\cdot\nabla\boldsymbol{w}_h - \overline{\rho^f}\boldsymbol{f}^n\cdot\boldsymbol{w}_h\right] \mathrm{d}\overline{\Omega}_{t_n} \qquad \forall\boldsymbol{w}_h \in \boldsymbol{W}_h^D(\Omega_{t_n}), \\
&\int_{\overline{\Omega}_{t_n}} \left[q_h\nabla\cdot\left(\frac{\xi\boldsymbol{\chi}_h}{\Delta t} + \xi\boldsymbol{z}_h\right)\right] \mathrm{d}\overline{\Omega}_{t_n} + J\left(\frac{\xi p_h}{\Delta t}, q_h\right) \\
&\quad = -\int_{\overline{\Omega}_{t_n}} \left[q_h\nabla\cdot(\overline{\boldsymbol{\chi}_{\Delta t,h}^n} + \overline{\boldsymbol{z}_h}) - g q_h\right] \mathrm{d}\overline{\Omega}_{t_n} + J\left(\overline{p_{\Delta t,h}}, q_h\right) \quad \forall q_h \in Q_h(\Omega_{t_n}).
\end{aligned}
\tag{6.10}
$$

## 6.4 Implementation details

### 6.4.1 Newton algorithm

We will now let $\mathfrak{u}_i^n := \{\boldsymbol{\chi}_i^n, \boldsymbol{z}_i^n, p_i^n\}$ denote the fully-discrete solution at the $i$th step within the Newton method at time $t^n$. To ease the notation, we have suppressed the lower case $h$, previously used to denote the spatial discretisation. To solve the nonlinear poroelastic problem using Newton's method at a particular time step, we perform the steps outlined in Figure 6.1. Further details about computational considerations are given in B.3.

$i = 0$

$\mathfrak{u}_0^n = \{\boldsymbol{\chi}^{n-1}, \boldsymbol{z}^{n-1}, p^{n-1}\}$

**for** $i = 1$ to $i \leq$ ITEMAX **do**

Assemble $\boldsymbol{R}(\mathfrak{u}_i^n, \mathfrak{u}^{n-1})$ and $\boldsymbol{K}(\mathfrak{u}_i^n)$ on $\Omega(t_n)_i$

Solve $\boldsymbol{K}(\mathfrak{u}_i^n)\delta\mathfrak{u}_{i+1}^n = -\boldsymbol{R}(\mathfrak{u}_i^n, \mathfrak{u}^{n-1})$

Compute $\mathfrak{u}_{i+1}^n = \mathfrak{u}_i^n + \delta\mathfrak{u}_{i+1}^n$

Update the mesh, $(\Omega_{t_n})_{i+1} = \boldsymbol{\chi}_i^n$

**if** $||\boldsymbol{R}(\mathfrak{u}_i^n, \mathfrak{u}^{n-1})|| <$ TOL & $||\mathfrak{u}_i^n - \mathfrak{u}_{i-1}^n|| <$ TOL **then**

Newton iteration has converged.

Break out of for loop.

**end if**

**end for**

Figure 6.1: Newton algorithm at $t_n$. The norms here are L2 norms scaled relative to the size of the solution.

At each Newton iteration we are required to solve the linear system

$$\boldsymbol{K}(\mathfrak{u}_i^n)\xi\mathfrak{u}_{i+1}^n = -\boldsymbol{R}(\mathfrak{u}_i^n, \mathfrak{u}^{n-1}), \tag{6.11}$$

where $\boldsymbol{K}(\mathfrak{u}_i^n)$ and $\boldsymbol{R}(\mathfrak{u}_i^n, \mathfrak{u}^{n-1})$ are the matrix and vector representations of $DG(\mathfrak{u}_i^n)$ and $G(\mathfrak{u}_i^n, \mathfrak{u}^{n-1})$, respectively. The system (6.11) can be expanded as

$$\begin{bmatrix} \boldsymbol{K}^e & 0 & \boldsymbol{B}^T \\ 0 & \Delta t\boldsymbol{M} & \Delta t\boldsymbol{B}^T \\ \boldsymbol{B} & \Delta t\boldsymbol{B} & -\boldsymbol{J} \end{bmatrix} \begin{bmatrix} \xi\boldsymbol{u}_{i+1}^n \\ \xi\boldsymbol{z}_{i+1}^n \\ \xi p_{i+1}^n \end{bmatrix} = - \begin{bmatrix} \boldsymbol{r}_1(\boldsymbol{\chi}_i^n, p_i^n) \\ \Delta t\boldsymbol{r}_2(\boldsymbol{\chi}_i^n, \boldsymbol{z}_i^n, p_i^n) \\ -\boldsymbol{r}_3(\boldsymbol{\chi}_i^n, \boldsymbol{\chi}^{n-1}, \boldsymbol{z}_i^n, p_i^n) \end{bmatrix}, \tag{6.12}$$

where the elements in the matrices in (6.12) are given by

$$
\begin{aligned}
\boldsymbol{k}^e_{kl} &= \int_{(\Omega_{t_n})_i} \boldsymbol{E}_k^T \boldsymbol{D}(\boldsymbol{\chi}^n_i)\boldsymbol{E}_l + (\nabla \boldsymbol{\phi}_k)^T \boldsymbol{\sigma}_e(\boldsymbol{\chi}^n_i)\nabla\boldsymbol{\phi}_l \,\mathrm{d}(\Omega_{t_n})_i, \\
\boldsymbol{m}_{kl} &= \int_{(\Omega_{t_n})_i} \boldsymbol{k}^{-1}(\boldsymbol{\chi}^n_i)\boldsymbol{\phi}_k \cdot \boldsymbol{\phi}_l \,\mathrm{d}(\Omega_{t_n})_i, \\
\boldsymbol{b}_{kl} &= -\int_{(\Omega_{t_n})_i} \psi_k \nabla \cdot \boldsymbol{\phi}_l \,\mathrm{d}(\Omega_{t_n})_i, \\
\boldsymbol{j}_{kl} &= \delta \sum_{K\in\mathcal{T}^h_i} \int_{\partial K \backslash \partial(\Omega_{t_n})_i} h_{\partial K}[\psi_k][\psi_l] \,\mathrm{d}s. \\
\boldsymbol{r}_{1l} &= \int_{(\Omega_{t_n})_i} (\boldsymbol{\sigma}_e(\boldsymbol{\chi}^n_i) - p^n_i \boldsymbol{I}) : \nabla\boldsymbol{\phi}_l - \rho(\boldsymbol{\chi}^n_i)\boldsymbol{\phi}_l \cdot \boldsymbol{f} \,\mathrm{d}(\Omega_{t_n})_i \\
&\qquad - \int_{(\Gamma_N)_i} \boldsymbol{\phi}_l \cdot \boldsymbol{t}_N(\boldsymbol{\chi}^n_i) \,\mathrm{d}(\Gamma_N)_i, \\
\boldsymbol{r}_{2l} &= \int_{(\Omega_{t_n})_i} \boldsymbol{k}^{-1}(\boldsymbol{\chi}^n_i)\boldsymbol{\phi}_l \cdot \boldsymbol{z}^n_i - p^n_i \nabla\cdot\boldsymbol{\phi}_l - \rho^f(\boldsymbol{\chi}^n_i)\boldsymbol{\phi}_l \cdot \boldsymbol{f} \,\mathrm{d}(\Omega_{t_n})_i, \\
\boldsymbol{r}_{3l} &= \int_{(\Omega_{t_n})_i} \psi_l \nabla\cdot\left(\boldsymbol{\chi}^n_i - \boldsymbol{\chi}^{n-1}\right) + \Delta t \psi_l \nabla \cdot \boldsymbol{z}^n_i - \Delta t \psi_l g \,\mathrm{d}(\Omega_{t_n})_i \\
&\qquad + \delta \sum_{K\in\mathcal{T}^h} \int_{\partial K\backslash\partial(\Omega_{t_n})_i} h_{\partial K}[\psi_l][p^n_i - p^{n-1}] \,\mathrm{d}s.
\end{aligned}
$$

Here $\boldsymbol{\phi}_k$ are vector valued linear basis functions such that the displacement vector at the $i$th iteration can be written as $\boldsymbol{\chi}^n_i = \sum_{k=1}^{n_\chi} \boldsymbol{\chi}^n_{i,k}\boldsymbol{\phi}_k$, with $\sum_{k=1}^{n_\chi} \boldsymbol{\chi}^n_{i,k}\boldsymbol{\phi}_k \in \boldsymbol{W}^E_h$. Similarly for the fluid flux vector we have $\boldsymbol{z}^n_i = \sum_{k=1}^{n_z} \boldsymbol{z}^n_{i,k}\boldsymbol{\phi}_k$, with $\sum_{k=1}^{n_z} \boldsymbol{z}^n_{i,k}\boldsymbol{\phi}_k \in \boldsymbol{W}^D_h$. The scalar valued constant basis functions $\psi_i$ are used to approximate the pressure, such that $\boldsymbol{p}^n_i = \sum_{k=1}^{n_p} p^n_{i,k}\psi_k$, with $\sum_{k=1}^{n_p} p^n_{i,k}\psi_k \in Q_h$. Also to aid the assembly of the fourth order tensor we have adopted the matrix voigt notation. In particular $\boldsymbol{D}$ is the matrix form of $\boldsymbol{\Theta}$, and $\boldsymbol{E}_k$ is the matrix version of $\nabla^S \boldsymbol{\phi}_k$, see equations (A.3) and (A.4) for details.

### 6.4.2 Fluid-flux boundary condition

When solving the equations for Darcy flow using the Raviart-Thomas element (RT-P0), the fluid-flux boundary condition is enforced naturally by this diver-gence free element. Unfortunately this is not possible using our proposed P1-P1-P0-stabilised element. However, solving the poroelastic equations (6.1) using a piecewise linear approximation for the deformation and Raviart-Thomas element for the fluid (P1-RT-P0) does not satisfy the discrete inf-sup condition and can yield spurious pressure oscillations, see Phillips and Wheeler (2008, 2009) for details.

To enforce the no-flux boundary condition $\boldsymbol{z}\cdot\boldsymbol{n} = q_D$ we introduce a Lagrange multiplier $\Lambda$ along the boundary $\Gamma_F$. Let $W^F = \{l \in H_{div}(\Gamma_F, \mathbb{R})\}$. The resulting modified continuous weak-form is now:

$$G((\boldsymbol{u},\boldsymbol{z},p)),(\boldsymbol{v},\boldsymbol{w},q)) + (\Lambda, \boldsymbol{w}\cdot\boldsymbol{n})_{\Gamma_F} = q_D \ \forall(\boldsymbol{v},\boldsymbol{w},q) \in \boldsymbol{W}^E(\Omega_t), \boldsymbol{W}^D(\Omega_t), \mathcal{L}(\Omega_t),$$
$$(\boldsymbol{z}\cdot\boldsymbol{n}, l)_{\Gamma_F} = q_D, \ \forall l \in W^F. \tag{6.13}$$

The discretisation and implementation of this additional constraint is straightforward and results in a linear system with additional degrees of freedom for every node on $\Gamma_F$. The terms $(\Lambda, \boldsymbol{w}\cdot\boldsymbol{n})_{\Gamma_F}$ and $(\boldsymbol{z}\cdot\boldsymbol{n}, l)_{\Gamma_F}$ are nonlinear since the normal is a function of the displacement. We have found that linearising these terms using a Picard type linearisation (lagging) does not affect the convergence of the Newton algorithm. Alternatively these terms could be linearised explicitly as has been described in detail for the traction boundary condition, see Wriggers (2008, section 4.2.5) and Ateshian et al. (2010).

## 6.5 Numerical results

We present three numerical examples to test the performance of the proposed stabilised finite element method. The first two examples are from mechanobiology and geotechnical applications and the third is a swelling example that undergoes significant large deformations. For the implementation we used the C++ library libMesh (Kirk et al., 2006), and the multi-frontal direct solver mumps (Amestoy et al., 2000) to solve the resulting linear systems. For the strain energy law we chose a Neo-Hookean law taken from Wriggers (2008, eqn. (3.119)), with the penalty term chosen such that $0 \leq \phi < 1$, namely

$$W(\boldsymbol{C}) = \frac{\mu}{2}(\mathrm{tr}(\boldsymbol{C}) - 3) + \frac{\lambda}{4}(J^2 - 1) - (\mu + \frac{\lambda}{2})\ln(J - 1 + \phi_0). \tag{6.14}$$

For further discussion on strain energy laws for porelasticity we refer to Chapelle and Moireau (2014) and Vuong et al. (2015). The material parameters $\mu$ and $\lambda$ in (6.14) can be related to the Young's modulus $E$ and the Poisson ratio $\nu$ by $\mu = E/(2(1+\nu))$ and $\lambda = (E\nu)/((1+\nu)(1-2\nu))$. Details of the effective stress tensor and fourth-order spatial tangent modulus for this particular law can be found in A.3. For the permeability law we chose

$$\boldsymbol{k}_0(\boldsymbol{C}) = k_0\boldsymbol{I}. \tag{6.15}$$

### 6.5.1 3D unconfined compression stress relaxation

This is the same test as previously described in section 5.6. The outer radius and height of the cylinder is $5mm$, whereas the axial compression is $0.01mm$. The parameters used for the simulation can be found in Table 6.1. The permeability has been chosen to be the same as in (Maas et al., 2012) where the same problem has also been used to test other large deformation poroelastic software such as FEBio (Maas et al., 2012). This permeability is comparable to the permeability used in section 5.6, which is one hundred times higher. This is because the duration of the simulation time in this example is one hundred times longer than in section 5.6. A smaller time step has also been used to allow for better convergence of the Newton method, see B.3.

| Parameter | Description | Value |
|---|---|---|
| $k$ | Dynamic permeability | $10^{-3}$ m$^3$ s kg$^{-1}$ |
| $\nu$ | Poisson's ratio | 0.15 |
| $E$ | Young's modulus | 1000 kg m$^{-1}$ s$^{-2}$ |
| $\Delta t$ | Time step used in the simulation | 4 s |
| $T$ | Final time of the simulation | 1000 s |
| $\delta$ | Stabilisation parameter | $10^{-3}$ |

Table 6.1: Parameters used for the unconfined compression test problem.

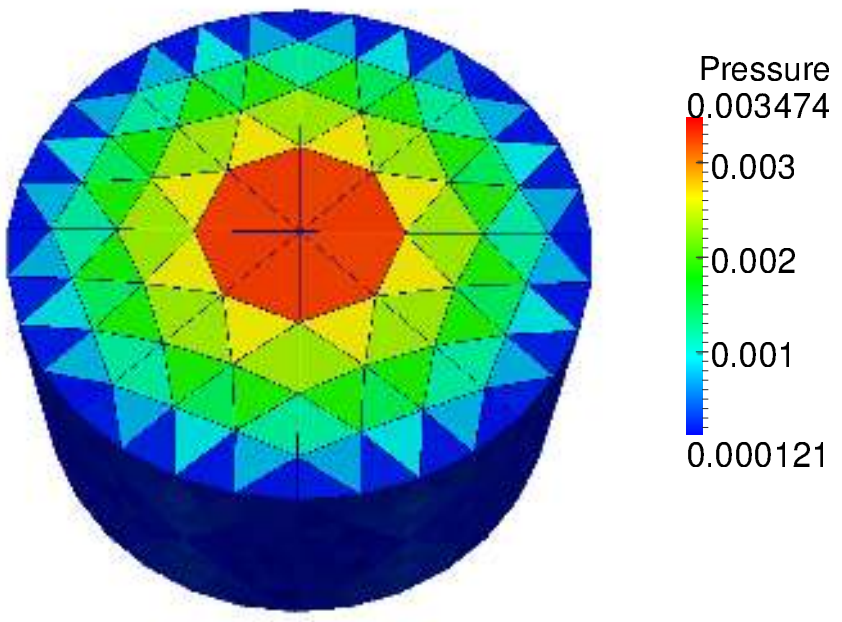

Figure 6.2: Pressure field at $t = 200s$ using a mesh with 3080 tetrahedra.

Figure 6.2 shows the pressure solution after 200 seconds. The computed radial displacement (Figure 6.3) shows good agreement with the analytical solution (5.3).

### 6.5.2 Terzaghi's problem

This is a classical geomechanics example with an analytical solution, and has been used to investigate finite element pressure oscillations, caused by overshooting of the numerical solution near the boundary (Murad and Loula, 1994; White and Borja, 2008). The domain consists of a porous column of unit height, bounded at the sides and bottom by rigid and impermeable walls. The top is free to drain ($p_D = 0$) and has a downward traction force, $p_0$, applied to it. The boundary

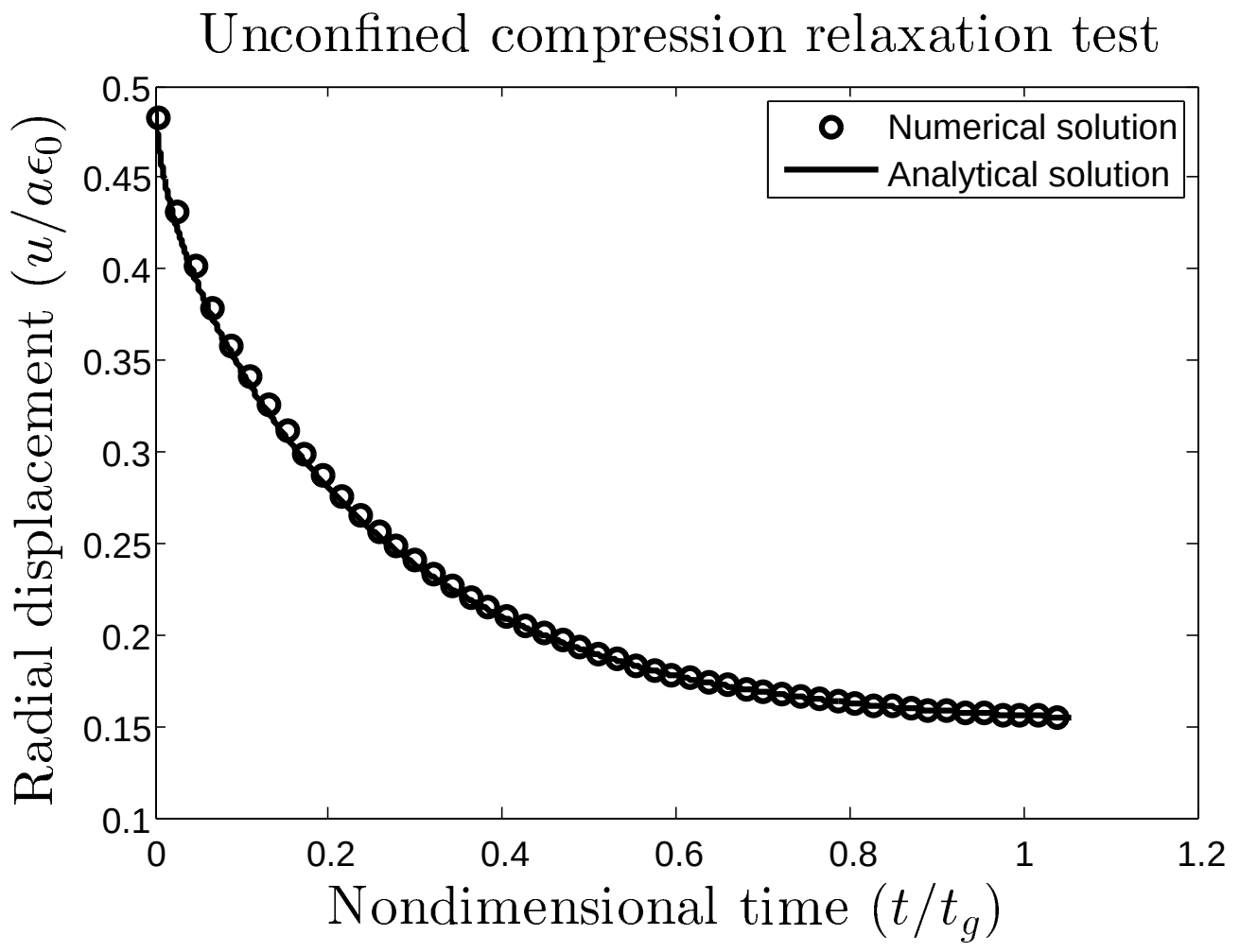

Figure 6.3: Radial expansion versus time comparing the analytical and numerical solutions with $\delta = 0.001$.

and initial conditions for this 1D problem can be written as

$$\begin{aligned} t_N = -p_0, \;\; p_D = 0 \;\; \text{on } x = 0, \\ u = 0, \;\; z = 0, \;\; \text{on } x = 1, \\ u = 0, \;\; z = 0, \;\; p = 0 \text{ in } (0,1). \end{aligned} \tag{6.16}$$

The analytical pressure solution, in nondimensional form is given by

$$p^* = \sum_n^\infty \frac{2}{\pi(n+1/2)} \sin(\pi(n+1/2)) \exp^{-\pi(n+1/2)(\lambda+2\mu)kt}. \tag{6.17}$$

When the poroelastic medium is subjected to the sudden loading, the saturating fluid undergoes an overpressurisation. Subsequently this overpressure progressively vanishes, owing to the diffusion process of the fluid towards the boundary at the top of the column, which remains drained. For a detailed explanation and derivation of this solution see Coussy (2004, section 5.2.2). We discretised the column using 60 hexahedral elements and solved the problem using the proposed stabilised low-order finite element method and a higher-order inf-sup stable fi-nite element method that uses a piecewise linear pressure approximation. The simulation results of the pressure for the two methods, taken at $t = 0.01s$ and $t = 1s$ are shown in Figure 6.4. The material parameters used for the simulation can be found in Table 6.2. At $t = 0.01s$ the piecewise linear (continuous) approximation suffers from overshooting due to the boundary layer solution (Figure 6.4a). The proposed method, which uses a piecewise constant pressure approximation does not suffer from this problem, and captures the pressure boundary layer solution reliably (Figure 6.4b). At $t = 1s$ the boundary layer has grown and both the piecewise linear (Figure 6.4c) and piecewise constant (Figure 6.4d) approximation yield satisfactory results.

| Parameter | Description | Value |
|---|---|---|
| $k_0$ | Dynamic permeability | $10^{-5}\,\mathrm{m^3\,s\,kg^{-1}}$ |
| $\nu$ | Poisson ratio | 0.25 |
| $E$ | Young's modulus | $100\,\mathrm{kg\,m^{-1}\,s^{-2}}$ |
| $\Delta t$ | Time step used in the simulation | $0.01\,\mathrm{s}$ |
| $T$ | Final time of the simulation | $1\,\mathrm{s}$ |
| $\delta$ | Stabilisation parameter | $2 \times 10^{-5}$ |

Table 6.2: Parameters used for Terzaghi's problem.

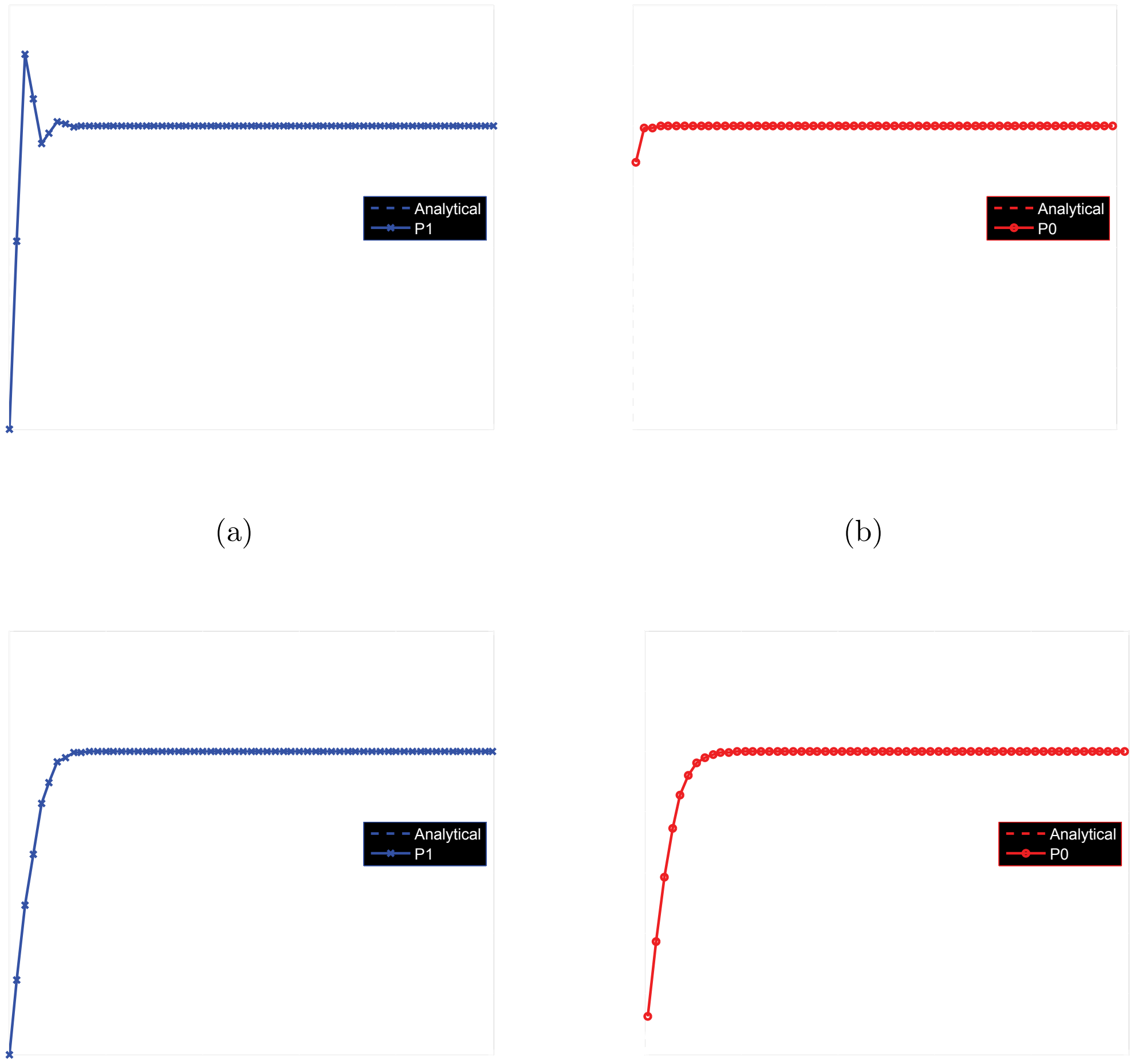

(a) (b) (c) (d)

Figure 6.4: (a) Pressure at $t = 0.01s$ using a continuous linear pressure approximation. (b) Pressure at $t = 0.01s$ using a discontinuous piecewise constant approximation. (c) Pressure at $t = 1s$ using a continuous linear pressure approximation. (d) Pressure at $t = 1s$ using a discontinuous piecewise constant approximation.

### 6.5.3 Swelling test

This problem is similar to the one in Chapelle et al. (2010) and highlights the method's ability to reliably capture jumps in the pressure solution due to changes in material parameters. Given a unit cube of material, a fluid pressure gradient is imposed between the two opposite faces at $X = 0$ and $X = 1$. The pressure $p_D$ on the inlet face $X = 0$ is increased very rapidly from zero to a limiting value of 10kPa, i.e., $p_D = 10^4(1 - \exp(-t^2/0.25))$ Pa). On the outlet face $X = 1$, the pressure is fixed to be zero, $p_D = 0$. There are no sources or sinks of fluid. A zero flux condition is applied for the fluid velocity on the four other faces ($Y = 0, 1,\ Z = 0, 1$). Normal displacements are required to be zero on the planes $X = 0$, $Y = 0$ and $Z = 0$. The permeability of the cube $0 < X < 0.5, 0.5 < Y < 1, 0 < Z < 0.5$ (1/8 of the volume of the unit cube) is smaller than in the rest of the domain by a factor of 500. The computational domain is shown in Figure 6.5a, highlighting the region of reduced permeability. The parameters chosen for this test problem are shown in Table 6.3.

Fluid enters the region from the inlet face and the material swells like a sponge, undergoing large deformation as shown in Figure 6.5b. The evolution of the pressure and the Jacobian at the points at $(0, 0, 1)$, $(0.5, 0, 1)$ and $(1, 0, 1)$ in the reference configuration are shown in Figures 6.6a and 6.6b respectively. These positions are indicated by the red, blue and green balls in Figure 6.5a. The pressure decreases roughly linearly with $x$, the increase in volume also follows a similar pattern. The pressure and volume change at the point $(0, 1, 0)$ (black ball in Figure 6.5a) is also shown in Figures 6.6a and 6.6b. Due to its reduced permeability this region is much slower to swell and achieve its ultimate equilibrium state and the fluid mainly flows around the area of reduced permeability, see Figure 6.5b. The steep pressure gradients at the boundary of the less permeable region seen in Figure 6.5b are well approximated by the piecewise constant (discontinuous) pressure space even on this relatively coarse discretisation. Continuous pressure spaces would require a much finer discretisation in this region.

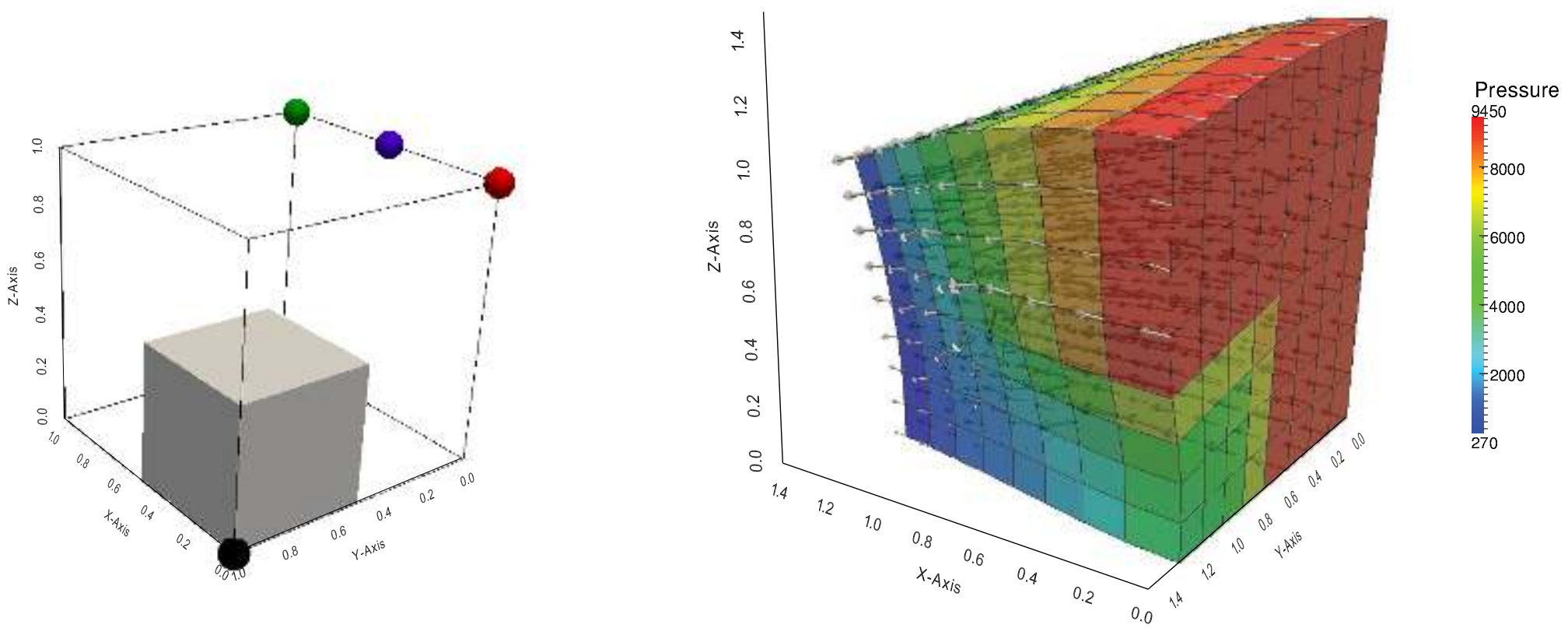

(a) (b)

Figure 6.5: (a) Initial simulation setup. The grey cube represents the area of reduced permeability. The colored balls highlight the position of the points used for tracking the pressure and volume change during the simulation, shown in Figures 6.6a and 6.6b. (b) The deformed cube after $1s$. The pressure solution is plotted and the jumps in pressure at the interface between the high and low permeability regions can clearly be seen. The arrows illustrate the fluid-flux profile.

| Parameter | Value |
|---|---|
| $k_0$ | $10^{-5}$ $\mathrm{m^3\,s\,kg^{-1}}$ |
| $\nu$ | 0.3 |
| $E$ | 8000 $\mathrm{kg\,m^{-1}\,s^{-2}}$ |
| $\Delta t$ | $0.02\,\mathrm{s}$ |
| $T$ | $20\,\mathrm{s}$ |
| $\delta$ | $10^{-4}$ |

Table 6.3: Parameters used for the swelling test problem.

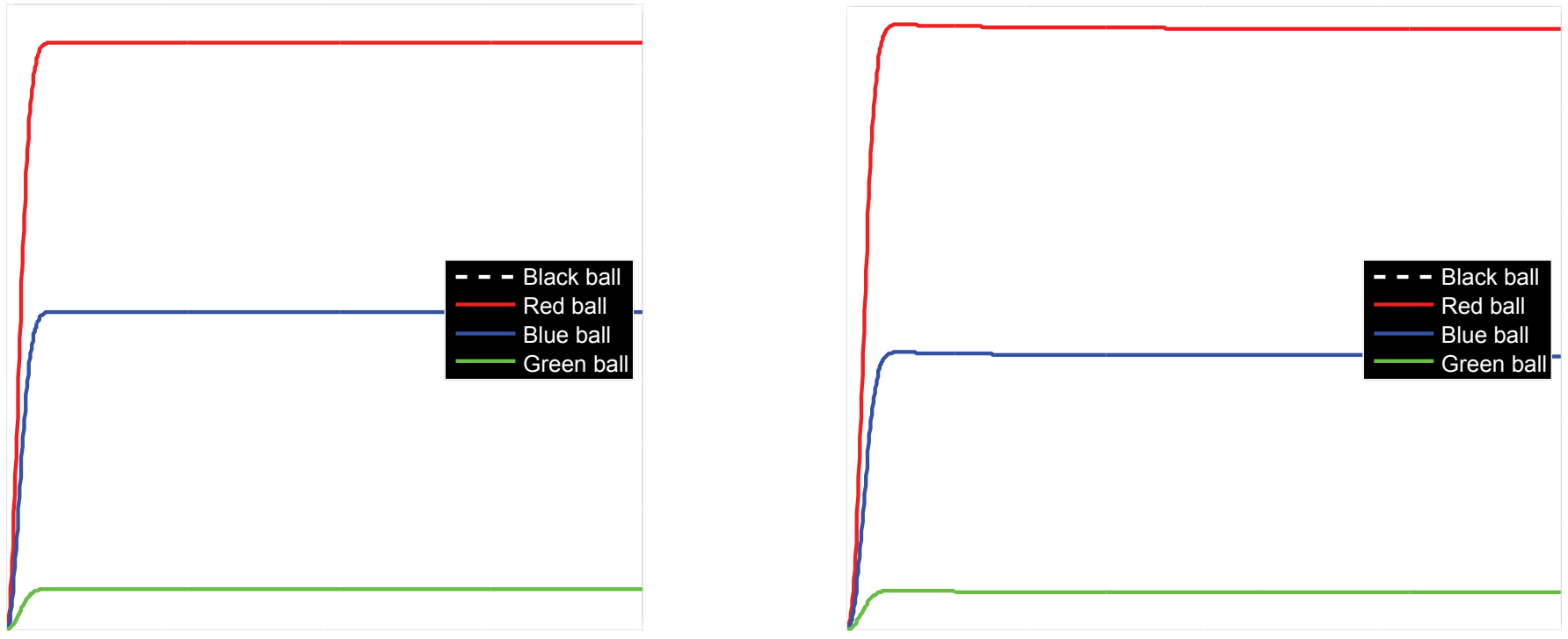

(a) (b)

Figure 6.6: Pressure (a) and volume change, $J$, (b) are plotted against time for four points, $(0, 0, 1)$ (red), $(0.5, 0, 1)$ (blue), $(1, 0, 1)$ (green), and $(1, 0, 1)$ (black) in the reference configuration. The position of these balls is also shown in Figure 6.5a.

## 6.6 Conclusion

The main contribution of this chapter has been to extend the local pressure jump stabilisation method (Burman and Hansbo, 2007), already applied to three-field linear poroelasticity in Chapter 4 to the large deformation case. Thus, the proposed scheme is built on an existing scheme, for which rigorous theoretical results about the stability and optimal convergence have been proven, and numerical experiments have confirmed its ability to overcome spurious pressure oscillations. Due to the discontinuous pressure approximation, sharp pressure gradients due to changes in material coefficients or boundary layer solutions can be captured reliably, circumventing the need for severe mesh refinement. Also, the addition of the stabilisation term introduces minimal additional computational work, can be assembled locally on each element using standard element information, and leads to a symmetric addition to the original system matrix, thus preserving any existing symmetry. As the numerical examples have demonstrated, the stabilisation scheme is robust and leads to high-quality solutions.

# Chapter 7

# A poroelastic-fluid-network model of the lung

The contents of this chapter closely follows the joint publication: L. Berger, R. Bordas, K. Burrowes, V. Grau, D. Kay, and S. Tavener; A poroelastic model coupled to a fluid network with applications in lung modelling *International Journal for Numerical Methods in Biomedical Engineering* 2015. L. Berger developed the coupling stratergy between the poroelastic medium and the fluid network, implemented the resulting algorithm, with guidance from D. Kay and R. Bordas, and wrote the original draft of the paper. The mesh of the lung lobes was provided by Materialise, and the airway tree mesh generated by R. Bordas, which was later pruned by L. Berger. The numerical tests were designed by L. Berger, R. Bordas and D. Kay, and were implemented by L. Berger. S. Tavener, K. Burrowes and V. Grau, assisted in improving the quality of the writing along with the other authors.

## 7.1 Introduction

The main function of the lungs is to exchange gas between air and blood, supplying oxygen during inspiration and removing carbon dioxide by subsequent expiration. Gas exchange is optimised by ensuring efficient matching between ventilation and blood flow, the distributions of which are largely governed by tissue deformation, gravity and branching structure of the airway and vascular trees. Understanding the interdependence between structure, and mechanical function in the lung has traditionally relied on direct measurement or medical imaging. Limitations of these approaches include difficulty in determining the contribution of specific subsystems to the function of the rest of the organ, making it hard to gain an indepth understanding of the underlying mechanics. A carefully constructed computational model provides the advantage of exact control over functional parameters and the geometry of the solution domain, allowing for investigations into complex functional mechanisms. The work developed in this chapter is part of a longer term aim to link detailed anatomic imaging to computational analysis of structure-function relationships in the integrated pulmonary system through computational modelling of the lung tissue and airway tree (Tawhai et al., 2006).

Previous work has typically focused on modelling either ventilation or tissue deformation in isolation. However evaluation of each component (i.e. tissue deformation and ventilation) separately does not necessarily give accurate ventilation predictions or provide a good indication of how the integrated organ works, this is because both components are interdependent. To gain a better understanding of the biomechanics in the lung it is therefore necessary to fully couple the tissue deformation with the ventilation. To achieve this tight coupling between the tissue deformation and the ventilation we propose a multiscale model that approximates the lung parenchyma by a biphasic (tissue and air, ig-noring blood) poroelastic model, that is then coupled to an airway fluid network model.

An integrated model of ventilation and tissue mechanics is particularly important for understanding respiratory diseases since nearly all pulmonary diseases lead to some abnormality of lung tissue mechanics (Suki and Bates, 2011). Chronic obstructive pulmonary disease (COPD) encompasses emphysema, the destruction of alveolar tissue, as well as chronic bronchitis, which can cause severe airway remodelling, bronchoconstriction and air trapping. All of the above affect tissue deformation since sections of lung are either not able to expand to inspire air, or to contract to release air. The effects of physiological changes occurring during disease, such as airway narrowing and changes in tissue properties, on regional ventilation and tissue stress are not well understood. For example, one hypothesis is that airway disease may precede emphysema (Galbán et al., 2012). An integrated model of ventilation and tissue mechanics can be used to investigate the impact of airway narrowing and tissue stiffness during obstructive lung diseases on tissue stresses, alveoli pressure and ventilation. Developing such a fully coupled model has to our knowledge not yet been achieved.

The rest of this chapter is organised as follows. In section 7.2 we give a brief overview of lung physiology, and in 7.3 we review the literature on computational ventilation models. In section 7.4 we outline the modelling assumptions, define the mathematical lung model in section 7.5, and describe its implementation in section 7.6. In section 7.7 we describe the generation of the computational lung geometry and boundary conditions, and in section 7.8 we present numerical simulations of tidal breathing, and investigate the effect of airway constriction and tissue weakening. Finally in section 7.9, we discuss limitations of the model, possible future directions and draw some conclusions.

## 7.2 Lung physiology

We will now give a basic review of lung physiology at an organ scale, focusing on the ventilation and tissue properties of the lung. A more complete introduction can be found in Cotes et al. (2009); West (2008).

### 7.2.1 Mechanics of breathing

During inspiration, the volume of the thoracic cavity increases and air is drawn into the lung by creating a sub-atmospheric pressure distribution. The increase in volume is brought about mainly by contraction of the diaphragm, which causes it to descend, and partly by the action of the intercostal muscles, which raise the ribs. The lung is elastic and subsequently returns passively to its preinspiratory volume during resting breathing. During expiration the intra-alveolar pressure becomes slightly higher than atmospheric pressure and gas flows out of the lungs (West, 2008).

### 7.2.2 Airway tree

The airway tree is divided into a conducting zone and a respiratory zone. Air passes through the upper respiratory tract to the trachea. From here the airway tree divides into right and left main bronchi, which in turn divide into lobar and then segmental bronchi. This process continues down to the terminal bronchioles, which are the smallest airways without alveoli. All of these bronchi make up the conducting airways. The terminal bronchioles, which appear at around generation 15-16, then continue to divide into respiratory bronchioles, which have occasional alveoli budding from their walls. Finally, we get to the alveolar ducts, which are completely lined with alveoli, see Figure 7.1a. This alveolated region of the lung where the gas exchange occurs is known as the respiratory zone (West, 2008). Table 7.1 documents the different flow characteristics found in the airway tree during slow and rapid breathing.

| Generation | Diameter cm | Length cm | Flow rate 10L/min Velocity (m/s) | Re | Flow rate 100L/min Velocity (m/s) | Re |
|---|---|---|---|---|---|---|
| Trachea | 1.80 | 12.0 | 65.8 | 775 | 658 | 7750 |
| 1 | 1.22 | 4.76 | 71.6 | 573 | 716 | 5730 |
| 5 | 0.35 | 1.07 | 53.6 | 123 | 536 | 1230 |
| 10 | 0.13 | 0.46 | 12.55 | 10.6 | 125 | 106 |
| 15 | 0.066 | 0.20 | 1.48 | 0.63 | 14.8 | 6.30 |
| 20 | 0.045 | 0.083 | 0.10 | 0.031 | 1.00 | 0.31 |

Table 7.1: Shows dimensions, velocity and the corresponding Reynolds number for different sections of the airway tree during slow and rapid breathing. These values have been taken from Pedley et al. (1970b).

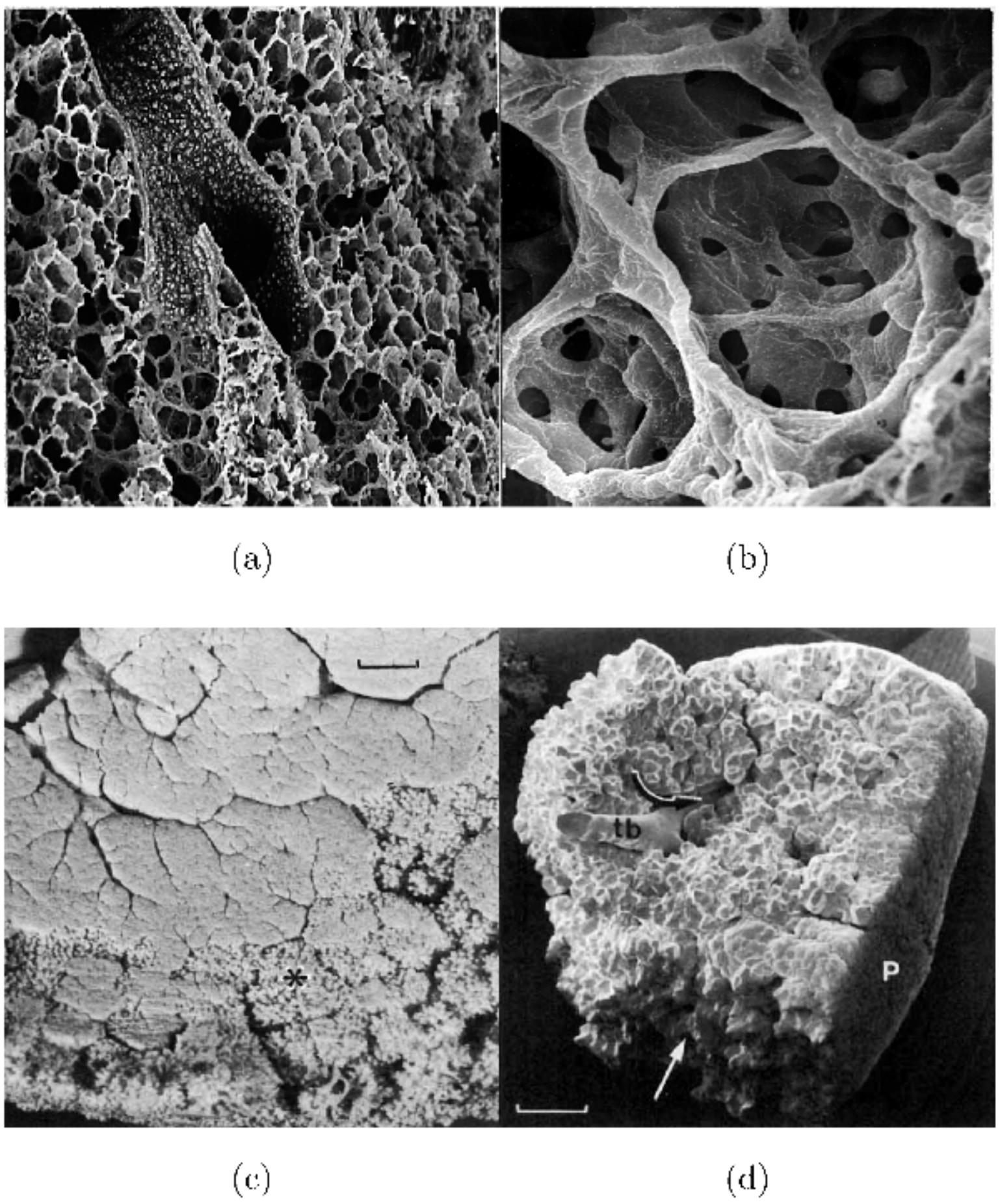

Figure 7.1: (a) Transition from terminal bronchiole to alveolar duct, from conducting airway to oxygen transfer area, diameter of terminal bronchiole is 0.5 mm. (b) A few alveoli in an alveolar duct. The dark round openings are pores between alveoli. The alveolar wall is quite thin and contains a network of capillaries. The average diameter of one alveoli is 0.2 mm. (c) Portions of silicone rubber casts of upper lobes of human lungs; asterisk marks incompletely filled regions. The outline of individual unfilled acinar units can also be seen. Scale marker, 5 mm. (d) Scanning electron micrograph of complete acinus with transitional bronchiole (tb) and surface abutting on pleura (P). Note the irregular surface where alveolar sacs of adjacent acini interdigitate (straight arrow). Scale marker, 1 mm. Images are reproduced from Lawrence Berkeley National Laboratory (1995).

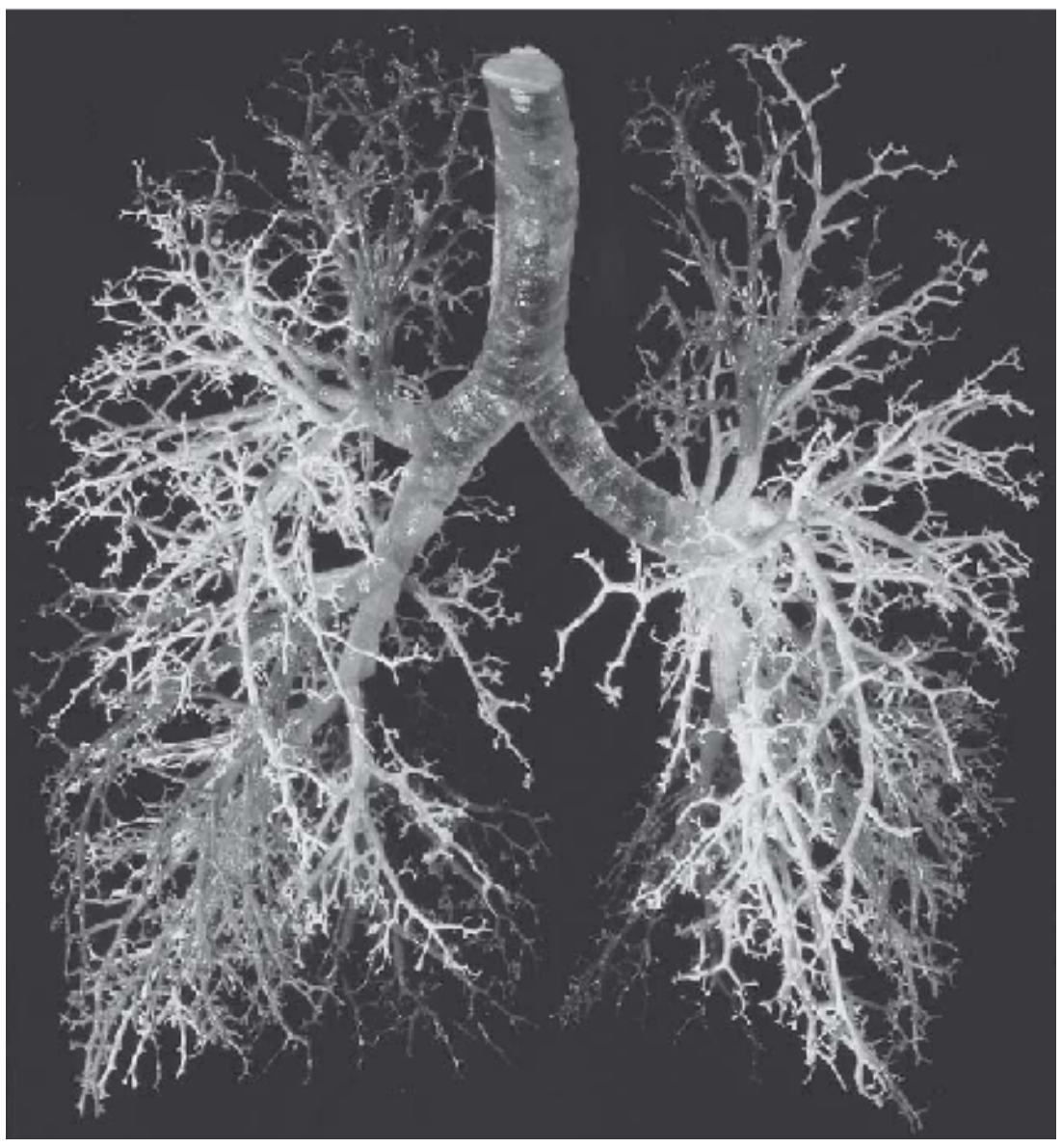

Figure 7.2: A rubber cast of the conducting airways of a human lung. The image is reproduced from West (2008).

### 7.2.3 Lung parenchyma

Lung parenchyma refers to the portion of the lung made up of the small air chambers (alveoli) participating in gas exchange. The alveoli are made up of collagen, elastin fibers and membranous structures containing the capillary network, see Figure 7.1b. Alveoli are arranged in sponge like structures and fill the entire volume of the lungs surrounding the conducting passages. Figure 7.1c shows a rubber cast of lung parenchyma, the dark lines outline the branching structure of the airways. The right and left lung are partitioned into three and two lobes, respectively. Lung segments of conic shape are then the first subdivision of these lobes. These structures are bounded by connective tissue such that surgical separation is often possible. In the right lung, there are usually ten segments whereas only nine can be found in the left lung. Within the segments, the bronchi branch about six to twelve times. The terminal bronchioles which appear after roughly $15 - 16$ branching generations then finally feed into approximately $30,000$ acini, see Figure 7.1d. These acini represent the largest lung units of which all airways are alveolated and thus participate in gas exchange (Weichert, 2011).

Also, the lung and lobes are surrounded by the pleura which is a membrane that folds back upon itself forming a double-layered structure between the lungs and the chest wall. The space in between the pleura is filled with fluid, allowing the lobe surfaces to slide over each other during the expansion and recoil of the chest wall, while maintaining the surface tension required to keep the lung in contact with the chest wall and thus inflated (Hedges, 2009).

### 7.2.4 The diseased lung

There exist numerous ways in which the mechanical function of the lung can be altered. In this section we will briefly describe pulmonary fibrosis, emphysema, and airway constriction.

Pulmonary fibrosis is a so-called restrictive disease. Here, abnormal deposition and organisation of connective proteins, particularly collagen, leaves lung tissue scarred and stiff with with compliance values decreasing to approximately 20% of normal values (Bates, 2009; Cotes et al., 2009).

Emphysema is characterised by an abnormal, permanent enlargement of air spaces distal to the terminal bronchioles and the destruction of their walls associated with loss of the elastic connective tissue. Large areas of lung tissue completely break down leaving big holes, see Figure 7.3. This results in a reduced area for gas exchange and a reduction in the elastic recoil of the lungs.

Airway constriction, which occurs in both asthma and COPD, changes airway resistance patterns. The level of airway resistance is sensitive to disease in the lungs. Narrowing of the airways can be caused through various mechanisms such as the airway inflammation or bronchoconstriction observed in asthma, mucous hyper-secretion and inflamed bronchi observed in chronic bronchitis, or the flaccid airways observed in emphysema (Hedges, 2009). This decrease in airway radius can significantly increase the resistance to flow.

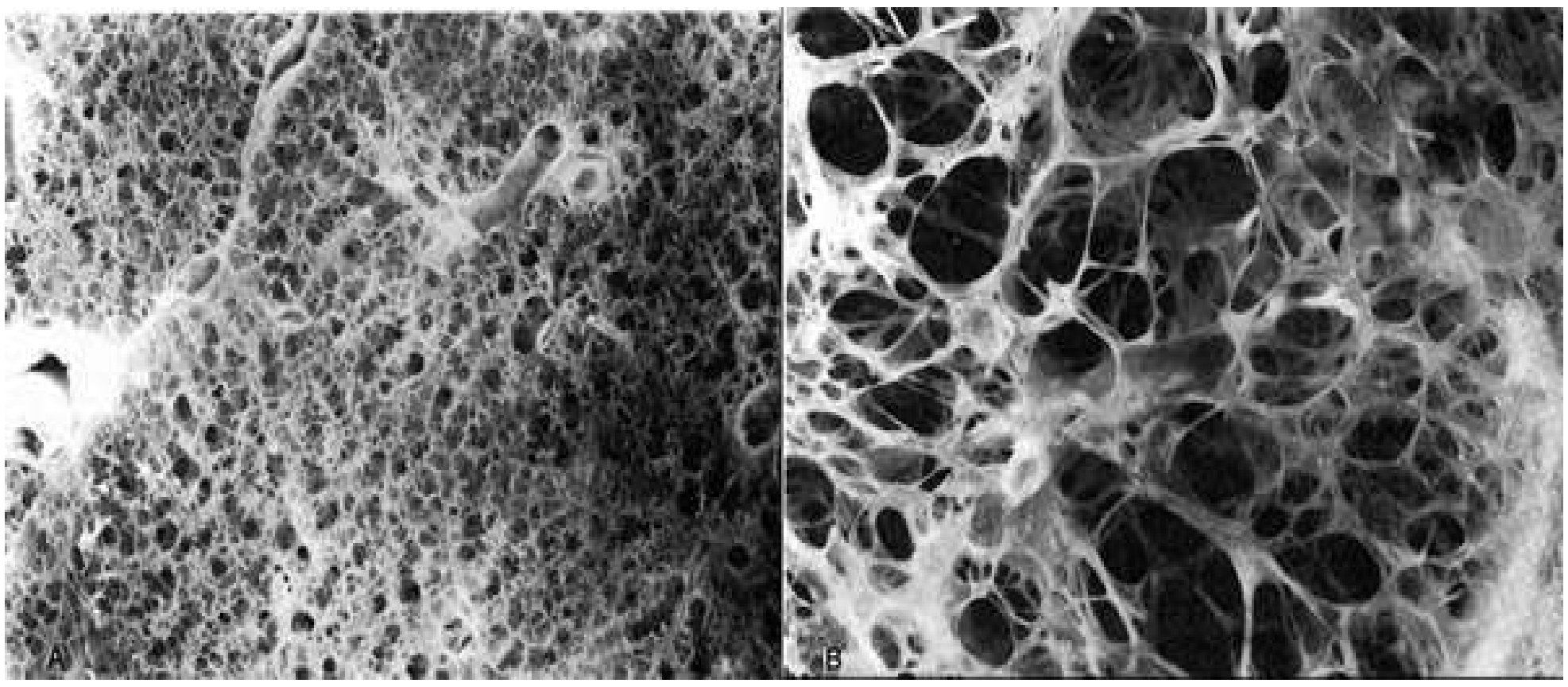

Figure 7.3: Left, a cross section of healthy parenchyma. Right, a cross section of diseased (emphysemic) lung parenchyma, with big holes appearing. Images are reproduced from G. Snell, ctsnet.org.

## 7.3 Computational lung models

There exist a large number of computational ventilation and deformation models for the lung. Some models are designed to model particular phenomena whilst others are more general. They also range in spatial complexity from 0D compartment type models to 3D models which are able to incorporate 'patient-specific' geometries extracted from CT images. In this brief review, we will focus on models that couple ventilation with tissue deformation and can be used as 'patient-specific' models. The term patient-specific is used very loosely here and only highlights that the geometry (computational mesh) used is extracted from an individual patient's CT scan. Unfortunatley a patient-specific model, that is able to produce clinically meaningful results, is currently not feasible. This would require more detailed information on the geometry of the lower airways and lobar segments, structure and elastic properties of the tissue, cardiac motion, and possibly cellular data. This information is currently not readily available and requires significant advances in experimental and imaging techniques. However a model, as presented in this work, that uses a basic geometry and has the capability to incorporate all the required data for a patient specific model is a good starting point for investigating general dynamics of the model, identifying key model parameters, and developing appropriate numerical methods that will allow a detailed patient specific model to be solved in the future.

One study that couples ventilation and tissue deformation using a one way coupling approach is described in Tawhai and Lin (2010). Here a mechanics model for lung tissue is used to provide flow boundary conditions at terminal branches for an airway model. This makes the resultant ventilation distribution dependent on the tissue deformation, for example due to gravity. Other sophisticated models of the whole lung that model ventilation and tissue deformation also exist (Ismail et al., 2013; Swan et al., 2012). Here the tissue is modelled by many independent elastic alveolar units. There is no clear way to conserve mass locally, so alveolar units can expand irrespectively of the size and position of neighbouring units. In reality the acini do not function as independent elastic balloons. They are physically coupled through fibrous scaffolding and shared alveolar septa. In our proposed model the tissue is modelled as one continuum, thus allowing us to conserve volume and couple neighbouring units. This is illustrated in Figure 7.4. Also, these lung models (Ismail et al., 2013; Swan et al., 2012) give information about the distribution of flow within the lung as a result of a pleural pressure boundary condition. However it is not possible to experimentally measure the pleural pressure in vivo using imaging or other apparatus. As part of the simulation protocol, the pleural pressure is therefore often tuned until physiological realistic flow rates are achieved. To overcome this issue, Yin et al. (2010, 2013) proposed to estimate the flow boundary conditions for full organ ventilation models by means of image registration. By solely relying on image registration to determine the ventilation distribution within the tissue one is not able to model the change in ventilation distribution due to progression of disease. We will build on Yin et al. (2010) by integrating image registration based boundary conditions within the proposed poroelastic model of lung deformation. In particular, we propose to register expiratory images to the inspiratory images, to yield an estimate of the deformation boundary condition for the lung surface, and drive the simulation through this deformation boundary condition. Thus the tissue deformation and subsequent flow boundary condition for tree branches inside the lung and ventilation distribution is not pre-determined, but calculated from the coupled poroelastic-airway-tree model.

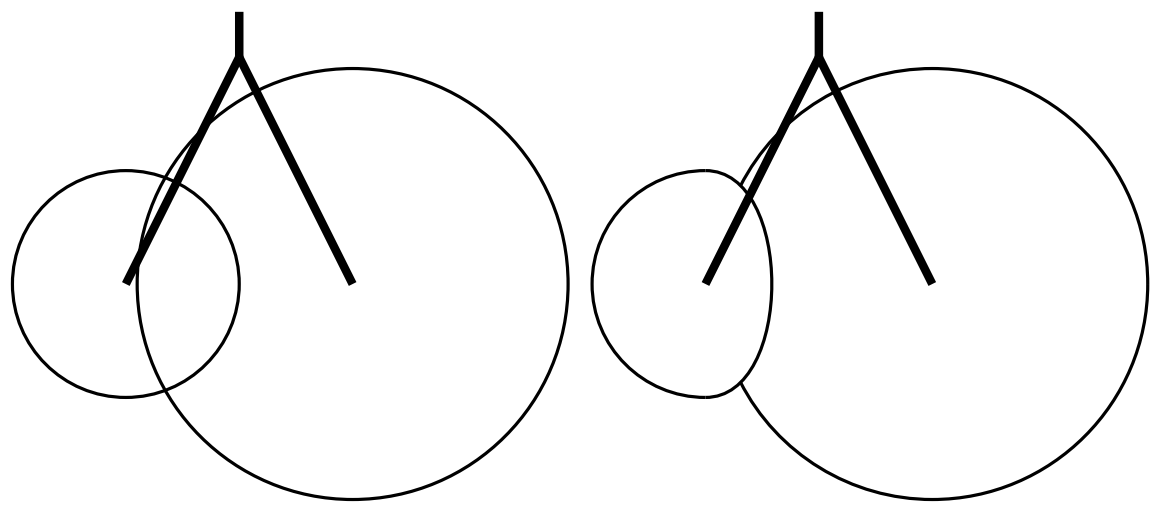

(a) (b)

Figure 7.4: Sketch of two balloon models where the right unit is more compliant, thus being able to expand more easily. (a) Balloon model with independent alveolar units. The overlap in the alveolar units illustrates that mass is not conserved. (b) Balloon model where the alveolar units are coupled. Here the inflation of each alveolar unit is compromised by the expansion of its neighbour.

## 7.4 Modelling assumptions

We will now give a brief commentary of the main modelling assumptions and how they might affect the proposed model's ability to predict deformation and ventilation within healthy and diseased lungs.

### 7.4.1 Approximating lung parenchyma using a poroelastic medium

**Averaging over the tissue:** One of the major assumptions is that we can approximate the lung parenchyma using a poroelastic continuum description. This makes our model computationally tractable and allows us to use the well established theory of poroelasticity to couple the air with the tissue.

The use of a continuum model can be further supported by looking at the different length scales and structures of the tissue. For the microscopic length scale denoted by $l$ of the parenchyma we will use the diameter of an alveolus that can be approximated to be 0.02 cm (Ochs et al., 2004). The macroscopic length scale $L$ can be taken to be the diameter of a segment which measures around 4 cm of tissue. So the ratio of the different length scales is small i.e $\epsilon := \frac{l}{L} \approx 0.005 \ll 1$. This along with the assumption that the structure of an acinus is porous (see Figures 7.1a and 7.1b) and approximately periodic supports the use of averaging techniques over the tissue to obtain a continuum description in the form of a poroelastic medium.

In Lewis and Owen (2001) a more rigorous approach has been used to derive macroscopic poroelastic equations for average air flows and tissue displacements in lung parenchyma using homogenisation theory. The resulting model is a system of ordinary differential equations that is used to investigate the effect of high-frequency ventilation on strain in the parenchymal tissue. To apply homogenisation theory the simplifying assumption that lung tissue at the alveolar level is comprised of an array of units of similar size and shape in a highly idealized form is made. This allows the authors to move from a microscopic to a macroscopic space scale, from a single alveolus to an acinus. Since it is only feasible to solve the resulting fluid interaction problem on this lung geometry for a small number of alveoli, not for the thousands which make up a single acinus, the approach is to treat the structure as an array of repeating cells, representing alveoli, and to consider the average flow and deformation in a cell neglecting the microscopic details. The mathematical details are technical and beyond the scope of this thesis.

To further simplify the poroelastic equations we assume that the poroelastic continuum can be described by a solid phase (blood and tissue) and a fluid phase (air), where both phases are assumed to be incompressible. The interaction between the fluid pressure and the deformation of the solid skeleton is assumed to obey the effective stress principle. Note that by averaging over the tissue we do not seek to model individual alveoli but introduce macroscopic parameters such as the permeability and elasticity coefficients. In general, lung diseases usually affect significant regions of alveoli (lung tissue), thus, by changing the macro-scale parameters over the affected area of tissue we are still able to model changes in the tissue due to disease.

**Ignoring blood flow:** Apart from collagen, fibers and air the other major component in the lung is blood. The volume taken up by collagen and elastin fibers is similar to the volume occupied by the capillaries filled with blood (illustrated in Figure 7.1b). In fact, the space not occupied by air is about 7% of the parenchymal volume and is made up of 50% capillary blood and 50% of collagen and elastin fibers (Weichert, 2011). Also the density of blood is similar to the density of tissue and much larger than that of air ($1060\,\mathrm{kg\,m^{-3}} \gg 1.18\,\mathrm{kg\,m^{-3}}$). Since the capillaries are constantly filled with blood and the density of blood is similar to that of alveolar tissue we will make the simplifying assumption that the blood is simply part of the tissue (solid phase) and thus ignore accelerations and any redistribution of blood during breathing.

**Assuming incompressibility of the solid and the fluid:** Blood and tissue can be assumed to be incompressible. Under physiological conditions, air can also be assumed to be incompressible (Ismail et al., 2013).

**Ignoring solid inertia forces:** Simple calculations considering the sinusoidal motion of tissue near the diaphragm during normal breathing yield an estimate of $0.02\,\mathrm{ms}^{-2}$ for the maximum acceleration of lung parenchyma. Compared to the acceleration of gravity this is negligible, and it is therefore reasonable to ignore the inertia forces in the tissue.

**Ignoring fluid inertia forces:** The fluid's Reynolds number in the lower airways that form part of the lung parenchyma, has been estimated to be around 1 to 0.01 (Pedley et al., 1970b). Due to this relatively low Reynolds number we choose to ignore fluid inertia forces in the poroelastic medium.

**Ignoring viscous forces in the fluid:** A dimensional analysis shows that the viscous stress in the fluid is small compared to the drag forces between the fluid and the porous structure, when the ratio of the different length scales is small (Markert, 2007). We will therefore neglect the fluid viscous stress implying that the fluid behaves more or less inviscid within the porous structure.

### 7.4.2 Approximating the airways using a fluid network model

In order to make the coupled model computationally feasible we assume that a simple laminar flow model can describe the air flow in the airways and we make the common Poiseuille flow assumption. This flow assumption is also made in Leary et al. (2014); Swan et al. (2012) where the air flow in a whole airway tree, from trachea down to the final bronchioles was assumed to be governed by Poiseuille flow. Diseases affecting the airway tree can be modelled effectively by changing resistance (airway radius) parameters in the network flow model.

## 7.5 Mathematical model

### 7.5.1 A poroelastic model for lung parenchyma

Having made the assumptions in section 7.4 for the tissue we are left with the large deformation quasi-static incompressible poroelastic model (2.38).

**Constitutive laws.**

To close the poroelastic model for the tissue (2.38) we need to choose constitutive laws for the permeability and strain energy. We will use the same permeability law that has already been proposed in Kowalczyk and Kleiber (1994) to model lung parenchyma,

$$\boldsymbol{k}_0 = k_0 \left( J \frac{\phi}{\phi_0} \right)^{2/3} \boldsymbol{I}. \tag{7.1}$$

For a summary on previous defined variables see Table 2.1. Exponential strain energy laws for lung parenchyma exist, for example the popular law by Fung (1975). However little is known about how the constants in these laws should be interpreted and altered to model weakening of the tissue in an diseased state. Further, the constants in these laws are thought to have no physical meaning (Tawhai et al., 2009). To make the interpretation of the elasticity constants and dynamics of the model as simple as possible we chose a Neo-Hookean law taken from Wriggers (2008), with the penalty term chosen such that $0 \leq \phi < 1$,

$$W(\boldsymbol{C}) = \frac{\mu}{2}(\text{tr}(\boldsymbol{C}) - 3) + \frac{\lambda}{4}(J^2 - 1) - (\mu + \frac{\lambda}{2})\ln(J - 1 + \phi_0). \tag{7.2}$$

The material parameters $\mu$ and $\lambda$ can be related to the more familiar Young's modulus $E$ and the Poisson ratio $\nu$ by $\mu = \frac{E}{2(1+\nu)}$ and $\lambda = \frac{E\nu}{(1+\nu)(1-2\nu)}$. The values of these constants for modelling lung tissue have been investigated in De Wilde et al. (1981); Werner et al. (2009); Zhang et al. (2004) and are shown in Table 7.2.

### 7.5.2 A network flow model for the airway tree

The flow rate $Q_i$ through the $ith$ segment in the airway network is given by the pressure-flow relationship

$$P_{i,1} - P_{i,2} = R_i Q_i, \tag{7.3}$$

where

$$R_i = \frac{8l\mu_f}{\pi r^4}, \tag{7.4}$$

is the Poiseuille flow resistance of a pipe segment, where $r$ and $l$ are the radius and length of the pipe, $\mu_f$ is the dynamic viscosity, and $P_{i,1}$ and $P_{i,2}$ are the pressures at the proximal and distal nodes of the pipe segment, respectively. Let $\mathcal{A}_i$ be the set of pipe segments emanating from the $ith$ pipe segment in the airway network. We can express the conservation of flow in the airway network as

$$Q_i = \sum_{j \in \mathcal{A}_i} Q_j. \tag{7.5}$$

The outlet pressure of the airway network is set using the boundary condition $P_0 = \hat{P}$.

**Coupling the airway network to the poroelastic model.**

We introduce subdomains to identify the region of the domain that is supplied with fluid from a specific branch of the airway network and returns fluid through that branch. For notational purposes we use the subscript $di$ to indicate the most distal branches that have no further conducting branches emanating from them, but which enter a group of acinar units approximated by the continuous poroelastic model. We construct a Voronoi tesselation based on the $N$ terminal locations $\boldsymbol{y}_{di}, i = 1, \ldots, N$ of the airway network. The $i$th subdomain $\Omega_t^i$ is the subset of $\Omega_t$ that is closer to the $i$th terminal location at $\boldsymbol{y}_{di}$ than to any of the other terminal locations, i.e,

$$\Omega_t^i := \{\boldsymbol{x} \in \Omega_t : ||\boldsymbol{x} - \boldsymbol{y}_{di}|| < ||\boldsymbol{x} - \boldsymbol{y}_{dj}||,\ j = 1, 2..., N\ , j \neq i\}\,, \quad i = 1, \ldots, N. \tag{7.6}$$

Obviously we have $\Omega_t = \bigcup \Omega_t^i$. A simple 2D examples is shown in Figure 7.5.

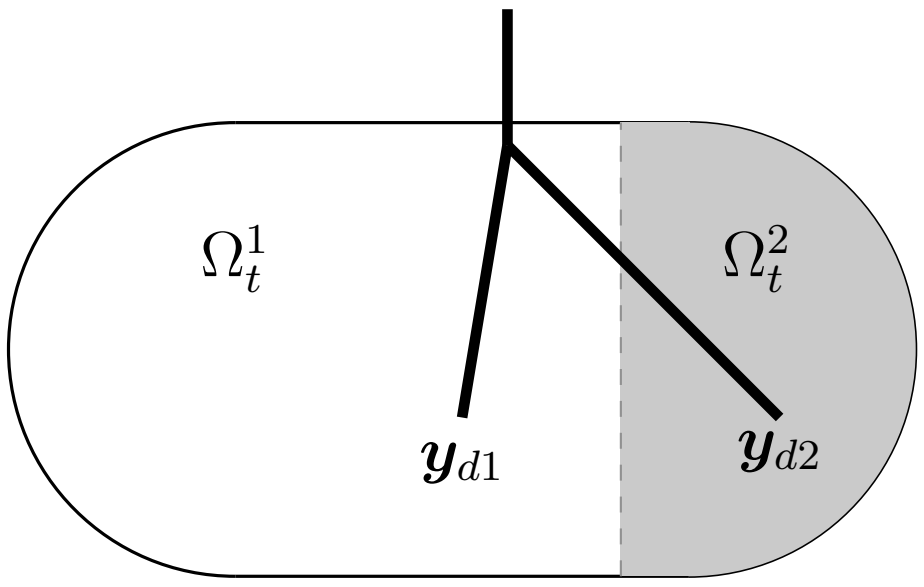

Figure 7.5: A simple example of a 2D domain being split into two subdomains according to (7.6).

We couple the airway network to the poroelastic domain in two ways. Firstly, the flux from each distal airway acts as a source term in the poroelastic mass conservation equation, namely

$$\nabla \cdot (\boldsymbol{\chi}_t + \boldsymbol{z}) = Q_{di} \quad \text{in } \Omega_t^i. \tag{7.7}$$

Secondly, the pressure at the distal airway $P_{di}$, determines the average pressure within subdomain $\Omega_t^i$, i.e.,

$$\frac{1}{|\Omega_t^i|}\int_{\Omega_t^i} p \,\mathrm{d}\Omega_t^i = P_{di}, \tag{7.8}$$

where $|\Omega_t^i|$ denotes the volume of the subdomian $\Omega_t^i$. Equation (7.8) enforces the assumption that the end pressure in a terminal bronchiole is the same as the alveolar pressure in the surrounding tissue.

### 7.5.3 The coupled lung parenchyma / airway model

To solve the coupled poroelastic-fluid-network lung model we need to find $\boldsymbol{\chi}(\boldsymbol{X}, t)$, $\boldsymbol{z}(\boldsymbol{x}, t)$, $p(\boldsymbol{x}, t)$, $P_i$ and $Q_i$ such that

$$\left.\begin{aligned}
-\nabla\cdot(\boldsymbol{\sigma}_e - p\boldsymbol{I}) &= \rho\boldsymbol{f} \quad \text{in } \Omega_t,\\
\boldsymbol{k}^{-1}\boldsymbol{z} + \nabla p &= \rho^f\boldsymbol{f} \quad \text{in } \Omega_t,\\
\nabla\cdot(\boldsymbol{\chi}_t + \boldsymbol{z}) &= Q_{di} \quad \text{in } \Omega_t^i,\\
\boldsymbol{\chi}(\boldsymbol{X}, t)|_{\boldsymbol{X}=\boldsymbol{\chi}^{-1}(\boldsymbol{x},t)} &= \boldsymbol{X} + \boldsymbol{u}_D \quad \text{on } \Gamma_D,\\
(\boldsymbol{\sigma}_e - p\boldsymbol{I})\boldsymbol{n} &= \boldsymbol{t}_N \quad \text{on } \Gamma_N,\\
\boldsymbol{z}\cdot\boldsymbol{n} &= q_D \quad \text{on } \Gamma_F,\\
p &= p_D \quad \text{on } \Gamma_P,\\
\boldsymbol{\chi}(\boldsymbol{X}, 0) &= \boldsymbol{X}, \quad \text{in } \Omega_0,\\
P_0 &= \hat{P},\\
P_{i,1} - P_{i,2} &= R_i Q_i,\\
Q_i &= \sum_{j\in\mathcal{A}_i} Q_j,\\
\frac{1}{|\Omega_t^i|}\int_{\Omega_t^i} p \,\mathrm{d}\Omega_t^i &= P_{di}.
\end{aligned}\right\} \tag{7.9}$$

## 7.6 Numerical solution of the coupled lung model

Since the system of equations (7.9) is highly nonlinear, its solution requires a scheme such as Newton's method. In Chapter 6 a finite element scheme using Newton's method for the solution of the poroelastic equations valid in large deformations (2.38) has already been presented. Here we adopt the same finite element scheme as presented in Chapter 6 for solving the poroelastic equations and expand the linear system (discretised linearisation) to include additional matrices required for solving the fluid network and its coupling to the poroelastic medium. This results in a monolithic coupling scheme that ensures good convergence even for problems with strong coupling interactions between the poroelastic medium and the fluid network. In section 7.6.1 we describe how to couple the fluid network to the discrete poroelastic model, and in section 7.6.2 we present details on how the stiffness matrix $\boldsymbol{K}$ (discretised linearisation of the full lung model (7.9)), and the residual vector $\boldsymbol{R}$ are built.

### 7.6.1 Discrete coupling of the fluid network to the poroelastic model

If we discretise the space using triangles and employ a piecewise constant pressure approximation (one node at the center of each element), the resulting coupling for the simple 2D example (Figure 7.5) is shown in Figure 7.6a. Once we refine the mesh (Figure 7.6b), the discretised division of subdomains tends to the subdivision of the original problem (Figure 7.5). The $i$th discretised subdomain $\Omega_t^i$ is defined as the set of elements $E_k$ whose centroids $\overline{\boldsymbol{x}_k}$ are closer to the distal end of the $i$th terminal branch than to the distal end of any other terminal branch, i.e.,

$$\Omega_t^i := \{E_k \in \Omega_t : ||\overline{\boldsymbol{x}_k} - \boldsymbol{y}_{di}|| < ||\overline{\boldsymbol{x}_k} - \boldsymbol{y}_{dj}||,\ j = 1, 2..., N\,, j \neq i\}\,. \qquad (7.10)$$

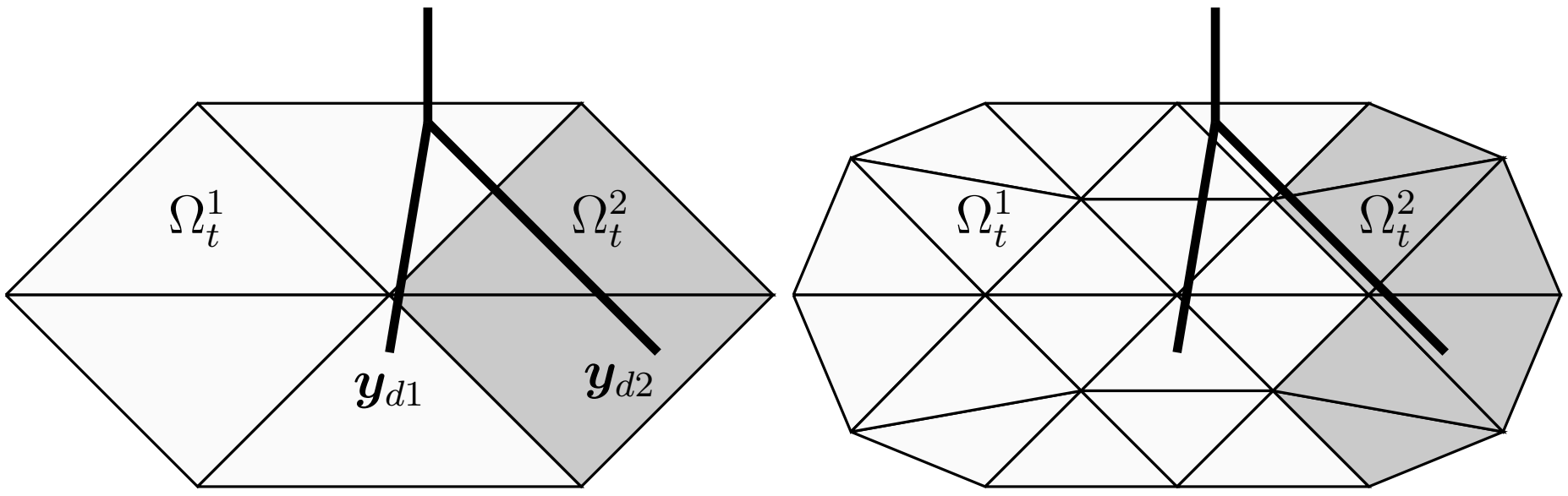

(a) (b)

Figure 7.6: (a) Coupling between the discretised domain and the fluid network using a piecewise constant pressure approximation for the example shown in Figure 7.5. (b) Coupling between the discretised domain and the fluid network after mesh refinement.

### 7.6.2 Finite element matrices

For the fully-coupled large deformation poroelastic fluid network model we need to solve the linear system $\boldsymbol{K}(\mathfrak{u}_i^n)\xi\mathfrak{u}_{i+1}^n = -\boldsymbol{R}(\mathfrak{u}_i^n, \mathfrak{u}^{n-1})$ at each Newton iteration. This can be expanded as

$$\begin{bmatrix} \boldsymbol{K}^e & 0 & \boldsymbol{B}^T & 0 & 0 & 0 & 0 & 0 \\ 0 & \boldsymbol{M} & \boldsymbol{B}^T & \boldsymbol{L}^T & 0 & 0 & 0 & 0 \\ -\boldsymbol{B} & -\Delta t\boldsymbol{B} & \boldsymbol{J} & 0 & 0 & 0 & 0 & -\Delta t\boldsymbol{G}^T \\ 0 & \boldsymbol{L} & 0 & 0 & 0 & 0 & 0 & 0 \\ 0 & 0 & 0 & 0 & \boldsymbol{T}_{11} & \cdots & \cdots & \boldsymbol{T}_{14} \\ 0 & 0 & 0 & 0 & \vdots & & & \vdots \\ 0 & 0 & 0 & 0 & \boldsymbol{T}_{31} & \cdots & \cdots & \boldsymbol{T}_{34} \\ 0 & 0 & \boldsymbol{G} & 0 & 0 & -\boldsymbol{X} & 0 & 0 \end{bmatrix} \begin{bmatrix} \xi\boldsymbol{\chi}_{i+1}^n \\ \xi\boldsymbol{z}_{i+1}^n \\ \xi\boldsymbol{p}_{i+1}^n \\ \xi\boldsymbol{\Lambda}_{i+1}^n \\ \xi\boldsymbol{P}_{i+1}^n \\ \xi\boldsymbol{P}_{d,i+1}^n \\ \xi\boldsymbol{Q}_{i+1}^n \\ \xi\boldsymbol{Q}_{d,i+1}^n \end{bmatrix} = - \begin{bmatrix} \boldsymbol{r}_{1,i} \\ \boldsymbol{r}_{2,i} \\ \boldsymbol{r}_{3,i} - \Delta t\boldsymbol{G}^T\boldsymbol{Q}_{d,i}^n \\ 0 \\ 0 \\ 0 \\ 0 \\ \boldsymbol{G}\boldsymbol{p}_i^n - \boldsymbol{X}\boldsymbol{P}_{d,i}^n \end{bmatrix}, \tag{7.11}$$

where

$$
\begin{aligned}
\boldsymbol{k}^e_{kl} &= \int_{(\Omega_{t_n})_i} \boldsymbol{E}_k^T \boldsymbol{D}(\boldsymbol{\chi}_i^n)\boldsymbol{E}_l + (\nabla\boldsymbol{\phi}_k)^T \boldsymbol{\sigma}_e(\boldsymbol{\chi}_i^n)\nabla\boldsymbol{\phi}_l \ \mathrm{d}(\Omega_{t_n})_i, \\
\boldsymbol{m}_{kl} &= \int_{(\Omega_{t_n})_i} \boldsymbol{k}^{-1}(\boldsymbol{\chi}_i^n)\boldsymbol{\phi}_k \cdot \boldsymbol{\phi}_l \ \mathrm{d}(\Omega_{t_n})_i, \\
\boldsymbol{b}_{kl} &= -\int_{(\Omega_{t_n})_i} \psi_k \nabla\cdot\boldsymbol{\phi}_l \ \mathrm{d}(\Omega_{t_n})_i, \\
\boldsymbol{j}_{kl} &= \delta \sum_{K\in\mathcal{T}_i^h} \int_{\partial K \setminus \partial(\Omega_{t_n})_i} h_{\partial K}[\psi_k][\psi_l] \ \mathrm{d}s. \\
\boldsymbol{r}_{1l} &= \int_{(\Omega_{t_n})_i} (\boldsymbol{\sigma}_e(\boldsymbol{\chi}_i^n) - p_i^n \boldsymbol{I}) : \nabla\boldsymbol{\phi}_l - \rho(\boldsymbol{\chi}_i^n)\boldsymbol{\phi}_l \cdot \boldsymbol{f} \ \mathrm{d}(\Omega_{t_n})_i \\
&\qquad\qquad - \int_{(\Gamma_N)_i} \boldsymbol{\phi}_l \cdot \boldsymbol{t}_N(\boldsymbol{\chi}_i^n) \ \mathrm{d}(\Gamma_N)_i, \\
\boldsymbol{r}_{2l} &= \int_{(\Omega_{t_n})_i} \boldsymbol{k}^{-1}(\boldsymbol{\chi}_i^n)\boldsymbol{\phi}_l \cdot \boldsymbol{z}_i^n - p_i^n \nabla\cdot\boldsymbol{\phi}_l - \rho^f(\boldsymbol{\chi}_i^n)\boldsymbol{\phi}_l \cdot \boldsymbol{f} \ \mathrm{d}(\Omega_{t_n})_i, \\
\boldsymbol{r}_{3l} &= \int_{(\Omega_{t_n})_i} \psi_l \nabla\cdot\left(\boldsymbol{\chi}_i^n - \boldsymbol{\chi}^{n-1}\right) + \Delta t \psi_l \nabla\cdot\boldsymbol{z}_i^n - \Delta t \psi_l g \ \mathrm{d}(\Omega_{t_n})_i \\
&\qquad\qquad + \delta \sum_{K\in\mathcal{T}^h} \int_{\partial K \setminus \partial(\Omega_{t_n})_i} h_{\partial K}[\psi_l][p_i^n - p^{n-1}] \ \mathrm{d}s, \\
\boldsymbol{l}_{kl} &= \int_{(\Omega_{t_n})_i} \epsilon_k \boldsymbol{\phi}_l \cdot \boldsymbol{n}(\boldsymbol{\chi}_i^n), \ \mathrm{d}(\Omega_{t_n})_i, \\
\boldsymbol{x}_{mn} &= \begin{cases} 1 & \text{if } ||\boldsymbol{y}_{dm} - \overline{\boldsymbol{x}_n}|| < ||\boldsymbol{y}_{dk} - \overline{\boldsymbol{x}_n}||, k = 1, 2..., N\ , k \neq m, \\ 0 & \text{otherwise}, \end{cases} \\
\boldsymbol{g}_{kl} &= \int_{(\Omega_{t_n})_i} \boldsymbol{x}_{kl} \frac{\boldsymbol{\phi}_l}{|E_l|} \ \mathrm{d}(\Omega_{t_n})_i,
\end{aligned}
$$

and $\boldsymbol{T}$ represents the matrix entries arising from equations (7.3) and (7.5). Here $\epsilon_k$ are scalar valued linear basis functions such that the Lagrangian multiplier vector at the $i$th iteration can be written as $\boldsymbol{\Lambda}_i^n = \sum_{k=1}^{n_\Lambda} \Lambda_{i,k}^n \epsilon_k$. Also, $\boldsymbol{P}^n$ and $\boldsymbol{Q}^n$ are the pressures at each junction and the fluid fluxes in each branch of airway network, except for the pressures at the distal end of, and the fluxes in, the most distal branches of the airway network which are given by $\boldsymbol{P}_d^n$ and $\boldsymbol{Q}_d^n$ respectively. Finally, $\overline{\boldsymbol{x}_n}$ denotes the centroid of the $n$th element. All other terms have already been defined in section 6.4.1.

## 7.7 Model generation

### 7.7.1 Mesh generation

We derive a whole organ lung model, of the right lung, from a high-resolution CT image taken at total lung capacity (TLC) and functional residual capacity (FRC). The bulk lung is first segmented from the CT data (slice thickness and pixel size 0.73 mm) using the commercially available segmentation software Mimics[1]. We then use the open-source image processing toolbox iso2mesh (Fang and Boas, 2009) to generate a Tetrahedral mesh containing 38369 elements. The conducting airways are also segmented from the CT data taken at TLC level, and a centerline with radial information is calculated. To approximate the remaining airways up to generation 8-13 we use a volume filling airway generation algorithm to generate a mesh of the airway tree containing 13696 nodes, with 2140 terminal branches (Bordas et al., 2016).

### 7.7.2 Reference state, initial conditions and boundary conditions

The poroelastic framework we have described requires a stress free reference state. In general, biological tissues do not possess a "reference state" where the material is free of both stress and strain, rather the cells that make up tissues are born into stressed states and live out their lives in these stressed states (Freed and Einstein, 2013). In order to define a stress-free reference state we scale the lung from FRC to a configuration in which the internal stresses and strains are assumed to be zero. The lung model is then uniformly inflated from the reference state to create a pre-stressed FRC configuration which has a mean elastic recoil of approximately $0.49 \times 10^3$ Pa, commonly understood to be a typical value (West, 2008). Thus the displacement of the boundary required to get from the reference state to FRC is given by

$$\boldsymbol{u}_{D,FRC} = (s-1)\boldsymbol{X}_{\partial\Omega}, \tag{7.12}$$

where $s$ is a scaling factor and $\boldsymbol{X}_{\partial\Omega}$ is the position of the lung surface in the reference state. From there we simulate tidal breathing. A similar approach has also been used in Lee et al. (1983).

[1] http://biomedical.materialise.com/mimics

We register the expiratory (FRC) segmentation to the segmentation at TLC using a very simple procedure that uses independent scalings $a_1, a_2$ and $a_3$ in the $x, y$ and $z$ direction, respectively, to map between the bounding boxes of the segmentations at FRC and TLC. This allows an estimate of the displacement for the lung surface from expiration to inspiration to be given by

$$\boldsymbol{u}_{D,TLC} = \begin{bmatrix} a_1 - 1 & 0 & 0 \\ 0 & a_2 - 1 & 0 \\ 0 & 0 & a_3 - 1 \end{bmatrix} (\boldsymbol{X}_{\partial\Omega} + \boldsymbol{u}_{D,FRC}) + \boldsymbol{b}, \tag{7.13}$$

where $\boldsymbol{b}$ is a translation vector to ensure that the top of the lung stays pinned throughout the simulation. To simulate tidal breathing we assume a sinusoidal breathing cycle and expand the lung surface from FRC to 40% of the displacement from FRC to TLC. Specifically,

$$\boldsymbol{u}_D(t) = \boldsymbol{u}_{D,FRC} + 0.2\left(1 + \sin(\frac{\pi}{2}(t+3))\right)\boldsymbol{u}_{D,TLC} \quad \text{on } \Gamma_D. \tag{7.14}$$

This results in a physiologically realistic tidal volume of 0.59 liters at a breathing frequency of 15 breaths per minute. We simulate breathing for a total of eight seconds or two breathing cycles. Due to the incompressibility of the poroelastic tissue, this also determines the total volume of air inspired/expired and the flowrate at the trachea, see Figure 7.8a and 7.8b respectively. We assume that no fluid escapes from the lung (except via the trachea) and impose zero flux boundary conditions at the lung surface. The outlet pressure of the airway network is set to zero (atmospheric pressure).

### 7.7.3 Simulation parameters

Several parameters for lung tissue elasticity and poroelasticity have been proposed (De Wilde et al., 1981; Lande and Mitzner, 2006; Lewis and Owen, 2001; Werner et al., 2009; Zhang et al., 2004). There is no consensus in the values in the literature. In this study we have chosen parameters from the literature, as shown in Table 7.2. These parameters are within range of existing models, and result in physiologically realistic simulation results (see section 7.8).

| Parameter | Value | Reference |
|---|---|---|
| $\phi_0$ | 0.99 | Lande and Mitzner (2006) |
| $\kappa_0$ | $10^{-5}$ m$^3$ s kg$^{-1}$ | Lande and Mitzner (2006) |
| $E$ | $0.73 \times 10^3$ Pa | De Wilde et al. (1981) |
| $\nu$ | 0.3 | De Wilde et al. (1981) |
| $\mu_f$ | $1.92 \times 10^{-5}$ kg m$^{-1}$ s$^{-1}$ | Swan et al. (2012) |
| $T$ | $8s$ | - |
| $\Delta t$ | $0.2s$ | - |
| $\delta$ | $10^{-5}$ | - |

Table 7.2: Parameters for breathing simulations.

## 7.8 Model exploration

We will now explore the behavior of the proposed model using a series of simulations to investigate the coupling between the airways and the tissue, dynamic hysteresis effects and how mass is conserved within the tissue.

In the subsequent analysis the total and elastic stress is calculated as $\sqrt{\lambda_1^2 + \lambda_2^2 + \lambda_3^2}$, where $\lambda_1, \lambda_2, \lambda_3$ are the three eigenvalues of the stress tensor, respectively. We define the relative Jacobian, denoted by $J_V$, as a measure for ventilation, which is calculated to be the volume ratio between the current state and FRC, i.e., $J_V = J/J_{FRC}$, and is a direct measure of tissue expansion. By running simulations over many breaths we have found that differences between the second breath and subsequent breaths were negligible, and therefore only results from the second breath, $t = 4s$ to $t = 8s$ are presented. The sagital slice shown in Figure 7.7a gives a good representation of the general dynamics within the tissue. Unless otherwise stated, all subsequent figures that do not show time courses are taken at $t = 5.8s$ just before peak inhalation of the second time breath in the simulation.

To solve the nonlinear poroelastic problem using Newton's method at a particular time step, we perform the the steps already described in Figure 6.1. We set the relative tolerance to be $\text{TOL} = 10^{-4}$. For the subsequent numerical results, a maximum of 5 Newton iterations were required to solve each time step.

### 7.8.1 Normal breathing

To simulate tidal breathing we apply the boundary conditions and simulation parameters previously discussed in sections 7.7.2 and 7.7.3, respectively.

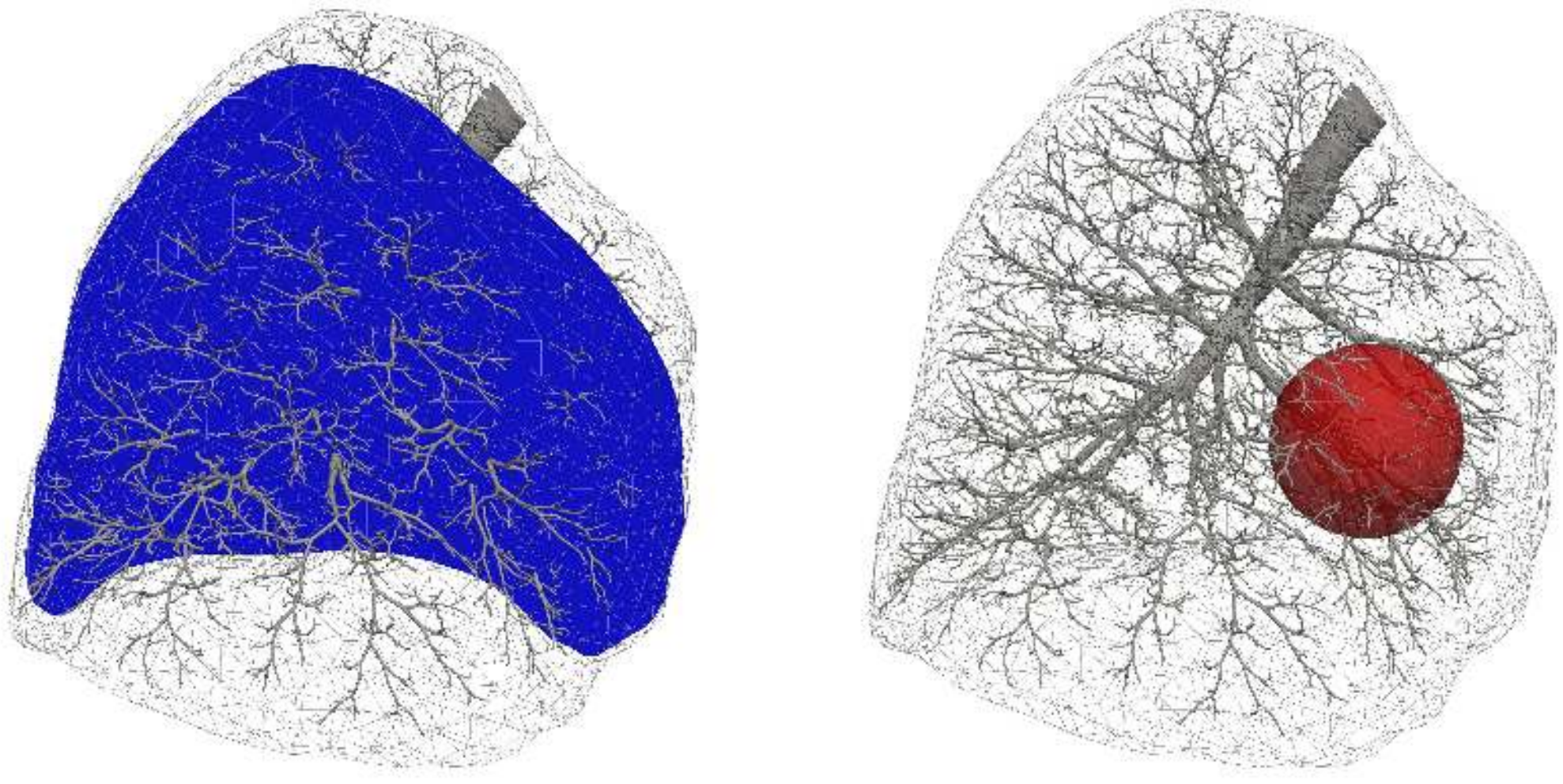

(a) (b)

Figure 7.7: (a) The blue sagital slice indicates the position of subsequent slices used for the data analysis of the tissue. (b) The red ball represents the structurally modified region, used to prescribe airway constriction and tissue weakening.

### Lung volume, flow and pressure drop

Figure 7.8 details the lung tidal volume, flow rate and pressure drop obtained from simulations of tidal breathing. Due to the incompressibility of the poroelastic medium and the fixed nature of the airway network, the lung tidal volume (Figure 7.8a) and flow rate (Figure 7.8b) follow a sinusoidal pattern that matches the from of the deformation boundary condition prescribed by equation (7.14). The mean pressure drop of the airways, is shown in Figure 7.8c, and agrees with previous simulation studies on full airway trees (Ismail et al., 2013; Swan et al., 2012).

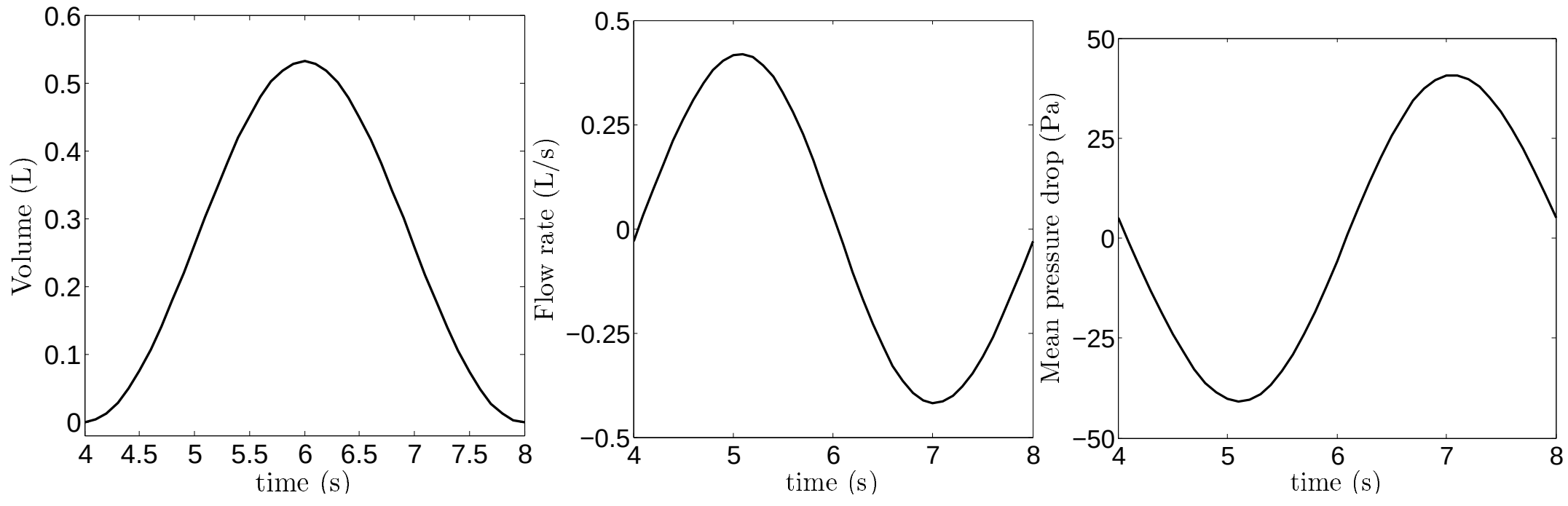

(a) (b) (c)

Figure 7.8: Simulated natural tidal breathing: (a) lung tidal volume (volume increase from FRC), (b) flow rate at the inlet, (c) mean pressure drop from the inlet to the most distal branches.

### Pathway resistance

The pathway resistance (Poiseuille flow resistance) from the inlet (right bronchus) to each terminal airway is shown in Figure 7.9a for the whole tree. In Figure 7.9b we show the pathway resistance of the terminal airways mapped onto the tissue.

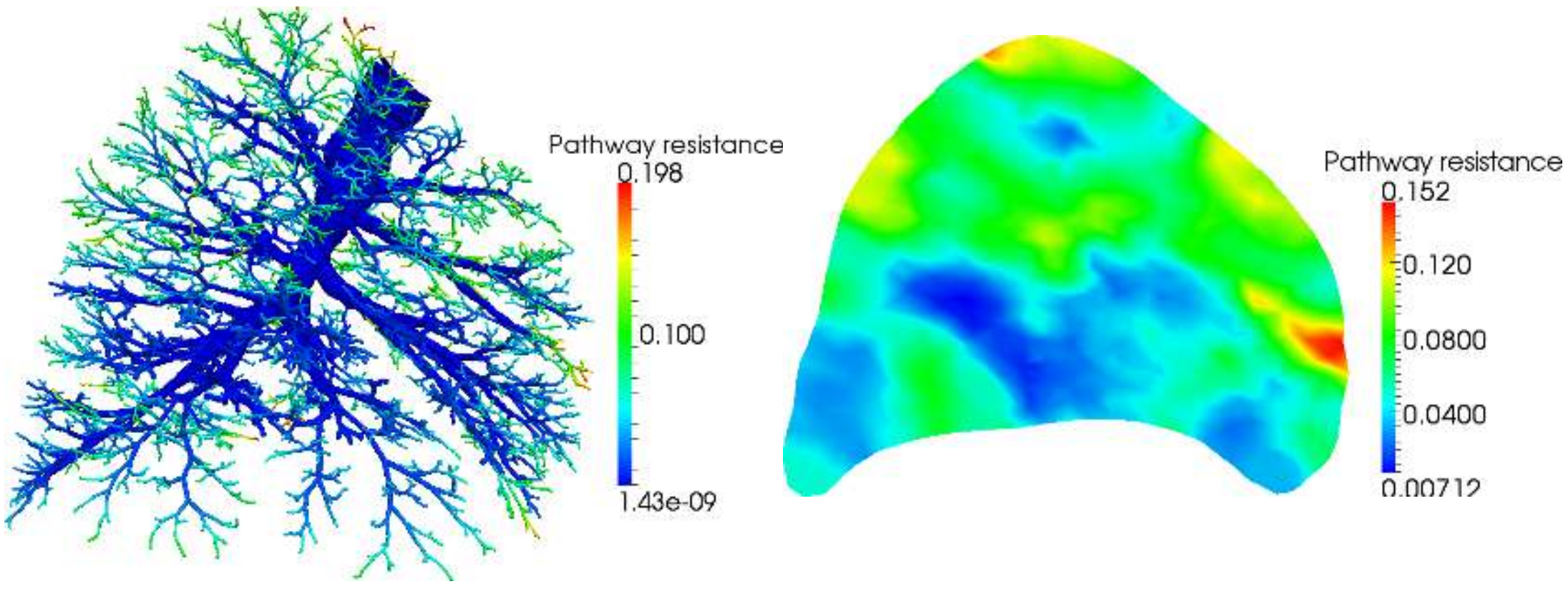

(a) (b)

Figure 7.9: (a) Pathway resistance ($\mathrm{Pa\,mm^{-3}s}$) from the inlet to the terminal branches in the airway tree. (b) Pathway resistance mapped onto a slice of tissue. The deformation of both the tree and the tissue in this figure correspond to the reference configuration.

### Airway tree-tissue coupling

In order to quantify the contribution of airway resistance to tissue expansion (ventilation), measured by $J_V$, the correlations between pathway resistance in the tissue and $J_V$ are plotted for each element in Figure 7.10a. There is a clear correlation between pathway resistance and tissue expansion, as is expected since the elastic coefficients are constant throughout the lung model. The Pearson correlation coefficients is $-0.55$, hence ventilation decreases as pathway resistance increases, with a p-value $< 0.0001$. Figure 7.10b shows there is also a strong correlation between the pathway resistance and pressure in the poroelastic tissue. Here the Pearson correlation coefficients is also $-0.55$, and pressure decreases (becomes more negative) with pathway resistance, with a p-value $< 0.0001$. Note that for a very few regions that are coupled to terminal branches with a low pathway resistance, positive pressures are possible. This results in a pressure gradient that pushes fluid from these well ventilated regions to neighbouring less ventilated regions (collateral ventilation). The distribution of pressure in the airway tree is shown in Figure 7.11a and the pressure inside the poroelastic tissue is shown in Figure 7.11b. Figure 7.11c shows the pressure on the lung surface. The patchy pressure field is well approximated by the piecewise constant pressure elements employed by the finite element method used to solve the poroelastic equations. Figure 7.11d shows the distribution of tissue expansion. Despite the heterogeneity in the airway tree the variations in tissue expansion are quite small, since the elastic coefficients are constant throughout the computational domain.

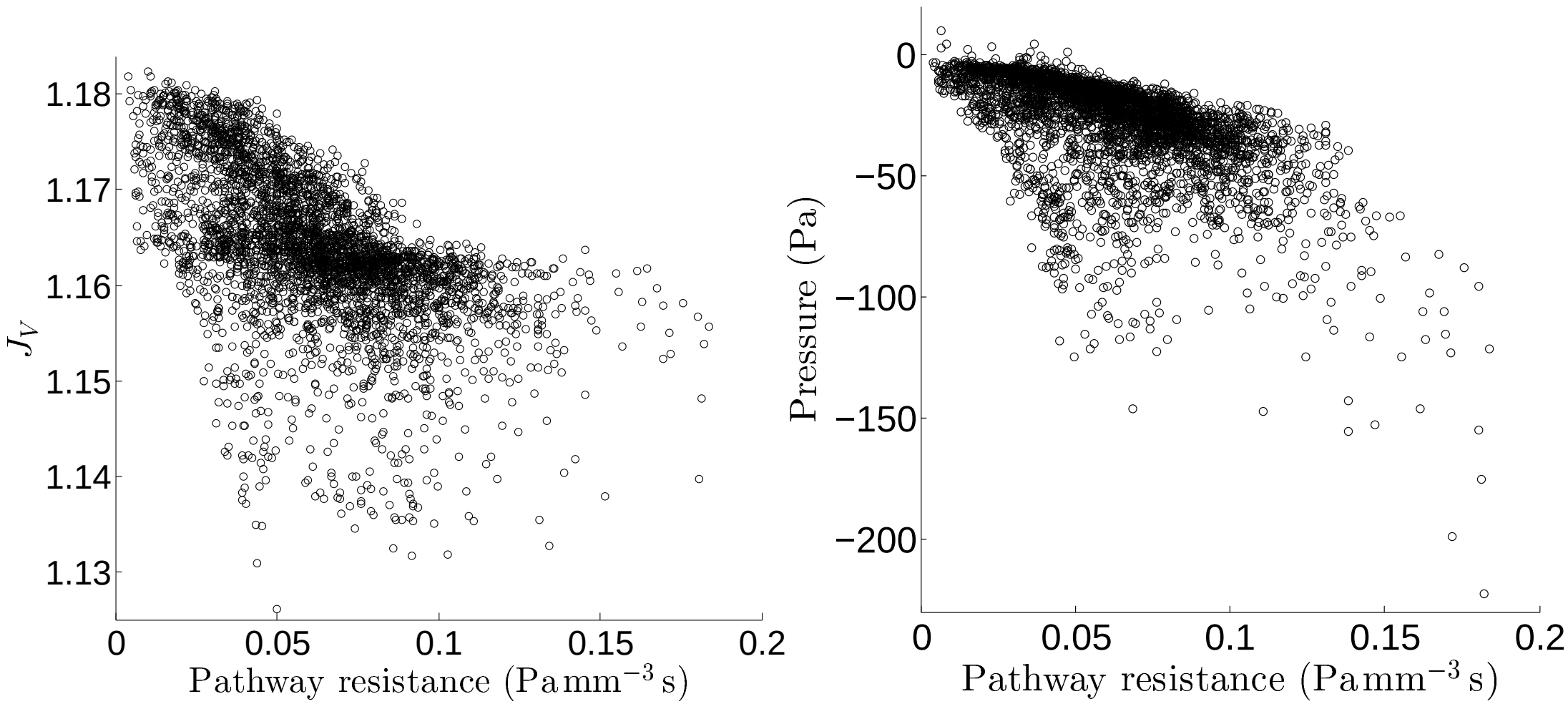

(a) (b)

Figure 7.10: (a) Correlation between tissue expansion (ventilation) and resistance of the pathways from the inlet to the terminal branch. (b) Correlation between pressure in the poroelastic medium (alveolar pressure) and pathway resistance.

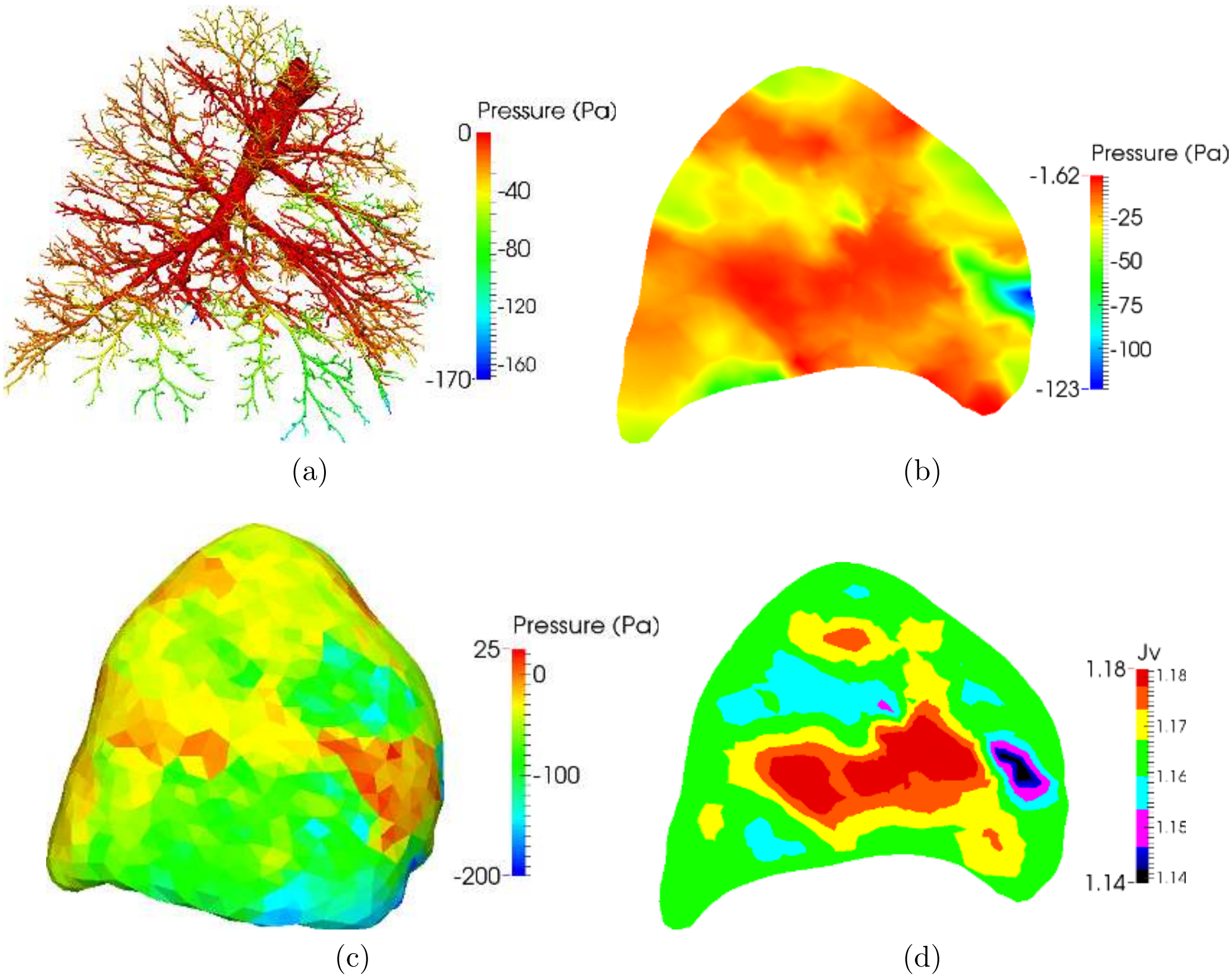

Figure 7.11: (a) Pressure in the airway tree. (b) Sagital slice showing pressure in the tissue using a linear interpolation. (c) Pressure on the lung surface. (d) Sagital slice showing tissue expansion from FRC.

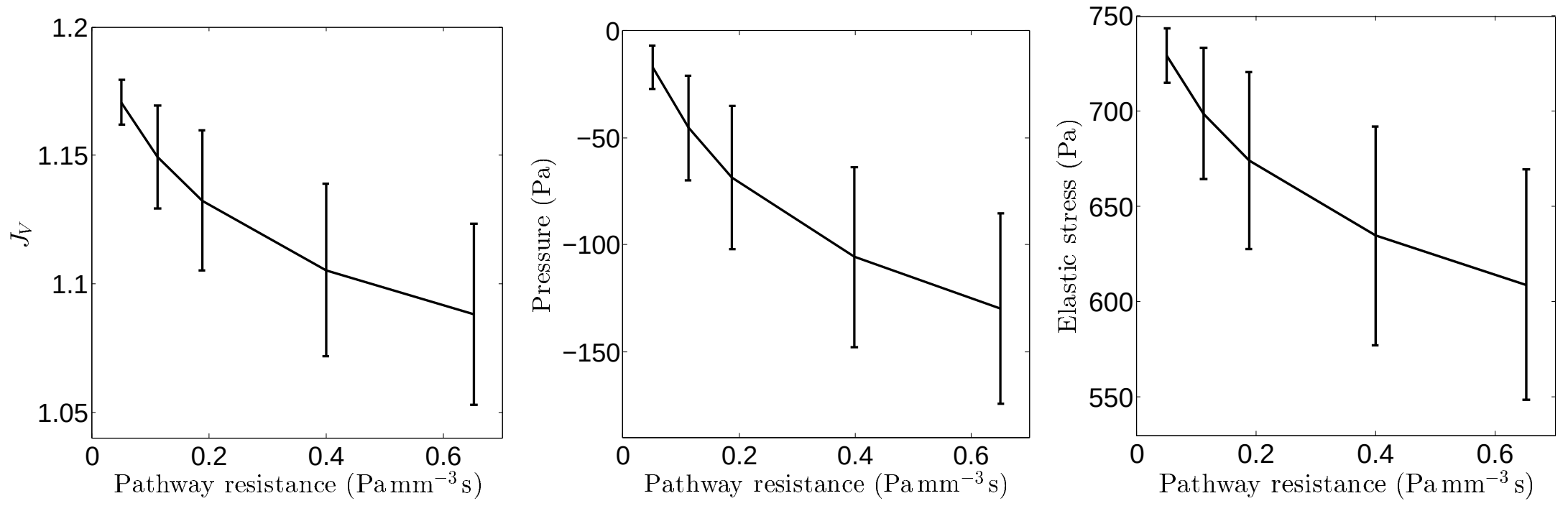

(a) (b) (c)

Figure 7.12: (a) Mean and standard deviations of the relative Jacobian from FRC, (b) pressure in the tissue and (c) elastic stress are plotted against increasing pathway resistance within the structurally modified region.

### 7.8.2 Breathing with airway constriction

We now simulate localised constriction of the airways by reducing the radii of the lower airways (with radius less than 4mm) within a ball near the right middle lobe. This region is represented by a red ball in Figure 7.7b. We reduce the radius of the aforementioned lower airways by $0\%, 40\%, 50\%, 60\%$ and $65\%$. This corresponds to a mean pathway resistance within the ball of $0.0507, 0.112, 0.188, 0.399$ and $0.651$ Pa mm$^{-3}$s, respectively. Figure 7.12 shows the changes in variables of physiological interest within the ball as the pathway resistance increases. The amount of tissue expansion during inspiration decreases as the airways become constricted (airway radius decreases and pathway resistance increases), as shown in Figure 7.12a. This is due to the reduced amount of flow in these airways. Further, the standard deviation increases because the pathway resistance of each branch increases by a different amount, depending on its original length and radius. Long and narrow branches will be affected most by the constriction. The pressure decreases with increasing pathway resistance as show in Figure 7.12b, since a larger pressure drop is needed to force the air down the constricted branches. Figure 7.12c shows the elastic stress in the tissue decreases as pathway resistance increases due to the decrease in tissue deformation (strain). However, as seen in Figure 7.13a, a large elastic stress appears near the boundary of the constricted region where the tissue is expanded by the surrounding tissue.

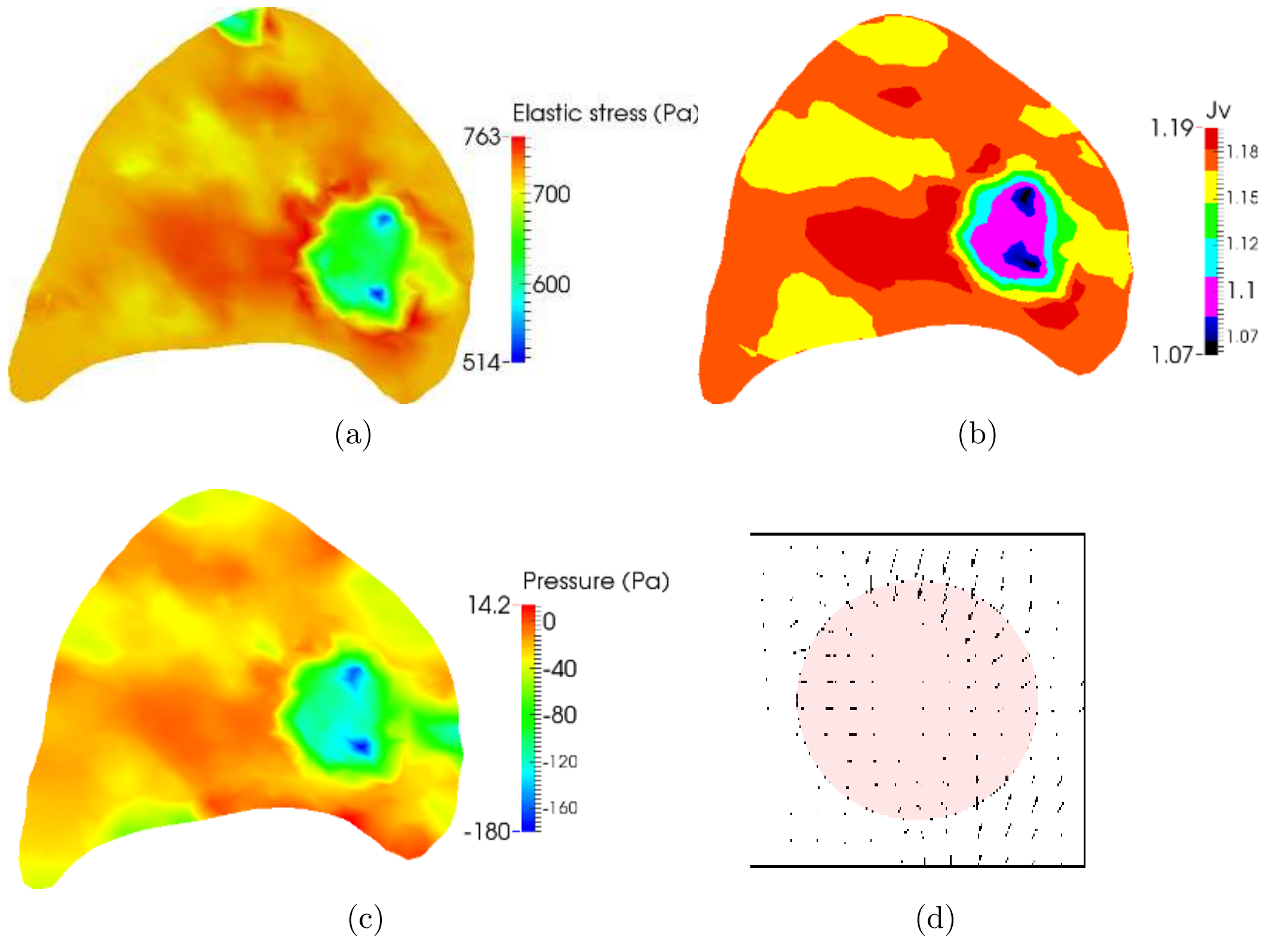

Figure 7.13: (a) Sagital slices showing the elastic stress, (b) Relative Jacobian, (c) pressure and (d) direction of the fluid flux near the structurally modified (constricted) region.

The simulation results shown in Figure 7.13 were performed with 65% airway constriction in the lower airways, applied within the structurally modified region. The volume conserving property (mass conservation) of the method is illustrated in Figure 7.13b where the tissue surrounding the constricted area is expanding to compensate for the reduction of tissue expansion due to the constriction within the structurally modified region. Figure 7.13c shows an increase in pressure near the boundary of this region. This facilitates a pressure gradient that allows for air to flow into the constricted region (collateral ventilation) to partially compensate for the reduced amount of ventilation, as is shown in Figure 7.13d. The magnitude of the maximum flow within the tissue is $8 \times 10^{-4}$ ms$^{-1}$, this is quite small and is due to the low permeability applied homogeneously within the model.

### 7.8.3 Breathing with locally weakened tissue

We now simulate localised weakening of the tissue by reducing the Young's modulus of the tissue within the structurally modified region represented by the red ball in Figure 7.7b. We reduce the Young's modulus by $0\%, 50\%, 75\%$ and $90\%$. This corresponds to a modified Young's modulus of $730, 365, 182.5$ and 73 Pa, respectively. Figures 7.14a-7.14c plot $J_V$, the pressure and the elastic stress within the modified region. As expected the local expansion increases as the tissue weakens, and the elastic stress decreases. Note that in all cases the range (heterogeneity) of local ventilation, pressure and elastic stress within the modified region increases dramatically as the stiffness of the modified region decreases.

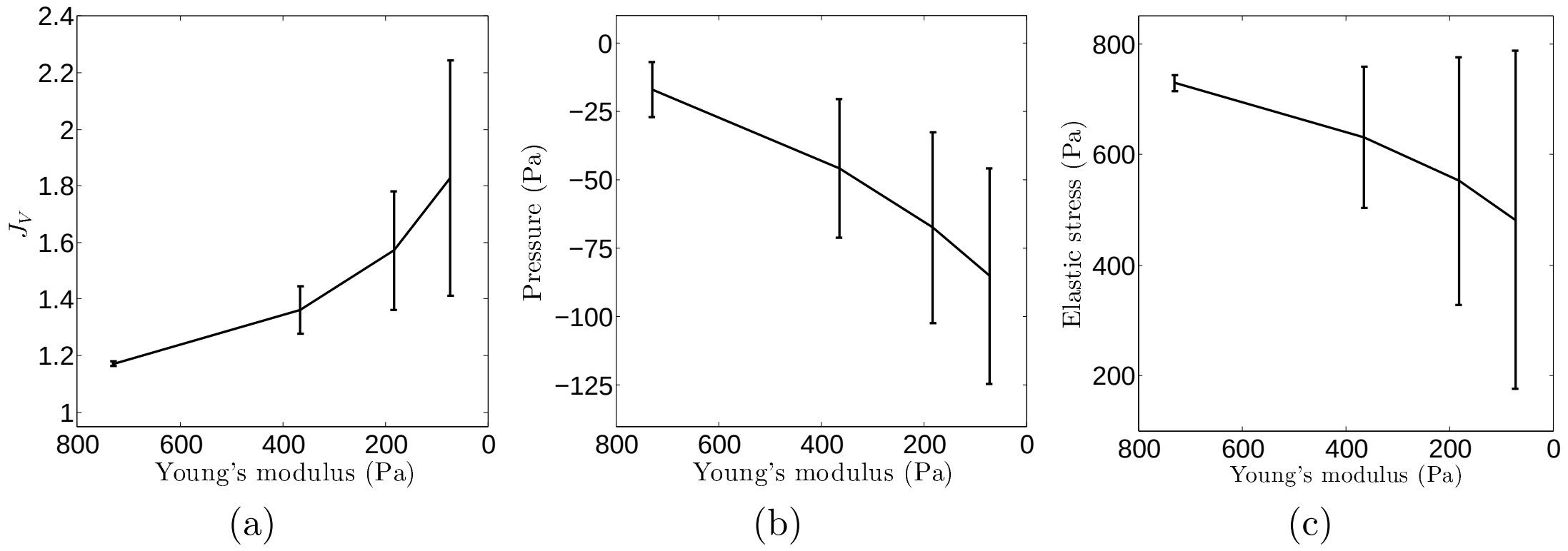

Figure 7.14: (a) Mean and standard deviations of the relative Jacobian from FRC, (b) pressure in the tissue and (c) elastic stress are plotted against Young's modulus within the structurally modified region.

Due to the large amount of tissue expansion within the structurally modified region, the tissue immediately surrounding this region is effectively squeezed between the expanded modified region and the surrounding tissue, and as a result expands the least, as seen in Figure 7.15.

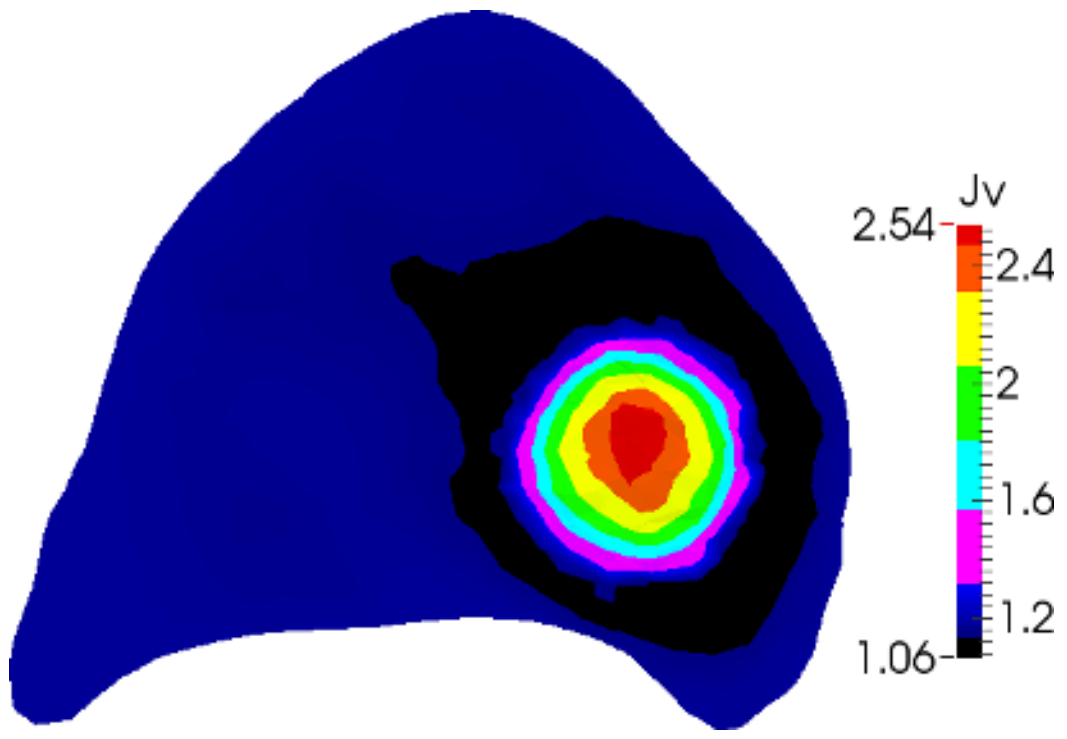

Figure 7.15: Slice showing the amount of tissue expansion ($J_V$) from FRC during inspiration with 90% localised tissue weakening.

### 7.8.4 Dynamic hysteresis

With the current choice of hyperelastic strain energy law (7.2) for the tissue mechanics, our model does not produce classic hysteresis effects, often attributed to surface tension within lung tissue (Kowalczyk and Kleiber, 1994). However, we are able to produce dynamic hysteresis effects, caused by delayed emptying and filling of parts of the lung.

Figure 7.16 shows the change in elastic recoil (total stress) with volume throughout the breathing cycle for three different breathing rates. This curve is commonly known as a dynamic pressure-volume (PV) curve, and shows the amount of dynamic hysteresis in the system. We will now explain the main features of this curve.

Figure 7.17a and 7.17b both show the distribution of pressure against pathway resistance within the tissue, shortly after inhalation. At this point the lung as a whole has started to exhale air. However some segments of the tissue have a negative pressure and are still filling up. These parts of the lung also tend to have a higher pathway resistance associated with them, which can explain the delayed filling. The reason that these parts of the lung continue to fill up, even during expiration, is that the continuum mechanics model of the tissue aims to achieve an energy minimum where the tissue is inflated evenly throughout

the lung, thus pulling open delayed segments of tissue. This is because the elasticity coefficients of the tissue have been parametrised homogeneously for these simulations. These negative pressures in the tissue, due to the delayed filling of parts of the lung, result in a larger total stress (elastic recoil), given by $\boldsymbol{\sigma} = \boldsymbol{\sigma}_e - p\boldsymbol{I}$. This effect is especially noticeable when transitioning from inspiration to expiration (and vice versa), causing the curve to shift right when moving from inspiration to expiration (due to delayed filling) and left when moving from expiration to inspiration (due to delayed emptying).

Also, we can clearly see an increase in the heterogeneity of the tissue's pressure distribution with increased breathing rate when comparing Figures 7.17a and 7.17b, for a four second and a one second breathing cycle, respectively. This increase in pressure heterogeneity is caused by the increased flow rates within the tree, and results in an increase in total stress. Therefore, a faster breathing rate causes an increasing amount of hysteresis (widening of the dynamic PV curve in Figure 7.16).

The increase of hysteresis in the dynamic PV curve and its shift as the breathing rate increases agrees with findings in the literature (Harris, 2005; Rittner and Döring, 2005). In the literature, hysteresis associated with dynamic PV curves is mostly hypothesized to be caused by flow-dependent resistances, pendelluft effects, chest wall rearrangement, and recruitment and derecruitment of lung units (Albaiceta et al., 2008; Harris, 2005; Ranieri et al., 1994).

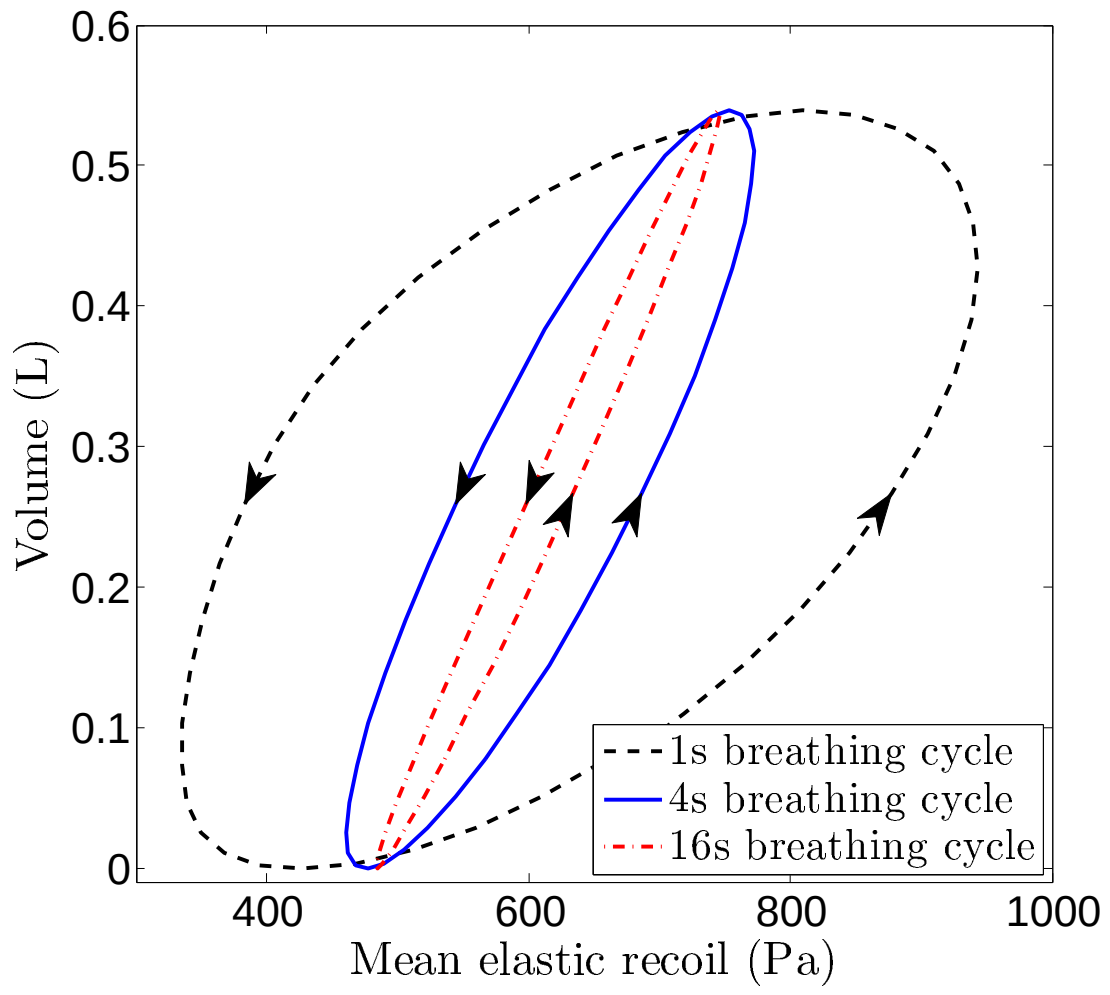

Figure 7.16: Dynamic pressure-volume curve: mean elastic recoil (total stress) against lung tidal volume during one full breathing cycle, for three different breathing rates. The arrows indicate the direction of time during the breathing cycle.

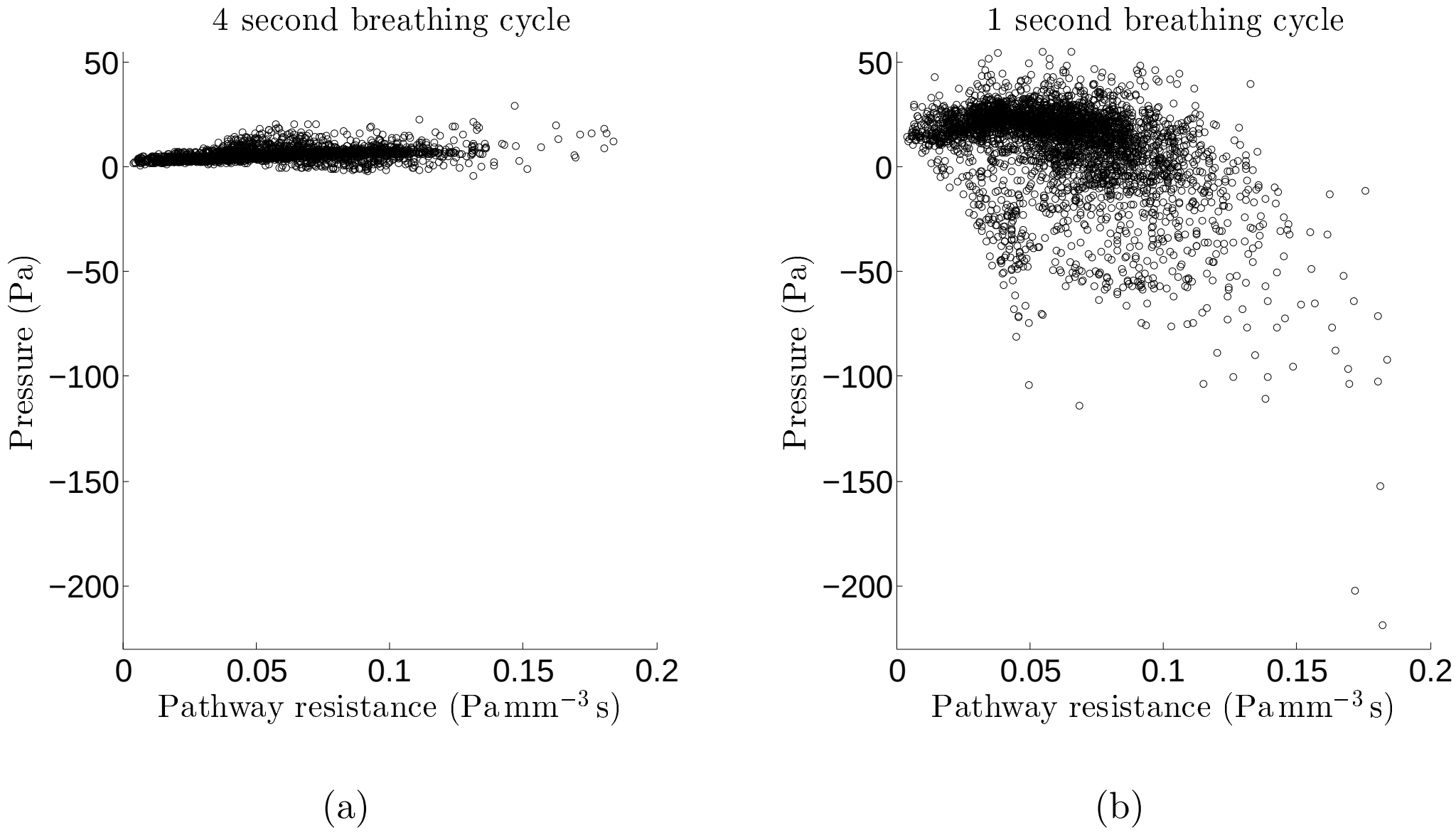

Figure 7.17: (a) Pathway resistance against pressure with a 4 second breathing cycle, 0.2 seconds after peak inhalation. (c) Pathway resistance against pressure with a 1 second breathing cycle, 0.05 seconds after peak inhalation.

## 7.9 Discussion

We have presented a mathematical model of the lung that tightly couples tissue deformation with ventilation using a poroelastic model coupled to a fluid network model. We have highlighted the assumptions necessary to arrive at such a model, and outlined its limitations. In comparison with previous ventilation models, the current approach models the tissue as a continuum and is therefore able to regionally conserve mass (which means conserve volume as the solid skeleton and fluid are both incompressible), and to model collateral ventilation. Further it is driven by deformation boundary conditions extracted from imaging data to avoid having to prescribe a pleural pressure which is impractical to be measured experimentally. In simulations of normal breathing, the model is able to produce physiologically realistic global measurements and dynamics. In simulations with altered airway resistance and tissue stiffness, the model illustrates the interdependence of the tissue and airway mechanics and thus the importance of a fully coupled model.

### 7.9.1 Contributors of airway resistance and tissue mechanics to lung function

We have found that there is a strong correlation between airway resistance and ventilation, see Figure 7.10a. Also, due to heterogeneity in airway resistance, hysteresis effects appear during breathing (Figure 7.16) and result in a complex ventilation distribution, caused by delayed filling and emptying of the tissue. Due to the Poiseuille law and the fourth power in airway radii that governs the resistance and flow through the airways (see equation (7.4)), small changes in airway radii can result in large changes in pathway resistance, which in turn can significantly affect the results of the coupled model. Thus, parameterising the airways correctly is very important. However this is notoriously difficult since CT data is only available down to the 5-6th generation, and small errors and biases in the segmentation, that get propagated by the airway generation algorithm, can have large influences in determining the simulation results. Changes in tissue elasticity coefficients also play an important role in determining the function of the lung model. This has been demonstrated in section 7.8.3 where are a reduction in the Young's modulus within a specified region causes significant changes in ventilation, pressure and stress.

The experiments performed in section 7.8.2 and section 7.8.3, begin to explore the sensitivity of changes in airway geometry and elastic parameters on the lung model's behaviour. However the changes in parameter were constrained to a small subregion of the model, making global inferences difficult. A more detailed sensitivity analysis should be performed to thoroughly investigate the importance of the airways and the tissue on lung function, as discussed in the future work section 8.2.2.

### 7.9.2 Limitations and future work

In order to move towards a more realistic model of the lung breathing, many steps need to be taken. We will list the main limitations that exist in the airway tree model, the poroelastic model, the boundary conditions and the geometry, and give indications on how these could be addressed in a future model.

**Airway tree limitations:** (1) The airway tree flow model currently implemented makes the Poiseuille flow assumption for the whole tree. The Poiseuille flow assumption requires flow to be fully developed and laminar. This may be true for the smaller airways where the Reynolds number is small but is certainly false for the larger upper airways where high Reynolds number flows occur. Such a model will therefore not be able to capture the high Reynolds number flows and turbulent effects that are known to exists in the upper airways. This could be improved by modifying the airway resistance at different generations according to the Reynolds number (Pedley et al., 1970a; Swan et al., 2012). Further improvements could be made by using a more sophisticated flow model for the airways, such as the 3D-0D model presented in Ismail et al. (2013). (2) The coupling of each terminal branch to the tissue currently assumes that there is no added resistance to air flowing from the terminal branch to each alveolar unit within the tissue. This could be improved by adding a simple resistive (impedance) model considering the volume of tissue that the terminal branch is feeding. This would also slightly increase the mean pressure drop of the lung model. (3) At the moment the airway tree is assumed to be static, and its configuration is not influenced by the deformation and stresses in the tissue. This could be improved by modelling the interaction of stresses and strains on the airway wall, opening up the airways during inspiration.

**Poroelastic tissue limitations:** (1) We have assumed a Neo-Hookean law for the strain-energy law to make the interpretation of the elasticity constants and dynamics of the model as simple as possible. However lung parenchyma is known to follow an exponential stress-strain relation, especially past tidal volume, where a law such as the one proposed by Fung (1975) might be more appropriate. Also little is known about the form of the strain-energy law during disease (e.g. fibrosis or emphysema). Similarly, for the permeability law little is known about its form for healthy or diseased tissue. Further experiments and modelling investigation would be needed to develop these. (2) Currently the tissue has been parameterised homogeneously to simplify the analysis of the results. Density information from CT images could be used to parameterise the initial porosity and elasticity coefficients. (3) We have ignored the effect of blood in the tissue. The inertia and gravity forces of blood acting on the tissue could be of importance when predicting deformation and ventilation in the lung. Due to the modular framework of the poroelastic theory it should be possible to include blood as a separate phase in a future version of the model. A vascular tree could also be generated from CT images and coupled to the poroelastic medium. (4) The airflow within the poroelastic tissue has been assumed to be inviscid. However, if we were to consider diseased states such as emphysema, where large areas of lung tissue completely break down leaving big holes, it could be argued that viscous forces could well play an important role, making it important to include them in our model. In a future version of the model the Darcy flow model could be replaced with a Brinkman, or even a Stokes flow model for big holes.

**Boundary condition limitations:** (1) The current registration should be updated to a more sophisticated nonlinear registration algorithm (e.g. Heinrich et al. (2013); Jahani et al. (2014); Yin et al. (2013)) that is able to account for the complicated deformation of the lung surface during breathing. (2) It is known that the lung surface is able to slide freely within the pleural cavity. This feature could be implemented using methods already presented in Kowalczyk and Kleiber (1994) and Ateshian et al. (2010).

**Geometry limitations:** (1) To model the complete organ and give a more accurate pressure drop, both the right and left lung, and the trachea and mouth should be included. (2) The airway tree generated in this work goes down to generations 8-13. More generations should be added to result in a fuller and more realistic tree. This would also require a finer mesh to approximate the lung tissue to resolve the coupling between each terminal branch and a subregion of lung tissue. (3) Cavities in the lung parenchyma due to large airways are currently not accounted for, i.e. it is assumed that the volume occupied by the airways is zero. To improve on this, a mesh of the lung with the larger upper airways removed would need to be generated. This new mesh could also incorporate a model of the cartilage found in the upper airways. (4) Additional no-flux boundaries should be introduced to represent the well defined and thought to be impermeable boundaries, between lobes (fissures) and lung segments.

**Validation:** No validation against experiments that contain spatial, mechanical or dynamic data has been made. Comparisons against information such as the vertical ventilation distribution from Single Photon Emission Computed Tomography (SPECT) data (Petersson et al., 2009) or pressure volume curve data obtained from experiments using the supersyringe method, the constant flow method, or ventilator method (Harris, 2005) should be made.

## 7.10 Conclusion

The model presented in this chapter can be used to investigate mechanical problems dependent on coupled deformation and ventilation in the lung. The numerical simulations are shown to be able to reproduce global physiologically realistic measurements. A fully nonlinear formulation permits the inclusion of various constitutive models, allowing investigation into different diseased states during various breathing conditions. A finite element method has been used to discretise the equations in a monolithic way to ensure convergence of the nonlinear problem, even under strong poroelastic-fluid-network coupling conditions. Due to the flexibility of the model, further improvements in its physiological accuracy are possible. It is hoped that the model presented here can form the basis for studies on the importance of airway and tissue heterogeneity on lung function, testing of mechanical hypotheses for the progression of disease, and investigations into phenomena such as hyperinflation, fibrosis and constriction.

The proposed lung model can also be used to validate and gain better insights into other types of computational lung models, such as zero dimensional compartment models that make extreme simplifying assumptions about the geometry of the lung (Bates, 2009). For example Whiteley et al. (2000) developed a multi-compartment ventilation model that is able to model inhomogeneous ventilation distributions in the lung. The 3D poroelastic lung model could be configured to mimic the simulations of this compartment model, and in an controlled simulation environment be used to confirm the effect of changes in parameters (e.g. resistances and compliances) on the resulting ventilation distributions within the lung.

# Chapter 8

# Conclusion

## 8.1 Review

In this thesis, we presented a low-order finite element method for solving the poroelastic equations valid in both small and large deformations. It has not been straightforward to arrive at the final formulation of the proposed stabilised finite element method. Only by performing detailed analysis of the error and stability of the discretised formulation were we able to determine the correct form of the stabilisation term that led to a stable and optimally converging method. This highlights the importance of rigorous analysis and testing when developing new numerical schemes. For the fully-discretised problem we proved existence and uniqueness, an energy estimate and an optimal a-priori error estimate. Numerical experiments performed in 2D and 3D illustrate the convergence of the method, and showed the effectiveness of the method to overcome spurious pressure oscillations. Due to the discontinuous pressure approximation, sharp pressure gradients due to changes in material coefficients or boundary layer solutions can be captured reliably, circumventing the need for severe mesh refinement. Thus, the proposed finite element method has made it possible to solve poroelastic models in biology previously not possible. As the numerical examples have demonstrated, the stabilisation scheme is robust and leads to high-quality solutions. A particularly nice feature is that in three dimensions only a very small value for $\delta$, the stabilisation parameter, is required to yield a stable solution, thus rendering the added mass effect of the stabilisation term negligible. This along with the method's simplicity compared to discontinuous and non-conforming finite element methods makes its implementation very appealing.

We also presented a mathematical (poroelastic) model of lung parenchyma that is coupled to a fluid network, modelling the airway tree. To the best of our knowledge, this is the first computational lung model built from patient specific imaging data that is able to capture the tight coupling between the tissue deformation and ventilation, as seen in Chronic Obstructive Pulmonary Diseases (COPD), such as emphysema. A numerical scheme to solve the coupled poroelastic fluid network has been presented and numerical software to simulate the lung model on patient specific lung geometries, extracted from imaging data has been implemented. Preliminary simulation results show physiologically realistic phenomena and have given some insights into the interdependence between ventilation and tissue deformation. The lung model appears to be a valid tool for solving the mechanical problem of tightly coupling lung deformation and ventilation during normal breathing and breathing with disease. We hope that due to the flexibility of the model, further improvements in its physiological accuracy, as outlined in section 7.9.2, will be made to yield an accurate whole organ lung model.

## 8.2 Future work

There are several areas which will pose interesting future research problems. These areas fall outside the scope of this work, but provide interesting challenges nonetheless.

### 8.2.1 Numerics

**Preconditioning:** By moving towards solving the poroelastic equations on more detailed 3D geometries the resulting linear system can grow to have several million degrees of freedom. For such problems direct solvers become impractical. To ensure robust and fast convergence of iterative methods such as the minimal residual method (MINRES), we need to precondition the linear system. An effective preconditioner for solving the Stokes problem using stabilised $P1 - P0$ elements has already been proposed in Wathen and Silvester (1993) and Silvester and Wathen (1994). This block preconditioning approach could be extended to the three-field poroelasticity case.

**A-posteriori error analysis:** A-posteriori error estimates could be derived for the finite element formulation of the linear porelasticity problem, which can be used for adaptive mesh refinement in space and time.

**Nonlinear elasticity:** There is a growing need for finite element methods of elasticity to capture steep pressure gradients due to material changes. For example changes in tissue types (fat, muscle and skin) when modelling the breast. To our knowledge there are currently no available finite element methods that use a simple to implement, low-order (discontinuous pressure) approximation to solve the incompressible nonlinear elasticity equations. It would be straightforward to extend the low-order method of nonlinear poroelastcicty to incompressible nonlinear elasticity.

### 8.2.2 Lung model

**Sensitivity analysis of airway geometry and elastic properties on lung function:** As shown in section 7.8, the proposed lung model has the capability to investigate the importance of airway resistance and tissue mechanics on lung function. A detailed sensitivity analysis should now be performed. The effect of changes in distribution of pathway resistance, upper and lower airway geometry, and distribution of elastic parameters within the tissue, on lung function should be investigated. A sliding boundary condition should be implemented, removing the need of having to prescribe deformation boundary conditions and subsequent flow rates. This would make it easier to relate the simulation results to global lung function, by being able to calculate physiologically meaningful measurements such as the force required by the diaphragm to achieve a given tidal volume.

**Constitutive laws for lung tissue:** Little is known about poroelastic constitutive laws for healthy and diseased lung tissue. Homogenisation theory (Lewis and Owen, 2001) and other modelling approaches such as spring models (Suki and Bates, 2011) could be used to derive new constitutive laws to better describe the elastic properties and fluid flow within the tissue.

**Validation:** For this model to be of practical use it is crucial that it is properly validated, this can be achieved by making use of different imaging modalities and phantom studies where model predictions can be tested. Computed tomography and 4D (dynamic) Magnetic resonance imaging (MRI) can be used to track displacements and calculate volume changes of lung structures. MRI of gases such as Hyperpolarised Xenon (Kaushik et al., 2011) and Helium 3 can be used to infer the flow and diffusion of gases, and with the use of elastography we are able to image stiffness and strain of lung tissue. Recently there has also been development in using Hyperpolarised Helium 3 MRI to estimate flow velocities and thus calculate pressure gradients (Patz et al., 2007).

**Surgical planning:** For patients with severe emphysema invasive surgical procedures such as lung volume reduction surgery (LVRS) and endobronchial valve placement are possible treatments. During LVRS part of the lung is excised in order to improve the configuration of the thoracic cavity, improve elastic recoil, and allow for improved lung inflation of the remaining and presumably better preserved tissue (Criner et al., 2011). Due to the high post-surgery mortality rate of around $5 - 10$ percent for LVRS and the fact that only some patients show an improvement with this therapy it is currently extremely challenging for doctors to select patients that will benefit from this invasive surgery. Boundary conditions allowing the lung surface to slide along the pleural cavity would have to be implemented, to allow for the removal of whole lobes in the model. A successful computational lung model would predict how much a particular patient will benefit from this high risk treatment, and help clinicians decide whether or not to perform surgery.

In addition to LVRS, various minimally invasive bronchoscopic approaches that also try to cure hyperinflation are being investigated. These include valves that reduce the air flow into the treated lobe during inspiration, stents that keep communications between pulmonary parenchyma and the segmental airways open, and lung volume reduction coils that aim to cause parenchymal compression and reduce the size of the hyperinflated tissue. More investigation into these techniques and which patients are best suited for a particular treatment is needed. A further developed computational lung model could be used to investigate these approaches and help surgeons plan for surgery by trialling different approaches in silico before the operation.

**Modelling other organs:** Finally, the proposed methodology for solving the lung model could also be adapted to model other biological tissues where blood vessels flow through and interact with a deforming tissue. For example, when modelling perfusion of blood flow in the beating myocardium (Chapelle et al., 2010; Cookson et al., 2012), modelling brain oedema (Li et al., 2010) or hydrocephalus (Wirth and Sobey, 2006), or microcirculation of blood and interstitial fluid in the liver lobule (Leungchavaphongse, 2013).

## 8.3 Final remarks

It is clear that there is a great requirement for effective simulation capabilities when it comes to modelling biological tissues. The possibility of robust and efficient simulations will enable researchers in the fields of medical device design, clinical treatment planning, and basic research. Although we have made some progress towards achieving this, still much research needs to be done, especially on how to implement models on high performance computers, to make detailed parameter studies possible.

The long term modelling aim of this project is to develop software which can accurately predict the ventilation and tissue deformation in the lungs. We have shown that, although such software would still be many years away from completion, requiring a great deal of work in the modelling, validation and biomechanical experimentation aspects, the aim is feasible and already computationally tractable.

# Appendix A

# Additional notation and workings

## A.1 Spatial tangent modulus

The spatial tangent modulus, fourth-order tensor, can be written as (see Bonet and Wood (1997, section 5.3.2) and Holzapfel et al. (2000, section 6.6))

$$\Theta_{ijkl} = \frac{1}{J} F_{iI} F_{jJ} F_{kK} F_{lL} \boldsymbol{C}_{IJKL}, \tag{A.1}$$

where $\boldsymbol{C}$ is the associated tangent modulus tensor in the reference configuration, given by

$$\boldsymbol{C}_{IJKL} = \frac{4\partial^2 W}{\partial C_{IJ} \partial C_{KL}} + pJ \frac{\partial C^{-1}_{IJ}}{\partial C_{KL}}. \tag{A.2}$$

## A.2 Matrix Voigt notation

To ease the implementation of the spatial tangent modulus we make use of matrix voigt notation. The matrix form of $\boldsymbol{\Theta}$ is given by $\boldsymbol{D}$, which can be written as (see Bonet and Wood (1997, section 7.4.2))

$$\boldsymbol{D} = \frac{1}{2}\begin{pmatrix} 2\boldsymbol{\Theta}_{1111} & 2\boldsymbol{\Theta}_{1122} & 2\boldsymbol{\Theta}_{1133} & \boldsymbol{\Theta}_{1112}+\boldsymbol{\Theta}_{1121} & \boldsymbol{\Theta}_{1113}+\boldsymbol{\Theta}_{1131} & \boldsymbol{\Theta}_{1123}+\boldsymbol{\Theta}_{1132} \\ & 2\boldsymbol{\Theta}_{2222} & 2\boldsymbol{\Theta}_{2233} & \boldsymbol{\Theta}_{2212}+\boldsymbol{\Theta}_{2221} & \boldsymbol{\Theta}_{2213}+\boldsymbol{\Theta}_{2231} & \boldsymbol{\Theta}_{2223}+\boldsymbol{\Theta}_{2232} \\ & & 2\boldsymbol{\Theta}_{3333} & \boldsymbol{\Theta}_{3312}+\boldsymbol{\Theta}_{3321} & \boldsymbol{\Theta}_{3313}+\boldsymbol{\Theta}_{3331} & \boldsymbol{\Theta}_{3323}+\boldsymbol{\Theta}_{3332} \\ & & & \boldsymbol{\Theta}_{1212}+\boldsymbol{\Theta}_{1221} & \boldsymbol{\Theta}_{1213}+\boldsymbol{\Theta}_{1231} & \boldsymbol{\Theta}_{1223}+\boldsymbol{\Theta}_{1232} \\ & \text{sym.} & & & \boldsymbol{\Theta}_{1313}+\boldsymbol{\Theta}_{1331} & \boldsymbol{\Theta}_{1323}+\boldsymbol{\Theta}_{1332} \\ & & & & & \boldsymbol{\Theta}_{2323}+\boldsymbol{\Theta}_{2332} \end{pmatrix}. \tag{A.3}$$

We also make use of the following implementation friendly matrix notation for $\nabla^S \boldsymbol{\phi}_k$,

$$\boldsymbol{E}_k = \begin{bmatrix} \phi_{k,1} & 0 & 0 \\ 0 & \phi_{k,2} & 0 \\ 0 & 0 & \phi_{k,3} \\ \phi_{k,2} & \phi_{k,1} & 0 \\ 0 & \phi_{k,3} & \phi_{k,2} \\ \phi_{k,3} & 0 & \phi_{k,1} \end{bmatrix}. \tag{A.4}$$

## A.3 Neo-Hookean strain energy

For the numerical examples we have used the following Neo-Hookean strain-energy law

$$W(\boldsymbol{C}) = \frac{\mu}{2}(\mathrm{tr}(\boldsymbol{C}) - 3) + \frac{\Lambda}{4}(J^2 - 1) - (\mu + \frac{\Lambda}{2})\ln(J - 1 + \phi_0). \tag{A.5}$$

Thus, the resulting effective stress tensor is given by

$$\boldsymbol{\sigma}_e = \frac{\Lambda}{2}\left(J - \frac{1}{J - 1 + \phi_0}\right)\boldsymbol{I} + \mu\left(\frac{\boldsymbol{C}^T}{J} - \frac{\boldsymbol{I}}{J - 1 + \phi_0}\right), \tag{A.6}$$

and the spatial tangent modulus tensor is given as

$$\boldsymbol{\Theta} = \boldsymbol{\Theta}_e + p(\boldsymbol{I} \otimes \boldsymbol{I} - 2\mathcal{Z}), \tag{A.7}$$

where

$$\begin{aligned}\boldsymbol{\Theta}_e = &\left[\Lambda J - 2\mu\left(\frac{1}{2(J - 1 + \phi_0)} - \frac{J}{2(J - 1 + \phi_0)^2}\right)\right]\boldsymbol{I} \otimes \boldsymbol{I} \\ &+ \left[\frac{2\mu}{J - 1 + \phi_0} - \Lambda(J - \frac{1}{J - 1 + \phi_0})\right]\mathcal{B},\end{aligned} \tag{A.8}$$

and

$$\mathcal{B}_{ijkl} = \frac{1}{2}(\delta_{ik}\delta_{jl} + \delta_{il}\delta_{jk}), \quad \mathcal{Z}_{ijkl} = \delta_{ik}\delta_{jl}, \quad \boldsymbol{I} \otimes \boldsymbol{I} = \delta_{ij}\delta_{kl}. \tag{A.9}$$

See Bonet and Wood (1997, chapter 5) and Wriggers (2008, chapter 3) for further details.

# Appendix B

# Computational considerations

## B.1 libMesh

All the numerical examples presented in this thesis were implemented using the C++ finite element library libMesh (Kirk et al., 2006). libMesh is an open-source library that has initially been developed at The University of Texas to provide a research platform for parallel adaptive finite element algorithms. The library has an active developer community, supports a range of standard and exotic elements in 2D and 3D, and has a good selection of example problems. Once the initial installation steps have been overcome the library is very accessible thanks to the detailed documentation.

## B.2 Linear solver

Another advantage of libMesh is that it interfaces with PETSc (Balay et al., 2015), the world's most widely used parallel numerical software library for partial differential equations.

### B.2.1 MUMPS

The Multifrontal Massively Parallel sparse direct Solver (MUMPS) is a direct method based on the LU factorization of sparse matrices (Amestoy et al., 2000), and available through the PETSc library. The solver handles both symmetric and nonsymmetric systems, allowing us to use this solver for all problems presented in this thesis. Because it is a direct method, no considerations about the convergence of the solver need to be taken, the solver will always produce the correct solution. Another advantage is the high parallelism of the method and its implementation. However parallelisation has not been investigated in this thesis.

### B.2.2 Alternatives for larger problems

The main disadvantage of direct solvers is that they require a lot of memory and can only fit 'small' problems into memory. Contrary to direct solvers, iterative methods approach the solution gradually, rather than in one large computational step. The big advantage for iterative solvers is that their memory usage is $O(N)$, allowing them to solve very large problems. The main disadvantage is that iterative solvers do not always converge. Different physics can require different solver settings and often need problem specific preconditioners to achieve convergence.

To solve the symmetric linear system of equations of linear poroelasticity (5.1), for large problems, a symmetric iterative solver such as the conjugate gradient method should be used. To solve the nonsymmetric large deformations problems of poroelasticity (6.12) and the lung (7.11), the generalized minimal residual method (GMRES) could be used. However further investigation into suitable preconditioners might well be required to obtain a solution. For a detailed explanation of these and other iterative solvers we refer to Elman et al. (2005).

## B.3 Nonlinear poroelasticity solver

The parameters for the Newton algorithm outlined in Figure 6.1 need to be chosen carefully. The maximum number of iterations, ITMAX, needs to be chosen large enough such that the Newton method can converge to the required tolerance at each time step. The relative tolerance, TOL, for the numerical experiments performed in section 6.5 has been chosen to be $10^{-4}$. An even lower tolerance could be chosen, however this would require more Newton iterations and will not necessarily result in a better approximation since the error due to the root finding is likely to be much smaller than the spatial and temporal finite element errors. Choosing a smaller time step can significantly reduce the number of required Newton steps, since each initial guess (the previous time step) is now much closer to the solution. This can result in an overall reduction in the computational time. Some experimentation in determining the optimal value for the tolerance and size of the time step is required, since these are heavily dependent on the problem under investigation.

Table B.2 shows the Newton convergence for the unconfined compression problem described in section 6.5.1 for the first timestep, which is the most demanding due to the initial displacement boundary condition. The resulting linear system contains 8162 degrees of freedom, takes $15.25s$ to assemble and $1.57s$ to solve, using one Intel Xenon CPU.

| Newton iteration | $\|\|\mathfrak{u}_i^n - \mathfrak{u}_{i-1}^n\|\|$ | $\|\|\boldsymbol{R}(\mathfrak{u}_i^n, \mathfrak{u}^{n-1})\|\|$ |
|---|---|---|
| 1 | 1.43852 | 0.00331502 |
| 2 | 0.553981 | 2.47657e-05 |
| 3 | 0.0149929 | 9.43506e-07 |
| 4 | 6.49539e-05 | 5.30122e-09 |

Table B.1: Convergence of the change in solution and residual for the unconfined compression test problem during the Newton iteration.

## B.4 Lung solver

Table B.2 shows the Newton convergence for the lung model as detailed in section 7.8 during the second timestep, which is a good representation of the other timesteps. The relative tolerance, TOL has been chosen to be $10^{-4}$. The resulting linear system contains 99009 degrees of freedom, takes $75.15s$ to assemble and $173.76s$ to solve, using one Intel Xenon CPU. Since the main nonlinearity is still the elasticity part, which has been linearized, the convergence is good and inline with the convergence of the nonlinear poroelasticity solver (Table B.2), provided a small enough time step is chosen.

| Newton iteration | $\|\|\mathfrak{u}_i^n - \mathfrak{u}_{i-1}^n\|\|$ | $\|\|\boldsymbol{R}(\mathfrak{u}_i^n, \mathfrak{u}^{n-1})\|\|$ |
|---|---|---|
| 1 | 0.42853 | 0.0273754 |
| 2 | 0.21266 | 0.0175368 |
| 3 | 0.000961185 | 0.000309249 |
| 4 | 9.90273e-05 | 4.02286e-05 |

Table B.2: Convergence of the change in solution and residual for the lung model during the Newton iteration.

### B.4.1 Data visulisation

All line plots presented in this thesis have been produced using MATLAB. The more complicated 2D and 3D visulisations have been produced using ParaView (Ahrens et al., 2005).